\pdfoutput=1
\RequirePackage{ifpdf}
\ifpdf 
\documentclass[pdftex]{sigma}
\else
\documentclass{sigma}
\fi

\numberwithin{equation}{section}

\newtheorem{Theorem}{Theorem}[section]
\newtheorem{Corollary}[Theorem]{Corollary}
\newtheorem{Lemma}[Theorem]{Lemma}
\newtheorem{Proposition}[Theorem]{Proposition}
 { \theoremstyle{definition}
\newtheorem{Definition}[Theorem]{Definition}
\newtheorem{Remark}[Theorem]{Remark} }

\input xy
\xyoption{all}

\def\C{{\mathbb C}}

\def\F{{\mathbb F}}
\def\H{\mathbb{H}}

\def\P{{\mathbb P}}

\def\R{\mathbb{R}}
\def\Z{\mathbb{Z}}
\def\N{\mathbb{N}}

\def\a{\alpha}
\def\bet{\beta}

\def\eps{\epsilon}

\def\Im{{\rm Im\,}}

\def\roots{{\Phi}}
\def\op{\hbox{ op }}

\def\wtil{\widetilde w}

\def\Wtil{\widetilde W}

\def\tW{{\widetilde{W}}}
\def\tw{{\widetilde{w}}}

\def\hV{\widehat V} 
\def\hfH{\widehat \fH} 
\def\hfk{\widehat \fk}
\def\hft{\widehat \ft}

\def\bm{\begin{bmatrix}}
\def\ebm{\end{bmatrix}}
\def\la{\langle}
\def\ra{\rangle}
\def\bi{{\bf i}}
\def\shf{\textrm{\scriptsize$\frac{1}{2}$\normalsize}}
\def\sbihf{\textrm{\scriptsize$\frac{\bi}{2}$\normalsize}}

\def\cA{\mathcal A}
\def\cB{\mathcal B}
\def\cC{\mathcal C}
\def\cF{\mathcal F}

\def\cI{\mathcal I}
\def\cL{\mathcal L}
\def\cP{\mathcal P}
\def\cS{\mathcal S}
\def\cT{\mathcal T}
\def\cCham{{\mathcal C}{\rm ham}}
\def\cU{\mathcal U}

\def\fb{\mathfrak{b}}
\def\g{\mathfrak{g}}
\def\b{\mathfrak{b}}
\def\fg{\mathfrak{g}}
\def\h{\mathfrak{h}}
\def\fh{\mathfrak{h}}
\def\fH{\mathfrak{H}}
\def\k{\mathfrak{k}}
\def\fk{\mathfrak{k}}
\def\fn{\mathfrak{n}}
\def\p{\mathfrak{p}}

\def\ft{\mathfrak{t}}
\def\fsl{\mathfrak{sl}}
\def\fsu{\mathfrak{su}}

\begin{document}
\allowdisplaybreaks

\newcommand{\arXivNumber}{1606.05638}

\renewcommand{\PaperNumber}{045}

\FirstPageHeading

\ShortArticleName{A Lightcone Embedding of the Twin Building of a Hyperbolic Kac--Moody Group}

\ArticleName{A Lightcone Embedding of the Twin Building \\ of a Hyperbolic Kac--Moody Group}

\Author{Lisa CARBONE~$^\dag$, Alex J.~FEINGOLD~$^\ddag$ and Walter FREYN~$^\S$}

\AuthorNameForHeading{L.~Carbone, A.J.~Feingold and W.~Freyn}

\Address{$^\dag$~Department of Mathematics, Rutgers University, Piscataway, New Jersey 08854, USA}
\EmailD{\href{mailto:lisa.carbone@rutgers.edu}{lisa.carbone@rutgers.edu}}
\URLaddressD{\url{https://sites.math.rutgers.edu/~carbonel/}}

\Address{$^\ddag$~Department of Mathematical Sciences, The State University of New York,\\
\hphantom{$^\ddag$}~Binghamton, New York 13902-6000, USA}
\EmailD{\href{mailto:alex@math.binghamton.edu}{alex@math.binghamton.edu}}
\URLaddressD{\url{http://people.math.binghamton.edu/alex/}}

\Address{$^\S$~Fachbereich Mathematik, Technical University of Darmstadt, Darmstadt, Germany}
\EmailD{\href{mailto:walter.freyn@googlemail.com}{walter.freyn@googlemail.com}}

\ArticleDates{Received July 23, 2019, in final form May 11, 2020; Published online May 29, 2020}

\Abstract{Let $A$ be a symmetrizable hyperbolic generalized Cartan matrix with Kac--Moody algebra $\mathfrak g=\mathfrak g(A)$ and (adjoint) Kac--Moody group $G = G(A)=\langle \exp({\rm ad}(t e_i)),\allowbreak \exp({\rm ad}(t f_i)) \,|\, t\in\C \rangle $ where $e_i$ and $f_i$ are the simple root vectors. Let $\big(B^+, B^-, N\big)$ be the twin $BN$-pair naturally associated to $G$ and let $\big(\mathcal B^+,\mathcal B^-\big)$ be the corresponding twin building with Weyl group $W$ and natural $G$-action, which respects the usual $W$-valued distance and codistance functions. This work connects the twin building $\big(\mathcal B^+,\mathcal B^-\big)$ of $G$ and the Kac--Moody algebra $\mathfrak g=\mathfrak g(A)$ in a new geometrical way. The Cartan--Chevalley involution,~$\omega$, of~$\mathfrak g$ has fixed point real subalgebra, $\mathfrak k$, the `compact' (unitary) real form of~$\mathfrak g$, and $\fk$ contains the compact Cartan $\mathfrak t = \mathfrak k \cap \mathfrak h$. We show that a real bilinear form $(\cdot,\cdot)$ is Lorentzian with signatures $(1, \infty)$ on $\mathfrak k$, and $(1, n -1)$ on~$\mathfrak t$. We define $\{k \in \fk \,|\, (k, k) \leq 0\}$ to be the {\it lightcone} of $\mathfrak k$, and similarly for $\mathfrak t$. Let $K$ be the compact (unitary) real form of $G$, that is, the fixed point subgroup of the lifting of $\omega$ to $G$. We construct a $K$-equivariant embedding of the twin building of~$G$ into the lightcone of the compact real form $\mathfrak k$ of $\mathfrak g$. Our embedding gives a geometric model of part of the twin building, where each half consists of infinitely many copies of a $W$-tessellated hyperbolic space glued together along hyperplanes of the faces. Locally, at each such face, we find an ${\rm SU}(2)$-orbit of chambers stabilized by~${\rm U}(1)$ which is thus parametrized by a Riemann sphere ${\rm SU}(2)/{\rm U}(1)\cong S^2$. For $n = 2$ the twin building is a~twin tree. In this case, we construct our embedding explicitly and we describe the action of the real root groups on the fundamental twin apartment. We also construct a spherical twin building at infinity, and construct an embedding of it into the set of rays on the boundary of the lightcone.}


\Classification{20G44; 20E42; 20F05; 51E24}

\section{Introduction}

Buildings were introduced by J.~Tits in the 1950's to provide a `geometric' interpretation of simple algebraic groups and finite groups of Lie type. Following ideas originally developed by F.~Klein in his Erlangen program, Tits aimed at understanding a~large class of groups, including simple algebraic groups and finite groups of Lie type, as the automorphism groups of carefully constructed geometric objects he called `buildings'. It turns out that the buildings of Tits correspond to groups which admit an additional structure, called a \emph{$BN$-pair} or equivalently a~\emph{Tits system}. This structure also gives rise to a~Bruhat decomposition for the corresponding group. Tits's approach was the geometric counterpart to the `algebraic' construction of these groups by C.~Chevalley using automorphisms of Lie algebras~\cite{Chevalley55}.

In general, the building $\cB$ of a group $G$ with a $BN$-pair is an abstract simplicial complex constructed from group theoretical data. The simplices in $\cB$ are in bijection with the union of all cosets~$G/P$, where $P$ runs through a set of representatives of the conjugacy classes of parabolic subgroups. These simplices satisfy incidence relations which can be phrased in terms of inclusions of these cosets.

For Kac--Moody groups, the closest infinite-dimensional analogue of simple algebraic groups, the structure of buildings and $BN$-pairs is richer. As Kac--Moody groups have two conjugacy classes of Borel subgroups, they admit the definition of two `opposite' $BN$-pairs which together form a `twin $BN$-pair'. Consequently the geometry associated to a Kac--Moody group $G$ naturally consists of two related components. This object is called a `twin building' $\cB=\cB^+\cup \cB^-$ associated to a twin $BN$-pair, \big($B^+, B^-, N$\big), where the subgroups $B^+$ and $B^-$ are the standard Borel subgroups constructed from the positive and negative roots of the Kac--Moody algebra respectively. Thus, by construction, the building is related to the combinatorial structure of its Kac--Moody group.

Kac--Moody algebras and groups fall naturally into three types, {\em finite type}, {\em affine type} and {\em indefinite type}. While Kac--Moody groups of finite type (simple Lie groups) and of affine type are well-understood, there are far fewer results known for the indefinite type. The most important subclass of indefinite type is the hyperbolic type, studied since the 1980s by various authors~\cite{Feingold80, FeingoldFrenkel83, KMW88, KangMelville95, LepowskyMoody79} but certainly of interest to physicists~\cite{DamourHenneauxNicolai02, DamourHillmann09, DamourKleinschmidtNicolai11, DamourSpindel13, Julia85, West01}.

\looseness=-1 While the algebraic properties of the hyperbolic Kac--Moody groups and algebras attracted some attention, there have been only a few mathematical results concerning their geometry, for example, the study of homogeneous or symmetric spaces associated to them~\cite{FreynHartnickHornKoehl17}.
These classes of objects are very well understood for finite-dimensional Lie groups and for affine Kac--Moody groups. It is hoped that the understanding of this geometry, called Kac--Moody geometry, will shed new light on algebraic and structural properties of these groups. As a~first step towards the goal of understanding the geometry of hyperbolic Kac--Moody algebras and groups, we show in this work that for
hyperbolic Kac--Moody groups~$G$ over~$\C$, the associated Tits building is not only an abstract simplicial complex admitting an action of~$G$, but that it admits a natural embedding inside the compact real form~$\k$ of the Kac--Moody algebra~$\fg$. As a consequence, the structure of the Tits building and of the compact real form of the Kac--Moody algebra are closely related.

Our results generalize work of Quillen and Mitchell in the finite-dimensional case, and Kramer and Freyn in the affine case. Mitchell, in a paper based on ideas of Quillen, used embeddings of spherical buildings associated to simple (real or complex) Lie groups into the tangent space of
the associated noncompact real symmetric space~\cite{Mitchell88}. In particular, the topological building of a real noncompact Lie group~$G$ with maximal compact subgroup~$K$ is canonically identified with a space homeomorphic to the unit sphere in the tangent space of the noncompact real symmetric space~$G/K$. For example, the topological building of type~$A_1$, isomorphic to~$S^1$, can be embedded into the unit circle in the tangent space of $\mathcal{H}^2={\rm SL}_2(\R)/{\rm SO}(2)$ which is $\R^2$.

Kramer gave a topological construction of the complex twin building of type $A_n^{(1)}$ and an equivariant embedding of this building into the associated affine Kac--Moody algebra $\g$~\cite{Kramer02}.
Freyn gave a 2-parameter family $\varphi_{\ell,r}$ of equivariant embeddings of affine `twin cities' \mbox{$\cB \!=\! \cB^+\!\cup\! \cB^-$} into the `$s$-representations' of affine Kac--Moody symmetric spaces~\cite{Freyn09, Freyn10b, Freyn10d, Freyn14, Freyn15b}. A `twin city' is the natural completion of a twin building, so affine twin cities correspond to completions of affine Kac--Moody groups as twin buildings correspond to minimal Kac--Moody groups. Twin cities carry a natural topology that is derived from the topology on the corresponding Kac--Moody group.

Denoting by $\g=\k \oplus \p$ the Cartan decomposition, and restricting Freyn's result to non-completed affine Kac--Moody algebras, yields an identification of the twin building with the intersection $\p_{\ell, r}$ of the sphere of squared length $\ell\in \R$ with horospheres parametrized by $r_d=\pm r\not=0$, where $r_d$ is the real coefficient of the derivation~$d$. The positive and negative components of the twin building, $\cB^+$ and $\cB^-$, are immersed into the two sheets of $\p_{\ell, r}$ described by $r_d<0$ respectively~$r_d>0$.

For affine Kac--Moody groups of the compact type, which are symmetric spaces called of \emph{`type~II'} following Helgason's classification (see~\cite{Freyn15, Helgason01}), this result includes embeddings of complex buildings into the compact forms of affine Kac--Moody algebras; this case is the affine counterpart to the embeddings constructed in this paper for symmetrizable hyperbolic Kac--Moody algebras. We refer the reader to~\cite{Freyn09, Freyn10b, Freyn10a, Freyn10d} for additional details.

For $\g$ a hyperbolic Kac--Moody algebra we are motivated in part by the appearance of $\fg/\fk$ and $G/K$ in coset models of certain supergravity theories, where coset spaces of the split real forms occur as parameter spaces for the scalar fields of the theory \cite{DamourHenneauxNicolai02,Julia85, West01}.

The ideas in this paper may be extended to the wider class of Lorentzian Kac--Moody algebras, but it remains to be seen how much the results are affected by the differences in the geometry of the lightcone and the Tits cone. Some work in this direction was started by A.~Tichai~\cite{Tich14}.

\subsection{Summary of results}\label{subsection:Summary of results}

The following paragraph summarizes the main results of this paper.
For $G$ of hyperbolic type, we construct a $K$-equivariant embedding of the twin building of $G$ into the lightcone of the compact
real form $\fk$. In $\ft$, $W$ acts on rays in the interior of the lightcone (and on certain rays on the nullcone), tessellating a copy of
hyperbolic space $H^{n-1}$ in each half of the lightcone (and some limit points in its boundary), giving a twin apartment. The $W$-images of the fundamental domain form the chambers of that apartment, so each face of any chamber is in a hyperplane fixed by a $W$-conjugate of a simple reflection~$w_i$. The associated conjugate of the compact subgroup ${\rm SU}(2)_i$ fixes that hyperplane pointwise, but rotates the rest of $\ft$ into a $K$-conjugate of~$\ft$, another apartment sharing that fixed chamber wall.
This gives our geometric model of part of the twin building where each half is infinitely many copies of a tessellated hyperbolic space glued together along hyperplanes of the faces. Locally, at each such face, a family of chambers meet, the orbit of an ${\rm SU}(2)$ with a ${\rm U}(1)$ stabilizer, so the family is indexed by a Riemann sphere, ${\rm SU}(2)/{\rm U}(1) \cong S^2$.

In rank $2$ the Weyl group is the infinite dihedral group, $D_\infty$, each apartment is a line, and each building $\cB^\pm$ is a tree. We model each apartment as a copy of the real line tessellated into unit intervals (chambers)
$C(n) = \big[n-\frac{1}{2},n+\frac{1}{2}\big]$ for $n\in\Z$, so the vertices are $\Z+\frac{1}{2}$.
At each vertex a family of intervals (chambers) is attached, each in a line (apartment) which is tessellated, and in each of those lines chambers are attached, on ad infinitum. We describe the action of the real root groups on the fundamental twin apartment. The family of chambers attached at any vertex is the projective space $P_1(\F)$, where $\F$ is the field over which the group is defined, so for $\F = \C$, we have $P_1(\C) = \hat\C$ is the Riemann sphere. Thus, in rank~$2$ each building~$\cB^\pm$ is a $\hat\C$-tree.
In that case, we also find a spherical twin building {\it at infinity}, and construct an embedding of it into the set of rays on the boundary of the lightcone.

In rank $3$ the Weyl group is a hyperbolic triangle group, each apartment is a copy of the Poincar\'e disk, $\cP$, tessellated into hyperbolic triangles by $W$. The boundary of each triangle is a segment in a hyperbolic geodesic. Along each geodesic segment we
have a $\hat\C$-family of attached triangles, each in a copy of~$\cP$, which is tessellated and has attached disks along each
geodesic, on ad infinitum. The geometrical embedding we found is only $K$-equivariant, so from a fundamental apartment we get only
to those apartments all of whose chambers are in some $K$-conjugate, $k\ft k^{-1}$. It means that in such apartments, pairs of triangles
attached to a geodesic edge are `balanced' opposite each other, corresponding to antipodal points in $\hat\C$. So under the action of~$K$, the fundamental $\cP$ can be rigidly rotated along any geodesic into another one, having only that geodesic line in common.
But the full Kac--Moody group~$G$ can leave a half-apartment fixed and rotate the other half up, creating a `hinge'. The complete apartment system would then be `hinged' copies of $\cP$ made up of pieces glued along geodesics.

\section{Kac--Moody algebras and Kac--Moody groups}

\subsection{Kac--Moody algebras}\label{subsection:Kac--Moody-algebra}

\noindent A Kac--Moody algebra $\mathfrak{g}_{\F}(A)$ over a field $\F$ may be constructed by generators and relations using a collection of data which includes a matrix $A=(a_{ij})_{i,j\in I}$ called a {\it generalized Cartan matrix} satisfying the following conditions for all $i,j\in I=\{1,\dots ,\ell\}$:
\[a_{ij}\in{\Z}, \qquad a_{ii}=2, \qquad a_{ij}\leq 0 \ \text{if} \ i\neq j, \qquad\text{and}\qquad a_{ij}=0 \iff a_{ji}=0. \]

A generalized Cartan matrix $A$ is {\it indecomposable} if there is no partition of the set $I = I_1\cup I_2$ into non-empty subsets
so that $a_{ij}=0$ for $i\in I_1$ and $j\in I_2$. The matrix $A$ is called {\it symmetrizable} if there is an invertible diagonal matrix
$D = \operatorname{diag}(d_1,\dots,d_\ell)$ such that $DA=(d_i a_{ij})$ is symmetric. One distinguishes various types of generalized Cartan matrices:
\begin{itemize}\itemsep=0pt
 \item[--] {\it Finite type:} $A$ is positive-definite. In this case $A$ is the Cartan matrix of a finite-dimensional semisimple Lie algebra and $\det(A)>0$.
 \item[--] {\it Affine type:} $A$ is positive-semidefinite, but not positive-definite, and all minors are positive definite. In this case $\det(A)=0$.
 \item[--] {\it Hyperbolic type:} $A$ is neither of finite nor affine type, but every proper, indecomposable submatrix is either of finite or of affine type. In this case $\det(A)<0$.
 \item[--] {\it Strictly hyperbolic type:} $A$ is hyperbolic type, but every proper, indecomposable submatrix is of finite type.
\end{itemize}

The terminology ``hyperbolic'' goes back to the original (independent) papers of Kac and Moody, and comes from the
geometry of the root systems for the corresponding Kac--Moody Lie algebras. The Weyl group orbits of roots lie on hyperbolas in the
rank~2 hyperbolic case and on hyperboloids in higher rank cases. We have included in Section 7 two figures illustrating the rank~2 root
systems which show this hyperbolic geometry clearly. The definitions of roots, root systems and Weyl groups are below.

A complex Kac--Moody algebra $\fg_\C(A)$ has at least two {\em real forms}, that is, real Lie algebras~$\fg_\R$ such that $\fg_\C(A) = \C\otimes \fg_\R$. The {\em split real form} of $\fg_\C(A)$ is $\fg_\R(A)$ (see~\cite{BBMR95}). {\it From this point on, all generalized Cartan matrices in this paper are assumed to be indecomposable, symmetrizable and of hyperbolic type.}

Given field $\F = \C$ or $\F = \R$:
\begin{itemize}\itemsep=0pt
 \item[--] a hyperbolic generalized Cartan matrix $A=(a_{ij})_{i,j\in I}$, and
 \item[--] a vector space $\h$ over $\F$ (which will play the role of a Cartan subalgebra) with $dim_\F(\h)=\ell$, and basis $\{h_i\, |\, i\in I\}$,
\end{itemize}
then there is a set of {\it simple roots} $\Pi=\{\alpha_j \, |\, j\in I\}\subseteq \h^*$ such that the pairing $\langle\cdot,\cdot\rangle\colon \h^*\times \h\to\F$ given by $\langle\alpha,h\rangle = \alpha(h)$ satisfies $\alpha_i(h_j) = a_{ij}$ for all $i,j\in I$, and
the hyperbolic Kac--Moody Lie algebra $\g=\g_\F(A)$ is generated by the elements $\{e_i, f_i, h_i\, |\, i\in I\}$, subject to the relations~\cite{GabberKac81, Kac90}:
\begin{itemize}\itemsep=0pt
 \item[--] $[h_i, h_j]=0$,
 \item[--] $[h,e_i]=\langle\alpha_i,h\rangle e_i$, $h\in \h$ and $[h,f_i]=-\langle\alpha_i,h\rangle f_i$, $h\in \h$,
 \item[--] $[e_i,f_j]=\delta_{ij}h_i$,
 \item[--] $({\rm ad}_{e_i})^{1-a_{ij}}(e_j)=0$, $i\neq j$ and $({\rm ad}_{f_i})^{1-a_{ij}}(f_j)=0$, $i\neq j$,
\end{itemize}

\noindent where ${\rm ad}_x(y)=[x,y]$. For each $i\in I$, let $\fsl_2^i$ be the Lie subalgebra (isomorphic to $\fsl_2(\F)$) with basis $\{e_i, f_i,h_i\}$, so that $\g$ is generated by these subalgebras. The abelian Lie subalgebra $\h$ with basis $\{h_i\, |\, i\in I\}$ is called the {\it standard Cartan subalgebra} of $\g$.

The algebra $\g=\g(A)$ is {\it infinite-dimensional} since $A$ is not positive definite, and it
admits an invariant symmetric bilinear form $(\hspace{2pt},\hspace{2pt})$ which is unique up to a global scaling factor
\cite[Section~II]{Kac90}, and which extends the form on $\h$ given by $2(h_i,h_j)/(h_j,h_j) = \alpha_i(h_j) = a_{ij}$. The nondegeneracy
of the pairing $\langle\cdot,\cdot\rangle$ between $\h^*$ and $\h$ determines a corresponding form on $\h^*$. This means that
$(\alpha_i,\alpha_j) = (h_i,h_j)$ and
\[a_{ji} (\alpha_i,\alpha_i)/2 = (\alpha_j,\alpha_i) = (\alpha_i,\alpha_j) = a_{ij} (\alpha_j,\alpha_j)/2, \]
so that we may take the diagonal matrix $D = \operatorname{diag}(d_1,\dots,d_\ell)$ with $d_i = 2/(\alpha_i,\alpha_i)$ and the symmetric matrix
$DA = (d_i a_{ij}) = (2 a_{ij}/(\alpha_i,\alpha_i))$. The standard way to choose
the global scaling factor is so that the longest square length of any simple root is~$2$.

The adjoint action of~$\h$ on~$\g$ is diagonalizable, and the simultaneous nonzero eigenspaces for that action,
\[\g_{\alpha} = \{x\in \g\,|\,[h,x]=\alpha(h)x,\, h\in \h\}\]
for $\alpha \neq 0$ are called {\it root spaces}.

The root system of $\g$ is the set $\roots= \{\alpha\in\h^*\, |\, \alpha \neq 0, \,\g_{\alpha}\neq 0\}$, and the $\Z$-span of $\roots$, called the
root lattice of $\g$, is denoted by $Q$. From the relations defining $\g$ we see that for each $i\in I$, $\g_{\alpha_i} = \C e_i$ and
$\g_{-\alpha_i} = \C f_i$, so that $\pm \alpha_i \in\roots$. In fact, we have
$Q=\Z\alpha_1\oplus\dots \oplus \Z\alpha_{\ell}$ is a free $\Z$-module.
Every $\alpha\in\roots$ can be written uniquely as $\alpha=\sum\limits_{i=1}^{\ell} k_i\alpha_i$ where either all $k_i\geq 0$, in which case $\alpha$ is called {\it positive}, or all $k_i\leq 0$, in which case $\alpha$ is called {\it negative}. The set of all positive roots is denoted $\roots^+$, and
the set of all negative roots is denoted by $\roots^-$. Any root is either positive or negative.

For each simple root $\alpha_i$, $i\in I=\{1,\dots, \ell\}$, we define the {\it simple root reflection}
\begin{gather}\label{equation:simple_reflection}
w_i(\alpha_j) := \alpha_j - \alpha_j(h_i)\alpha_i.
\end{gather}
The set $S = \{w_i\,|\, 1\leq i\leq\ell\}$ generates a group $W=W(A)$ of orthogonal transformations of $\h^*$, called the {\it Weyl group} of $A$. The non-degenerate pairing between $\h$ and $\h^*$ gives a corresponding action of $W$ as orthogonal transformations on $\h$. A root $\alpha\in\roots$ is called a {\it real root} if there exists $w\in W$ such that $w\alpha$ is a~simple root. A root $\alpha$ which is not real is called {\it imaginary}. We denote by $\roots^{\rm re}$ the set of all real roots and $\roots^{\rm im}$ the set of all imaginary roots.

The multibracket, $[e_{i_1},e_{i_2},\dots,e_{i_n}] = {\rm ad}_{e_{i_1}} {\rm ad}_{e_{i_2}} \cdots {\rm ad}_{e_{i_{n-1}}} e_{i_n}$ is in the root space~$\fg_\a$ for $\alpha = \sum\limits_{j=1}^n \alpha_{i_j} \in \roots^+$, while a similar multibracket with each $e_{i_j}$ replaced
by $f_{i_j}$ is in $\fg_{-\a}$. Therefore, $\g$~has a root space decomposition
\cite[Theorem~1.2]{Kac90}
\[\g = \h \oplus \bigoplus_{\alpha\in\roots} \g_{\a} = \h \oplus \g^+ \oplus \g^-,\]
where
\[\g^+ =\bigoplus_{\alpha\in\roots^{^+}}\g_{\alpha},\qquad
\g^- = \bigoplus_{\alpha\in\roots^{^-}}\g_{\alpha}.\]

The standard positive Borel subalgebra $\b\equiv\b^+$ is defined by $\b^+=\h\oplus \g^+$ and the standard negative Borel subalgebra by $\b^-=\h\oplus \g^-$.

\subsection{Kac--Moody groups}

There are various ways to define abstract Kac--Moody groups
(see for example~\cite{KacPeterson85c, Marquis18, Remy02,Tits87}).
The main point of the abstract approach is to give a flexible definition of Kac--Moody groups, allowing the construction of groups whose adjoint action is not faithful. There are indeed important examples of that kind, the smallest one being the finite type Kac--Moody group ${\rm SL}(2, \C)$, where the two matrices $\pm Id$ both act as the identity operator in the adjoint representation. Hence the adjoint representation of this Kac--Moody group is actually the group ${\rm PSL}(2, \C)$, and that is what we get by using the definition of the adjoint Kac--Moody group given in equation~(\ref{equation:adjoint_Kac_Moody_group}) for the Cartan matrix $A = [2]$ of type $A_1$.

By the definition of the abstract Kac--Moody group, there is a surjective group homomorphism: $Ad: G\longrightarrow G^{\rm ad}$ from an abstract Kac--Moody group onto the adjoint Kac--Moody group whose kernel is exactly the center of G (see~\cite[Proposition~9.6.2]{Remy02}). As we will see in Section~\ref{subsection:BNpair_Tits_building}, for the subgroups $B^\pm$ and $N$ defined there, we have
$B^\pm\cap N$ is abelian, and the center of $G$ is the kernel of the action of $G$ on the twin building (see~\cite[Lemma~1.7]{Caprace09}). Hence, without loss of generality, to understand the action on twin buildings we can work with the adjoint Kac--Moody group.
Our references for this section are~\cite{ KacPeterson85c,Kumar02, Marquis18, Moody95}.

Let $\g$ be a symmetrizable Kac--Moody algebra over $\C$, $L$ be a complex vector space and let
$\phi\colon \fg\to {\rm End}(L)$ be any integrable representation, so that
all $\phi(e_i)$ and $\phi(f_i)$ are locally nilpotent on $L$ and the linear operators
\[\chi^\phi_{{\a}_i}(t) = \exp(\phi(te_i))\qquad\hbox{and}\qquad \chi^\phi_{{-\a}_i}(t) = \exp(\phi(tf_i)),\quad \text{for} \ t\in\C,\]
are well-defined in ${\rm GL}(L)$. In fact, for any $x\in \g_{\a}$, $\a\in\roots^{\rm re}$ the operator $\phi(x)$ is locally nilpotent on~$L$, so
$\chi^\phi_x = \exp(\phi(x))$ is well-defined and these give the real root groups $U_\a^\phi$. In particular, the adjoint representation
${\rm ad}\colon \fg \to {\rm End}(\g)$, is integrable, and for all $x\in \g_{\a}$, $\a\in\roots^{\rm re}$ we have well-defined operators $\chi^{\rm ad}_x = \exp({\rm ad}_x)\in {\rm GL}(\g)$ giving the real root groups $U_\a^{\rm ad}$.

\begin{Definition}[minimal Kac--Moody groups] \label{definition:Kac--Moody_groups} Let the {\it minimal Kac--Moody group} associated to an integrable representation
$\phi$ be the group generated by these operators,
\[G^\phi = G^\phi(\C)=\big\la \chi^\phi_{{\a}_i}(t), \chi^\phi_{{-\a}_i}(t)\,|\, i\in I,\ t\in\C\big\ra\leq {\rm GL}(L).\]
In particular, this defines the {\it minimal adjoint Kac--Moody group}
\begin{gather}\label{equation:adjoint_Kac_Moody_group}
G^{\rm ad} = G^{\rm ad}(\C)=\big\la \chi^{\rm ad}_{{\a}_i}(t), \chi^{\rm ad}_{{-\a}_i}(t)\,|\, i\in I,\, t\in\C\big\ra\leq {\rm GL}(\g).
\end{gather}
These generators act on $\g$ as Lie algebra automorphisms, so we have $G^{\rm ad}\leq {\rm Aut}(\g)$.
\end{Definition}

Since ${\rm ad}_{e_i}$ and ${\rm ad}_{f_i}$ are locally nilpotent, $\g$ is the direct sum of finite-dimensional $\fsl_2(\C)$-modules for each of the subalgebras $\fsl_2^i$. For fixed $i\in I$, on each such summand the exponentials above generate the group ${\rm SL}_2^i$ isomorphic to ${\rm SL}_2(\C)$, so $G$ is also generated by the subgroups~${\rm SL}_2^i$, $i\in I$.

The operators $\phi(h)\in {\rm End}(L)$ for $h\in\h$ are semisimple so they can also be exponentiated to give a commutative group of operators $T_\C^\phi = T_\C^\phi(G) = \big\{\chi^\phi_h(t) = \exp(\phi(th))\, |\, h\in\h, t\in\C\big\}\leq G^\phi$ which is called the standard maximal torus of $G^\phi$. We also define the standard Borel subgroups $(B^\phi)^\pm = T_\C^\phi\big\la U^\phi_\a \, |\, \a\in(\Phi^{\rm re})^\pm\big\ra$ and the normalizer of $T_\C^\phi$ denoted by $N_\C^\phi$.

It is well known that the operators
\begin{gather*} \wtil^{\rm ad}_i = \exp({\rm ad}_{e_i})\exp({\rm ad}_{-f_i})\exp({\rm ad}_{e_i}) = \exp({\rm ad}_{-f_i})\exp({\rm ad}_{e_i})\exp({\rm ad}_{-f_i}) ,\qquad 1\leq i\leq\ell,\end{gather*}
in $G^{\rm ad}$, generate a subgroup
$\Wtil^{\rm ad}$ in $G^{\rm ad}$ such that the restriction of $\wtil^{\rm ad}_i$ to the standard Cartan subalgebra $\fh$ equals the simple Weyl group reflection~$w_i$ and $\wtil^{\rm ad}_i(e_i) = -f_i$. It means that $W$ is a homomorphic image of $\Wtil^{\rm ad}$.
Note that $\Wtil^{\rm ad}$ is a subgroup of $N_\C^{\rm ad}$ and that $N_\C^{\rm ad}/T_\C^{\rm ad} \cong W$. In Theorem~\ref{theorem:newWeylFormula}
we prove a formula in any integrable representation $\phi$ for
\[\wtil^{\phi}_i = \exp(\phi(e_i)) \exp(\phi(-f_i) \exp(\phi(e_i)) = \exp(\phi(\pi(e_i - f_i)/2)).\]

\subsection[Twin $BN$-pair and twin Tits building of a minimal Kac--Moody group]{Twin $\boldsymbol{BN}$-pair and twin Tits building of a minimal Kac--Moody group}\label{subsection:BNpair_Tits_building}

Our references for this section are~\cite[Sections 6.2 and 6.3]{AbramenkoBrown08} and~\cite[Chapters~5 and~11]{Ronan89}.
\begin{Definition}[BN-pair]\label{definition:BN-pair} A group $G$ is said to have a $BN$-pair if $G$ has subgroups $B$ and $N$ such that
\begin{itemize}\itemsep=0pt
\item[T1:] $G=\la B,N\ra$, $T = B\cap N \lhd N$, $W = N/T$ is generated by a set $S$.
\item[T2:] For $s\in S$ and $w\in W$ we have $sB w\subseteq B s wB \cup B wB\,$.
\item[T3:] For $s\in S$ we have $sB s^{-1}\not\subseteq B\,$.
\end{itemize}
The group $W$ is called the {\it Weyl group} of the $BN$-pair, and $(G,B,N,S)$ is called a {\it Tits system}. Furthermore, $(W,S)$ is a {\it Coxeter system} and there is a length function $\ell\colon W\to\N$.
\end{Definition}
\begin{Definition}[twin-BN-pair] \label{definition:Twin-BN-pair}
A group $G$ is said to have a twin $BN$-pair with Weyl group $W$ if $G$ has subgroups $B^+$, $B^-$ and $N$ such that
\begin{itemize}\itemsep=0pt
\item[TW1:] $\big(G,B^\pm,N,S\big)$ is a Tits system.
\item[TW2:] If $\ell(sw) < \ell(w)$ for $s\in S$ and $w\in W$, then $B^\pm s B^\pm w B^\mp = B^\pm sw B^\mp$.
\item[TW3:] $B^+ s \cap B^- = \varnothing$.
\end{itemize}
In this case, $(G,B^+,B^-,N,S)$ is called a {\it twin Tits system}.
\end{Definition}

For a hyperbolic adjoint Kac--Moody group $G = G^{\rm ad}_\C(A)$ we have standard Borel subgroups,
$B^{\pm} = (B^{\rm ad})^\pm$, the standard maximal torus, $T = T^{\rm ad}_\C$, and its normalizer in $G$,
$N = N_\C^{\rm ad}$. Thus the group $T = N\cap B^{\pm}$ is a normal subgroup of $N$.
The group $W = N_G(T)/T$ generated by $S = \big\{\wtil^{\rm ad}_i\,|\, 1\leq i\leq\ell\big\}$, is isomorphic to our
earlier definition in Section~\ref{subsection:Kac--Moody-algebra}
of the Weyl group $W$ generated by simple reflections as a group of orthogonal transformations of $\h^*$ given by formula~(\ref{equation:simple_reflection}). Thus, we have a twin $BN$-pair or twin Tits system, for the hyperbolic adjoint Kac--Moody group~$G$.

We have the (positive and negative) standard Borel subgroups, corresponding to the standard Borel subalgebras, $B^{\pm}=T(G)U^{\pm}$ where $U^{+}$ is generated by all positive real root groups and~$U^-$ is generated by all negative real root groups. The $BN$-pairs $\big(B^+,N\big)$ and $\big(B^-,N\big)$ have Birkhoff and Bruhat decompositions:
\[G = \coprod_{w\in W} B^{\pm}wB^{\mp} = \coprod_{w\in W}B^{\pm}wB^{\pm}.\]
In these double cosets, $w\in W = N/T$ is a coset $nT$ for $n\in N$, but for any two representatives of the same coset, $w = nT = n'T$, we have
\[B^{\pm}nB^{\pm} = B^{\pm}nTB^{\pm} = B^{\pm}n'TB^{\pm} = B^{\pm}n'B^{\pm},\]
so we can label a double coset by $w\in W$.

We can use the Bruhat decomposition to define a $W$-valued distance function on $G/B^\pm$,
\[\delta^\pm \colon \ G/B^\pm \times G/B^\pm \to W\]
by $\delta^\pm\big(g_1 B^\pm, g_2 B^\pm\big) = w$ when $g_1^{-1} g_2 \in B^{\pm}wB^{\pm}$. Similarly, we can use the
Birkhoff decomposition to define a $W$-valued codistance function
\[\delta^* \colon \ G/B^\pm \times G/B^\mp \to W\]
by $\delta^*\big(g_1 B^\pm, g_2 B^\mp\big) = w$ when $g_1^{-1} g_2 \in B^{\pm}wB^{\mp}$.

A proper subgroup $P^{\pm}$ of $G$ is called {\em parabolic} when it contains a conjugate of a Borel subgroup $B^{\pm}$, and it is called \emph{positive} or \emph{negative}, depending on the sign. For each subset $J\subset I$ define the subgroup
 $W_J = \la w_j \, |\, j\in J\ra$ of $W$ and the corresponding subgroups of $G$,
\[P_J^\pm\quad = \quad \coprod_{w\in W_J}B^{\pm}wB^{\pm} .\]
Note that $P_I^\pm = G$ is not parabolic since it is not proper,
$P_{\varnothing}^\pm = B^{\pm}$, and we write $P_i^\pm = P_{\{i\}}^\pm$. For $J\subsetneq I$ we
call $P_J^\pm$ a {\em standard parabolic subgroup}, and these form a complete set of representatives of the conjugacy classes of parabolic subgroups, so there are $2\cdot (2^\ell-1)$ conjugacy classes of parabolic subgroups. A parabolic subgroup $P^{\pm}$ is called \emph{maximal} if there is no parabolic subgroup $P'^{\pm}$ such that $P^{\pm}\subsetneq P'^{\pm}$. For each $i\in I$,
$P_{[i]}^\pm = P_{I\backslash\{i\}}^\pm$ is a maximal standard parabolic subgroup, so there are $2\ell$ conjugacy classes of maximal parabolic subgroups, $\ell$ positive and $\ell$ negative.

\begin{Definition}[Tits building] \label{definition:Tits_Building}A {\em Tits building} of type~$(W,S)$ consists of a simplicial complex~$\cB$ together with a collection~$\cA$ of subcomplexes,
each of which is called an apartment, such that
\begin{enumerate}\itemsep=0pt
 \item[1)] each apartment is a Coxeter complex for the Coxeter system $(W,S)$,
 \item[2)] each pair of {\em chambers}, i.e., simplices of maximal dimension in $\cB$, is contained in a common apartment,
 \item[3)] for two apartments $A$ and $A'$ there is an isomorphism $\varphi\colon A\rightarrow A'$, fixing the intersection $A\cap A'$.
\end{enumerate}
\end{Definition}
A Coxeter complex for $(W,S)$ is a simplicial complex on which there is a simply transitive action of a Coxeter group $W$ on
the simplices of maximal dimension $\ell -1$ (chambers). The simplices of dimension $\ell -2$ are called panels, and are
each labeled by a~generator, $s\in S$. We say that two chambers~$C_1$ and~$C_2$ are $s$-adjacent when $sC_1 = C_2$, which
means their intersection is an $s$-panel. Each element $w\in W$ is a product of generators from~$S$, so we define the {\it length}~$|w|$ to be the minimal number of generators in an expression for~$w$.
Suppose $C_i$ for $0\leq i\leq d$ is a~sequence of chambers such that $C_{i-1}$ and $C_i$ are
$s_i$-adjacent for $1\leq i\leq d$, so that $C_i = r_i C_0$ for $r_i = s_1\cdots s_i$. Then we define the $W$-valued distance function $\delta(C_0,C_d) = r_d = s_1\cdots s_d$, so for any $w_1,w_2\in W$ and any chamber $C$, we have $\delta(w_1 C,w_2 C) = w_1^{-1} w_2$. Since any two chambers in the building $\cB$ are in a common apartment, we have defined the $W$-valued distance function $\delta\colon \cC\times\cC\to W$ where $\cC = \cC(\cB)$ is the set of all chambers of~$\cB$.
One may choose a~pair~$(A,C)$ consisting of an apartment $A$ and a chamber $C$ in $A$, which we call {\it fundamental}, so that the chambers of~$A$ are uniquely labeled by the elements of~$W$.

From \cite[Proposition 4.84]{AbramenkoBrown08}, we have the following properties of $\delta$. For any chambers $C,C',D\in \cC(\cB)$ we have
\begin{enumerate}\itemsep=0pt
 \item $\delta(C,D) = 1$ iff $C = D$.
 \item $\delta(D,C) = \delta(C,D)^{-1}$.
 \item If $\delta(C',C) = s\in S$ and $\delta(C,D) = w\in W$, then $\delta(C',D) = sw$ or $\delta(C',D) = w$. If, in addition,
 $|sw| = |w| + 1$ then $\delta(C',D) = sw$.
 \item If $\delta(C,D) = w$ then for any $s\in S$, there exists a chamber $C'$ such that $\delta(C',C) = s$ and $\delta(C',D) = sw$.
 If $|sw| = |w| -1$ then there exists a unique such $C'$.
\end{enumerate}

From \cite[Definition 5.133]{AbramenkoBrown08}, we have the following definition of a twin building.
\begin{Definition}[twin building]\label{definition:Twin_Building}
A {\em twin building} of type $(W,S)$ is a triple $(\cB^+, \cB^-, \delta^*)$ consisting of two buildings $(\cB^+, \delta^+)$ and $(\cB^-, \delta^-)$, each of type $(W,S)$ and each with its own $W$-valued distance function, along with a {\it codistance} function
\[\delta^* \colon \ \big(\cC^+ \times\cC^-\big) \cup \big(\cC^- \times\cC^+\big) \to W\]
where $\cC^\pm$ is the set of chambers of $\cB^\pm$, satisfying the following conditions for each $\eps\in\{+,-\}$, any $C\in\cC^\eps$,
and any $D\in\cC^{-\eps}$, where $w = \delta^*(C,D)$:
\begin{enumerate}\itemsep=0pt
 \item[(Tw1)] $\delta^*(C,D) = \delta^*(D,C)^{-1}$.
 \item[(Tw2)] If $C'\in\cC^\eps$ satisfies $\delta^\eps(C',C) = s\in S$ and $|sw| < |w|$ then $\delta^*(C',D) = sw$.
 \item[(Tw3)] For any $s\in S$ there exists a chamber $C'\in\cC^\eps$ with $\delta^\eps(C',C) = s$ and $\delta^*(C',D) = sw$.
\end{enumerate}
\end{Definition}

Let us state Lemma 5.139 from~\cite{AbramenkoBrown08} for later reference:

\begin{Lemma}\label{lemma:longer_element_chamber_unique}
With the notation above:
\begin{enumerate}\itemsep=0pt
\item[$1.$] $\delta^*(C',D) \in \{w, sw\}$ for all $C'\in\cC^\eps$ with $\delta^\eps(C',C) = s$.
\item[$2.$] If $l(sw) > l(w)$, then there exists precisely one chamber $C'\in\cC^\eps$ satisfying $\delta^\eps(C',C) = s$ and $\delta^*(C',D) = sw$.
\end{enumerate}
\end{Lemma}

This leads to the following definition of the {\it opposition} relation between $\cC^+$ and $\cC^-$.
\begin{Definition}[opposition relation]\label{definition:Opposition_Relation}
For $C\in\cC^\eps$ and $D\in\cC^{-\eps}$ we say $C$ and $D$ are opposite, denoted by $C\ {\rm op}\ D$, when $\delta^*(C,D) = 1$.
\end{Definition}

Finally we define a {\it twin apartment} in a twin building as in \cite[Definition~5.171]{AbramenkoBrown08}.
\begin{Definition}[twin apartment]\label{definition:Twin_Apartment}
A twin apartment in a twin building $\big(\cB^+, \cB^-, \delta^*\big)$ is an ordered pair $\big(A^+, A^-\big)$ where $A^\pm\in \cA^\pm$ is an apartment in $\cB^\pm$ such that every chamber in $A^+\cup A^-$ is opposite to precisely one chamber in $A^+\cup A^-$.
\end{Definition}

From now on let $\cB$ denote the twin building associated to a Kac--Moody group $G$.
The simplices of $\cB$ are in bijection with parabolic subgroups in such a way that simplices in~$\cB^+$ correspond to positive parabolic subgroups and simplices in $\cB^-$ correspond to negative parabolic subgroups. The vertices ($0$-simplices) of the twin building $\cB$ are in bijection with maximal parabolic subgroups in~$G$. Chambers are in bijection with positive and negative Borel subgroups.
 We will denote these simplices by their corresponding parabolic subgroups.
Since parabolic subgroups are self-normalizing, the simplices in the building can be equivalently indexed by the coset spaces $G/P_J^\pm$, where $\big\{P_{J}^\pm\, |\, J\subsetneq I\big\}$ is a complete set of representatives of the conjugacy classes of parabolic subgroups.

The incidence relation on the set of vertices is given by intersections of parabolic subgroups as follows. For $0\leq r\leq \ell-1$ the $r+1$ vertices $P_{[i_1]}^\pm,\dots , P_{[i_{r+1}]}^\pm$ span an $r$-simplex if and only if the intersection
$P_{[i_1]}^\pm\cap \dots \cap P_{[i_{r+1}]}^\pm = P_{[i_{1},\dots, i_{r+1}]}^\pm$ is a parabolic subgroup, so~$\cB$ is a simplicial complex of dimension $ \dim(\cB) =\ell -1$.
For non-finite type Kac--Moody groups such as the hyperbolic type we consider, $B^+$ and $B^-$ are not conjugate in~$G$, so the intersection of a positive parabolic subgroup with a negative parabolic subgroup never contains a Borel subgroup.
Hence, from the twin $BN$-pair $\big(B^+,B^-,N\big)$ we get two buildings $\big(\cB^+,\delta^+\big)$ and $(\cB^-,\delta^-)$, each equipped with a~$W$-valued distance function, $\delta^\pm$,
each of type $(W,S)$, and a codistance function $\delta^*$, yielding a~twin building $\cB=\big(\cB^+, \cB^-, \delta^*\big)$. In $\cB^\pm$ the $(\ell-1)$-simplex $P_{[1,\dots, {\ell}]}^\pm = P_{[I]}^\pm = P_\varnothing^\pm = B^\pm$ is called the {\em fundamental chamber} of $\cB^\pm$. Each fundamental chamber has boundary consisting of the simplices $\Delta_J^\pm = P_J^\pm$, and has closure
\[\Delta^\pm = \bigcup_{J\subsetneq I} P_J^\pm.\]
We also choose a {\it fundamental apartment} $A^\pm_{\rm fund}$ in $\cB^\pm$ to be the one whose chambers are
\[\big(NB^\pm\big)/B^\pm = \big\{wB^\pm\,|\, w\in W\big\}.\]
Note that $\delta^\pm\big(w_1B^\pm,w_2B^\pm\big) = w_1^{-1} w_2$ and $\delta^*\big(w_1B^\pm,w_2B^\mp\big) = w_1^{-1} w_2$, so that $\big(A^+_{\rm fund},A^-_{\rm fund}\big)$ is a twin apartment.

Using the property that the simplices in each building are in bijection with the union of the coset spaces
$\bigcup_{J\subsetneq I} G/P_{J}^\pm$, we describe the buildings $\cB^+$ and $\cB^-$ associated to a twin $BN$-pair, $\big(B^+, B^-,N\big)$ for a Kac--Moody group $G$ as follows:
\begin{gather}\label{equation:buildingforKacMoodygroup}
\cB^{\pm}:=\big(G/B^{\pm}\times \Delta^\pm\big)/{\sim}.
\end{gather}

The equivalence relation $\sim$ is defined by $\big(fB^{\pm}, \Delta_J^\pm\big) \sim \big(gB^{\pm}, \Delta_{J'}^\pm\big)$ in $\big(G/B^{\pm}, \Delta^\pm\big)$ if and only if $\Delta_J^\pm = \Delta_{J'}^\pm$ (so $J = J'$) and $f^{-1}gP_J^{\pm}\subset P_J^\pm$.

Hence on the chambers $\Delta_{\varnothing}^\pm$, the equivalence relation $\sim$ is trivial, while on simplices in the boundary it is nontrivial.

Let $\phi\colon G\rightarrow G$ be an involution centralizing the Weyl group and such that $\phi\big(B^\pm\big)=B^\mp$. Then $\phi$ induces a twin building involution as follows (for details see~\cite{deMedtsGramlichHorn09, Horn09})
\begin{gather}\label{equation:twin_building_involution}
\phi\big(gB^\pm, \Delta_J^\pm\big)=\big(\phi(g)B^\mp, \Delta_J^\mp\big).
\end{gather}

One can give a {\em geometric realization} of a building as follows. Let $\{e_1,\dots,e_\ell\}$ denote the standard orthonormal basis of~$\R^\ell$. Each $r$-simplex is identified with a copy of the standard simplex
\[\Delta^r := \left\{x = \sum_{i=1}^{r+1} a_i e_i \in \R^{r+1}\, \big|\, 0\leq a_i\leq 1, \, \sum_{i=1}^{r+1} a_i=1 \right\}\]
which inherits the topology from $\R^{r+1}$. Appropriate identifications must be made among the copies of these standard simplices in order to reflect the incidence structure among the simplices in the building. For details see~\cite{AbramenkoBrown08}.

For the buildings associated with the hyperbolic Kac--Moody groups we wish to study, the geometric realization of apartments in the buildings $\cB^+$ and $\cB^-$ can be chosen to be isometric to hyperbolic spaces tessellated by the action of the hyperbolic Weyl group~$W$. In the case when the Cartan matrix is strictly hyperbolic, that tessellation is by compact simplices, but otherwise these simplices have ideal vertices stabilized by affine type subgroups of~$W$.

\subsection{Compact real forms of Kac--Moody algebras and groups}\label{section:Campact_real_forms_of_KM_algebras_and_groups}

Let $\g=\g_{\C}(A)$ be a complex Kac--Moody algebra and let $\h$ be the standard Cartan subalgebra. The {\it Cartan involution}
\[{\omega}_0\colon \ \g\longrightarrow \g\]
is the automorphism of $\g$ determined by ${\omega}_0(e_i)=-f_i$, $\omega_0(f_i)=-e_i$ and ${\omega}_0(h_i)= -h_i$. Composing ${\omega}_0$ with complex conjugation, we obtain a conjugate linear involution $\omega$, called the {\it Cartan--Chevalley involution}. Then $\k=\{x\in\g\,|\, \omega(x)=x\}$ is a Lie algebra over $\R$ called the {\it compact real form} of~$\g$ \cite[p.~243]{Kumar02}.

Note that $\omega_0$ and $\omega$ both centralize the Weyl group, so they each induce twin building involutions via formula~(\ref{equation:twin_building_involution}).

We may give generators for the compact real form $\k$ as follows.
For each $j=1,\dots ,\ell$, let $\g_j = \fsl_2^j$ be the Lie subalgebra of $\g$ isomorphic to $\fsl_2(\C)$ with basis $\{e_j, f_j, h_j\}$, so that
$\g$ is generated by the subalgebras $\g_j$ and $\omega(\g_j)= \g_j$. Then, using the notation $\bi = \sqrt{-1}\in\C$,
the real Lie algebra of fixed points of $\omega$ on $\g_j$, $\k_j=\fsu_2^j$ has basis
\[x_j = \frac{1}{2}(e_j - f_j),\qquad y_j = \frac{\bi}{2}(e_j + f_j), \qquad z_j = \frac{\bi}{2}(h_j)\]
with brackets $[x_j,y_j] = z_j$, $[y_j,z_j] = x_j$, $[z_j,x_j] = y_j$, and the compact real form $\k$ is generated by all of the subalgebras
$\k_j$, $j=1,\dots ,\ell$ (see \cite[Proposition~1]{Berman85}). A Cartan subalgebra in the compact real form $\k$ is an abelian subalgebra whose complexification is a~Cartan subalgebra in~$\g$. The standard Cartan subalgebra $\ft =\h\cap \k$ in $\k$ has real basis $\{z_j\, |\, 1\leq j\leq \ell\}$.

Let $G=G_\C(A)$ be the complex adjoint Kac--Moody group associated to $\g$. The involution $\omega$ of $\g$ lifts to a unique involution of $G$, also denoted by $\omega$, exchanging positive and negative real root groups since $\omega(\fg_\a) = \fg_{-\a}$. We have the following more general lemma about the action of any Lie algebra automorphism, which we will apply to $\omega$ as well as to $\wtil^{\rm ad}_i\in\Wtil^{\rm ad}$.

\begin{Lemma}\label{lemma:auto_action_on_real_root_groups}
For any $\a\in\roots^{\rm re}$, $e_\a\in\fg_\a$ and any $\phi\in {\rm Aut}(\fg)$, we have
$\phi\circ \exp({\rm ad}_{e_\a}) \circ\phi^{-1} = \exp({\rm ad}_{\phi(e_\a)})$ so that, in particular,
$\omega U^{\rm ad}_\a \omega^{-1} = U^{\rm ad}_{-\a}$ and $\wtil^{\rm ad}_i U^{\rm ad}_{\a_i} \big(\wtil^{\rm ad}_i\big)^{-1} = U^{\rm ad}_{-\a_i}$ for $1\leq i\leq\ell$.
\end{Lemma}

\begin{proof}For any $x\in\fg$, since $\phi$ is a Lie algebra automorphism, we have
\[\phi(\exp({\rm ad}_{e_\a})x) = \phi\left( \sum_{k\geq 0} \frac{1}{k!} ({\rm ad}_{e_\a})^k (x) \right)
= \sum_{k\geq 0} \frac{1}{k!} ({\rm ad}_{\phi(e_\a)})^k (\phi(x)) = \exp({\rm ad}_{\phi(e_\a)}) (\phi(x)).\]
The formula with $\phi = \omega$ gives $\omega U^{\rm ad}_\a \omega^{-1} = U^{\rm ad}_{-\a}$ since $\omega(\fg_\a) = \fg_{-\a}$, and with $\phi = \wtil^{\rm ad}_i$ gives $\wtil^{\rm ad}_i U^{\rm ad}_{\a_i} \big(\wtil^{\rm ad}_i\big)^{-1} = U^{\rm ad}_{-\a_i}$ since $\wtil^{\rm ad}_i(e_i) = -f_i$.
\end{proof}

We set $K=\operatorname{Fix}_G(\omega)$. Then $K$ is called the {\it unitary form} or {\it compact real form} of $G$. We will use the latter by analogy with the finite-dimensional case, even though $K$ is not compact. The group $K$ is generated by subgroups $K_j$ such that $\k_j = \operatorname{Lie}(K_j)$ for each $j=1,\dots ,\ell$ \cite{Caprace09, KacPeterson85c, Tits86b}.

For each $v = a + b\bi\in\C$ and $1\leq j\leq\ell$ we write a generator of $T_\C^{\rm ad}$ as
\[\exp({\rm ad}_{vh_j}) = \exp({\rm ad}_{ah_j}) \exp({\rm ad}_{b\bi h_j}) = \exp({\rm ad}_{ah_j}) \exp({\rm ad}_{b2 z_j}).\]
This gives the decomposition $T_\C^{\rm ad} = T_\R^{\rm ad}\, T$ where
\begin{equation}\label{equation:split_real_torus}
T_\R^{\rm ad} = \la \exp({\rm ad}_{ah_j}) \, |\, a\in\R,\, 1\leq j\leq\ell\ra
\end{equation}
is the split real torus and $T = \la \exp({\rm ad}_{b z_j}) \, |\, b\in\R, \, 1\leq j\leq\ell\ra$ is the compact real torus.

It is clear from the two expressions for $\wtil^{\rm ad}_i$ that $\Wtil^{\rm ad}\leq K$, but it is not so obvious that these operators can be expressed as
a single exponential $\wtil^{\rm ad}_i = \exp({\rm ad}_{\pi\, x_i})$. This is a special case of the formula $\wtil^{\phi}_i = \exp(\phi(\pi\, x_i))$ for any
integrable representation $\phi\colon \fg\to {\rm End}(V)$, where the inner $\pi$ is in $\R$, proven in Theorem \ref{theorem:newWeylFormula}
and first found in certain important cases by~\cite{DamourHillmann09}.
See also \cite{HainkeKoehlLevy15,KleinschmidtNicolai13}.

\begin{Proposition}
\label{proposition:Cartancompactconjugate}
All Cartan subalgebras of $\k$ are conjugate under the action of $K$.
\end{Proposition}

This result follows from \cite[Proposition~8.1(iii)]{Caprace09},. See also \cite[Corollary 5.33]{KacWang92} and~\cite[Proposition 3.5]{KacPeterson87}.

\section[Tits cone and lightcone of hyperbolic Kac--Moody algebras of compact type]{Tits cone and lightcone of hyperbolic Kac--Moody algebras\\ of compact type}

Let $\g = \g_\C(A)$ be a hyperbolic Kac--Moody algebra and let $\h$ be its standard Cartan subalgebra.
The split real form $\g_\R = \g_\R(A)$ of $\g$ contains the $\R$-subalgebra $\h_{\R}=\h\cap \g_{\R}$ which is just the $\R$-span of $\{h_i \,|\, i\in I\}$. We call $\h_{\R}$ the standard Cartan subalgebra of $\g_\R$.

The Weyl group $W$ has been defined above as the group of orthogonal transformations of~$\h^*$ generated by the simple reflections. The isomorphism between $\h^*$ and $\h$ gives the action of~$W$ on~$\h$, where the formula for simple reflections is just $w_i(h_j) = h_j - \alpha_j(h_i)h_i$. This same formula restricted to $\h_{\R}$ gives the action of $W$ and it also determines the action of $W$ on $\ft$ by $w_i(z_j) = z_j - \alpha_j(h_i)z_i$. These operators are orthogonal with respect to the restriction of the bilinear form $(\cdot,\cdot)$ to $\h_\R$ and to $\ft$, and therefore $W$ preserves each surface of constant
square length, $(\h_\R)_r = \{x\in\h_\R \, |\, (x,x) = r\}$ and $\ft_r = \{x\in\ft \, |\, (x,x) = r\}$.

Since we are assuming that the Cartan matrix $A$ is hyperbolic, the form $(\cdot,\cdot)$ is Lorentzian on~$\h_\R$ and on~$\ft$. The surface where $r=0$ is called the {\it nullcone}, each surface where $r < 0$ is called {\it timelike}, and each surface where $r > 0$ is
called {\it spacelike}. The set of all timelike points has two connected components, one called {\it forward} and denoted $TL^+$ and the other called {\it backward} and denoted $TL^-$. We have $TL^- = - TL^+$.

Each of these components is preserved by the linear action of $W$, which acts consistently on rays, ${\rm Ray}_x = \{rx\, | \, 0 < r\in\R\}$ since
$w(rx) = r w(x)$. A fundamental domain for the action of $W$ on each of the timelike components is defined by
\[\cC^\pm = \{x\in TL^\pm\, |\, \alpha_i(x) \geq 0,\, 1\leq i\leq \ell \}.\]
The union
\[X^\pm = \bigcup_{w\in W} w\big(\cC^\pm\big)\]
is called the {\it positive $($respectively negative$)$ Tits cone} and $X = X^+\cup X^-$ is called the {\it Tits cone}.

Clearly, we have $X^- = - X^+$. For $\g$ hyperbolic, $X^\pm \supseteq TL^\pm$, since it is possible that $\cC^\pm$, and therefore $X^\pm$, contains rays on the nullcone. This happens for the rank $3$ hyperbolic Cartan matrices
\[\bmatrix 2&-2&0 \\ -2&2&-1 \\ 0&-1&2 \endbmatrix \qquad\hbox{and}\qquad \bmatrix 2&-2&-2 \\ -2&2&-2 \\ -2&-2&2 \endbmatrix\]
corresponding to the hyperbolic Kac--Moody algebra $\cF$ \cite{FeingoldFrenkel83} whose Weyl group is the hyperbolic triangle group
$T(2,3,\infty)$, and the ``ideal'' hyperbolic Kac--Moody algebra $\cI$ whose Weyl group is the hyperbolic triangle group
$T(\infty,\infty,\infty)$, respectively.
Another such an example is $E_{10}$ because it contains the affine Kac--Moody algebra $E_{9} = E_8^{(1)}$.
For $A$ strictly hyperbolic, that is, whose principal minors are of finite type, we have $X^\pm = TL^\pm$.

We have the following description of the closure of the Tits cone (see \cite[equation~(5.10.2)]{Kac90})
\[\overline{X} = \overline{X^+}\cup \overline{X^-} = \{h\in \h_{\R}\,|\, (h, h)\leq 0\}.\]
We introduce the notations
\begin{gather*}
\cL_{\h_{\R}}=\{h\in \h_{\R}\,|\, (h, h)\leq 0\},\qquad
\cL_{\h_{\R}}^0=\{h\in \h_{\R}\,|\, (h, h)< 0\}\qquad \hbox{and}\\
\partial\cL_{\h_{\R}}=\{h\in \h_{\R}\,|\, (h, h)= 0\}.
\end{gather*}

\begin{Proposition} Let $\g$ be a hyperbolic Kac--Moody algebra over $\C$, $\k$ its compact real form and~$\g_\R$ its split real form. Let $\h$, $\h_\R$, $\ft$ be the standard Cartan subalgebras of $\g$, $\g_\R$, $\k$, respectively. Then $\h_{\R}\cong \R^{\ell-1,1}$ has signature $(\ell-1,1)$ and
$\ft= {\bi}\h_{\R} = \{{\bi}x\,|\, x\in \h_{\R}\}$, where ${\bi}^2=-1$. The signature of $\ft \cong \R^{1,\ell-1}$ is $(1,\ell-1)$.
\end{Proposition}
\begin{proof} The definition of a hyperbolic Cartan matrix gives that the bilinear form on the split real form of the Cartan subalgebra (the real span of the generators~$h_i$) is Lorentzian, so the statements in the proposition are clear.
\end{proof}

We use a sign convention on $\k$ adopted from the theory of finite-dimensional Riemannian symmetric spaces~\cite{Helgason01} and set
\begin{gather}\label{equation:compactmetric}
(\cdot,\cdot)_{\k} = - (\cdot,\cdot)|_{\k}.
\end{gather}
Note that in the affine case $(\cdot,\cdot)_{\k}$, this sign convention naturally occurs in the loop group realizations~\cite{Freyn09, PressleySegal86}. With respect to $(\cdot,\cdot)_{\k}$, the Cartan subalgebra $\ft$ has Lorentzian signature $(\ell-1,1)$.

\begin{Lemma}
The invariant symmetric bilinear form $(\cdot,\cdot)_{\k}$ on $\k$ is Lorentzian with signature $(\infty,1)$.
\end{Lemma}

\begin{proof} As a vector space, a complex hyperbolic Kac--Moody algebra $\g$ has a basis consisting of the Cartan subalgebra
generators, $h_i, i=1,\dots, \ell$, and certain multibrackets $[e_{i_1},e_{i_2},\dots,e_{i_n}]$ and $[f_{i_1},f_{i_2},\dots,f_{i_n}]$.
The action of $\omega$ on such multibrackets is simply $\omega([e_{i_1},e_{i_2},\dots,e_{i_n}]) = (-1)^n [f_{i_1},f_{i_2},\dots,f_{i_n}]$
and $\omega(h_i) = -h_i$ so an $\R$-basis for the compact real form $\k$ consists of the compact Cartan basis elements
$z_i =\frac{\bi}{2}h_i$, $i=1,\dots, \ell$, and the elements obtained from basis multibrackets above
\begin{gather*}
\frac{1}{2} ([e_{i_1},e_{i_2},\dots,e_{i_n}] + [f_{i_1},f_{i_2},\dots,f_{i_n}] ) \qquad \text{and}\\
 \frac{\bi}{2} ([e_{i_1},e_{i_2},\dots,e_{i_n}] - [f_{i_1},f_{i_2},\dots,f_{i_n}] )\qquad \text{for $n$ even},\\
\frac{1}{2} ([e_{i_1},e_{i_2},\dots,e_{i_n}] - [f_{i_1},f_{i_2},\dots,f_{i_n}] ) \qquad \text{and}\\
 \frac{\bi}{2} ([e_{i_1},e_{i_2},\dots,e_{i_n}] + [f_{i_1},f_{i_2},\dots,f_{i_n}] )\qquad \text{for $n$ odd}.
\end{gather*}
In particular, for $n = 1$ these include the elements $x_i=\frac{1}{2}(e_i-f_i)$ and
$y_i=\frac{\bi}{2}(e_i+f_i)$.

Following our sign convention, the Cartan subalgebra $\ft$ spanned by the elements $z_i$, $i=1,\dots, \ell$, has signature $(\ell-1,1)$. Using our sign convention and the characterization of the ${\rm ad}$-invariant scalar product in~\cite[equation~(2.2.1)]{Kac90}, applied to the basis vectors~$y_i$ and~$z_i$, we find that the bilinear form on each subalgebra $\k_i=\fsu_2^i$ is positive definite.

Furthermore the ${\rm ad}$-invariant bilinear form has the following properties on root spaces (see~\cite[Sections~2.1 and~2.2]{Kac90}). The root spaces $\g_{\alpha}$ and $\g_{-\alpha}$, which are interchanged by $\omega$, have dual bases
 \[\big\{e_{\alpha}^{j}\, |\, 1\leq j\leq \dim (\g_{\alpha})\big\} \qquad \text{and} \qquad
 \big\{f_{\alpha}^{j} = e_{-\alpha}^{j}\, |\, 1\leq j\leq \dim (\g_{\alpha}) \big\}, \qquad \text{respectively}\]
 such that
\[\big( f_{\alpha}^j, e_{\beta}^m\big) =\delta_{j,m}\delta_{\alpha, \beta} .\]

As a consequence, for positive roots $\alpha$, the basis elements of $\k$, $x_{\alpha}^j = \frac{1}{2}\big(e_{\alpha}^j-f_{\alpha}^j\big)$ and
$y_{\alpha}^j=\frac{{\bi}}{2}\big(e_{\alpha}^j+f_{\alpha}^j\big)$ satisfy
\[\big(x_{\alpha}^j, x_{\beta}^m\big) = -\frac{1}{2}\delta_{j,m}\delta_{\alpha, \beta} = \big(y_{\alpha}^j,y_{\beta}^m\big) \qquad \text{and} \qquad
\big(x_{\alpha}^j, y_{\beta}^m\big) = 0.\]
Hence it follows from~(\ref{equation:compactmetric}) that
\[\big(x_{\alpha}^j, x_{\beta}^m\big)_\k = \frac{1}{2}\delta_{j,m}\delta_{\alpha, \beta} = \big(y_{\alpha}^j,y_{\beta}^m\big)_\k \qquad \text{and} \qquad
\big(x_{\alpha}^j, y_{\beta}^m\big)_\k = 0.\]

The Cartan subalgebra has signature $(\ell-1,1)$ and for each positive root $\alpha$ the subspaces
\[\big\la x_{\alpha}^{j}\, |\, 1\leq j\leq \dim (\g_{\alpha})\big\ra \qquad \text{and} \qquad
 \big\la y_{\alpha}^{j} \, |\, 1\leq j\leq \dim (\g_{\alpha}) \big\ra\]
in $\k$ each have positive signature and are orthogonal to each other. For different $\alpha$ these spaces are orthogonal to each other and to the Cartan subalgebra $\ft$. Thus $\k$ has Lorentzian signature $(\infty,1)$.
\end{proof}

\begin{Remark} Note that the bilinear form $(\cdot,\cdot)$ on $\g_\R$ is indefinite with signature $(\infty,\infty)$ because, while the split Cartan $\h_{\R}$ has the signature $(\ell-1,1)$, each pair of dual root vectors $\big\{e_{\alpha}^{j}, f_{\alpha}^{j}\big\}$ forms a hyperbolic plane, and for different positive $\alpha$ and distinct $j$ these planes are orthogonal.
\end{Remark}

We recall that the ${\rm ad}$-invariant bilinear form on $\h_{\R}$ extends $\C$-bilinearly to the complexification $\h_{\C}$. In particular, for $x\in\h_{\R}$ we have
\begin{gather*}
({\bi}x,{\bi}x)=-(x,x)=(x,x)_{\k}
\end{gather*}
so the map $\varphi\colon \h_\R \longrightarrow \ft$ given by $\varphi(x) = {\bi}x$ is an isometry. In addition, for $w\in W$ we have $\varphi(wx)=w\varphi(x)$, hence $\varphi$ is $W$-equivariant.

We introduce notations for certain subsets in $\k$:
\begin{alignat*}{3}
& \cL_{\k}=\{x\in \k\,|\, (x,x)_{\k}\leq 0\},\qquad && \cL_{\ft} = \cL_{\k}\cap \ft = \varphi(\cL_{\h_\R}),&\\
& \cL_{\k}^0=\{x\in \k\,|\, (x,x)_{\k}< 0\}, \qquad && \cL_{\ft}^0 = \cL_{\k}^0\cap \ft = \varphi\big(\cL_{\h_\R}^0\big),& \\
& \partial\cL_{\k}=\{x\in \k\,|\, (x,x)_{\k}= 0\}, \qquad && \partial\cL_{\ft} = \partial\cL_{\k}\cap \ft = \varphi(\partial\cL_{\h_\R}).&
\end{alignat*}

The ones inside $\ft$ are related by the $W$-invariant isometry $\varphi$ to the corresponding subsets defined earlier in $\h_\R$.
Furthermore, we have
\begin{gather*} \varphi\big(TL^\pm\big) = TL_\ft^\pm,\qquad\! \varphi\big(\cC^\pm\big) = \cC_\ft^\pm,\qquad\! \varphi\big(X^\pm\big) = X_\ft^\pm,\qquad\!
\varphi\big(\overline{X}^\pm\big) = \overline{X}_\ft^\pm,\qquad\! \varphi(\overline{X}) = \overline{X}_\ft,
\end{gather*}
which correspond to those subsets defined earlier in $\h_\R$. We call $X_{\ft}$ the Tits cone of $\ft$ and note that $\overline{X}_\ft = \cL_{\ft}$.

\begin{Remark}\label{remark:lightconecompactform}
We will refer to $\cL_{\k}$ as the {\em lightcone} of $\k$ and $\cL_{\ft}$ as the {\em lightcone} of $\ft$.
\end{Remark}

\section[Group actions of the compact real form $K$]{Group actions of the compact real form $\boldsymbol{K}$}

\subsection[The adjoint action of $K$ on $\k$]{The adjoint action of $\boldsymbol{K}$ on $\boldsymbol{\k}$}\label{subsection:globalstructureAdjointaction}

Recall that $K$ denotes the compact form of the complex adjoint Kac--Moody group $G$. We set
\[\fH = \bigcup_{k\in K} k \ft k^{-1} = \big\{k x k^{-1}\in \fk \,|\, k\in K, x \in \ft\big\}.\]
By Proposition~\ref{proposition:Cartancompactconjugate} all Cartan subalgebras of $\k$ are $K$-conjugate, thus the definition of $\mathfrak{H}$ is independent of the choice of $\ft$.

\begin{Proposition}\label{proposition:tangentspace}Let $\ft$ be any Cartan subalgebra of $\k$ and let $z = \bi h\in\ft$ for $h\in\h_\R$ satisfying $\alpha(h)\neq 0$ for all $\alpha\in\Phi^+$.
Then the subspace $[\k,z]$ of $\k$ has basis
\[\bigcup_{\alpha\in\Phi^+} \big\{ x_{\alpha}^{j}\, |\, 1\leq j\leq \dim (\g_{\alpha})\big\} \cup
 \big\{ y_{\alpha}^{j} \, |\, 1\leq j\leq \dim (\g_{\alpha}) \big\}.\]
\end{Proposition}

\begin{proof}For $\alpha\in\Phi^+$ and $1\leq j\leq \dim (\g_{\alpha})$ we have
\begin{gather*}
\big[x_\a^j, z\big] = - \big[\bi h, \shf(e_\a^j-f_\a^j)\big] = -\sbihf \big[h, \big(e_\a^j-f_\a^j\big)\big] = -\sbihf \a(h) \big(e_\a^j + f_\a^j\big) = -\a(h) y_\a^j,\\
 \big[y_\a^j, z\big] = - \big[\bi h, \sbihf \big(e_\a^j+f_\a^j\big)\big] = \shf \big[h, \big(e_\a^j+f_\a^j\big)\big] = \shf \a(h) \big(e_\a^j - f_\a^j\big) = \a(h) x_\a^j
\end{gather*}
and $[\ft,z] = 0$, so no basis vectors of $\ft$ are in $[\k,z]$.
\end{proof}

By analogy with the finite-dimensional and affine cases, the following proposition shows that in the hyperbolic case, the $K$-orbits, $K\cdot z = \big\{k z k^{-1}\in \fk \,|\, k\in K\big\}$ for each $z\in\ft$, intersect each Cartan subalgebra orthogonally. For $z\in \ft$ let $T_z (K\cdot z)$ be the {\it tangent space} of the submanifold~$K\cdot z$ at the point~$z$.

\begin{Proposition}\label{proposition:orbits orthogonal} Let $\ft$ be any Cartan subalgebra of $\k$.
The orbits $K\cdot z$ for $z\in \ft$ are orthogonal to $\ft$ with respect to the ${\rm ad}$-invariant bilinear form $(\cdot, \cdot)_{\k}$, that is, $(T_z(K\cdot z),\ft)_{\k}=0$.
\end{Proposition}

\begin{proof} For $z\in \ft$ we have
\[T_z (K\cdot z)=[\k,z] .\]
By definition, for $w\in T_z(K\cdot z)$, there is some $y\in \k$ such that $w=[y,z]$. Since the form $(\cdot, \cdot)_{\k}$ is ${\rm ad}$-invariant, for $z'\in \ft$ we obtain
\begin{gather*}
(w,z')_{\k}=([y,z],z')_{\k}=(y,[z,z'])_{\k}=(y,0)_{\k}=0.\tag*{\qed}
 \end{gather*} \renewcommand{\qed}{}
\end{proof}

\begin{Definition}[surface notation]\label{definition:surface_notation} Let $X\subset V$ be a subset of a real vector space $V$ equipped with a bilinear form $(\cdot, \cdot)$. For any real number $r$ we define a `sphere' of radius $r$ in $X$ by
\[X_{r}:=\{x\in X\, |\, (x,x)=r\}.\]
\end{Definition}

If $(\cdot, \cdot)$ is Lorentzian and $r<0$ then $X_r$ is a two-sheeted hyperboloid of constant curvature $\kappa=-\frac{1}{r^2}$ whose connected components we call $X_r^+$ and $X_r^-$.
In particular we will use $\ft_{r}\subset \cL_\ft^0 \subset\ft$, $\k_r\subset \cL_\k^0\subset \k$ and $\fH_{r}\subset \fH$ and note that $\k_0=\partial\cL_\k$ and $\ft_0=\partial\cL_\ft$. For $\ell > 2$ we find $\ft_{-1}^+\cong\ft_{-1}^-\cong \H^{\ell-1}$, $(\ell-1)$-dimensional hyperbolic space.

\begin{Definition}[lightlike closure]\label{definition:lightlike_closure} \looseness=-1 Let $V$ be a possibly infinite-dimensional real vector space equipped with a Lorentzian form $(\cdot, \cdot)$ and let
$\partial\cL_V = \{x\in V\, |\, (x,x)=0\}$ denote its nullcone. The {\em boundary at infinity} $B^\infty(\partial\cL_V)$ of the nullcone consists of all rays ${\rm Ray}_x$ for $0\neq x\in\partial\cL_V$. If $\dim(V)=2$ it consists of four points. If $\dim(V)>2$ it consists of two components, corresponding to the future timelike boundary $B^\infty(\partial\cL_V)^+$ and the past timelike boundary $B^\infty(\partial\cL_V)^-$, each of which can be identified with a sphere of dimension $\ell-2$ We define the {\em lightlike closure} of~$V$ by
\[\hV = V\cup B^\infty(\partial\cL_V).\]
\end{Definition}

For each $r<0$, define $B^\infty(V_r^\pm)$, the boundary at infinity of $V_r^\pm$, to be equivalence classes of geodesic rays, where two rays are equivalent if the distance between them is finite at all points. For details, see~\cite{Eberlein96}. Note that for each $r<0$,
$B^\infty(V_r^\pm)$ can be identified with $B^\infty(\partial\cL_V)^\pm$. Using this observation we define the lightlike closure of
$V_r^\pm$ to be
\begin{equation}
\hV_r^\pm= V_r^\pm\cup B^\infty(\partial\cL_V)^\pm.
\end{equation}

Since the bilinear form $(\cdot,\, \cdot)_{\k}$ is ${\rm ad}$-invariant, the adjoint action of $K$ on $\k$ preserves the surfaces $\fH_r\subset \fH$ defined in Definition~\ref{definition:surface_notation}.
For $k\in K$ we have $k\cdot \ft_r=(k\cdot \ft)_r$. There is an induced action of $K$ on $B^\infty(\partial\cL_\ft)$ and for each $r<0$ on
$B^\infty(\partial\cL_{\ft_r})^\pm$, as well as on ${\widehat \ft}$ and on ${\widehat \ft}_r^\pm$.

Hence, for each $r<0$, the adjoint action of $K$ on $\k$ induces a well--defined action on the surface:
\[K\colon \ \fH_{r}\longrightarrow \fH_{r},\qquad k\cdot x = {\rm Ad}_k(x) = kxk^{-1}\]
as well as on the lightlike closure $\hfH_r$.

Since all Cartan subalgebras are conjugate, we can define the following surjective map
\begin{displaymath}
\psi\colon \ K\times \ft_{r}\longrightarrow \fH_{r},\qquad \psi(k,x) = kxk^{-1}.
\end{displaymath}

Note that the choice of the standard Cartan subalgebra $\ft$ in the definition of $\psi$ does not restrict the generality.

Let $T = \exp(\ft) = T\big(G^{\rm ad}\big)\cap K$ be the torus associated to $\ft$ and let
\[N = N_K(T) = N_{G^{\rm ad}}\big(T\big(G^{\rm ad}\big)\big)\cap K\]
be the normalizer of $T$ in $K$. For $k\in K$, $t\in T$ and $u\in \ft$ we have:
\[{\rm Ad}_{kt} (u) = ktut^{-1}k^{-1} = kuk^{-1}={\rm Ad}_k(u).\]
Thus the adjoint action on $\ft$ is $T$-invariant and factors to the quotient space $K/T$ yielding a~surjective map
\begin{displaymath}
\psi\colon \ K/T\times \ft_r\longrightarrow \fH_{r},\qquad \psi(kT,u) = kuk^{-1}.
\end{displaymath}

For any Cartan subalgebra $\ft' = k\ft k^{-1}$ of $\k$, $\cL_{\ft'}$ coincides with the closure of the Tits cone of~$\ft'$,~${\overline X}_{\ft'}$. Hence ${\overline X}_{\ft'}$ is the closure of the cone $\{su \in \ft' \, |\, s>0, \, u\in \ft'_r \}$, for any fixed $r<0$, which includes the boundary
$\partial\cL_{\ft'}$. We distinguish two cases:
\begin{enumerate}\itemsep=0pt
\item When every proper Cartan submatrix of $A$ is of finite type (i.e., $A$ is $(\ell-1)$-spherical or {\em strictly hyperbolic}) then we have $X_{\ft'} = \cL_{\ft'}^0$ and there is a bijection between each surface~$(\ft'_r)^\pm$ for $r<0$ and the set of rays in $\cL_{\ft'}^0$. In this case
we do not need to consider the lightlike closure $\hfH_r$ in order to embed the building in it.
\item When $A$ contains an affine Cartan submatrix the fundamental chambers $\cC^\pm$ contain rays that accumulate at rays on the lightcone~$\partial\cL_{\ft'}$. So the corresponding points on the surface $(\ft'_{-1})^\pm$ accumulate at points on the boundary of the
$(\ell-1)$-dimensional hyperbolic space~$\H^{\ell-1}$. Therefore, in these cases we do have to consider the lightlike closure in order to
embed the building. In the example of $\cF$, the surfaces $(\ft'_{-1})^\pm\cong\H^2$ are isometric to
the Poincar\'e disk, and the tessellation by the hyperbolic Weyl group~$T(2,3,\infty)$ includes chambers which have an ideal vertex on the boundary. In such an example, the building would have $0$-simplices corresponding to those ideal vertices, so to achieve our goal of embedding the twin building of a hyperbolic algebra inside the lightcone of~$\k$, we must use~$\hfH_r$.
\end{enumerate}

Since $W=N/T\leq K/T$, the restriction to $\ft_r$ of the adjoint action of $N$ on $\k_r$ coincides with the Weyl group action of $W$ on $\ft_r$.

Restriction of the second coordinate of the domain of $\psi$ to either fundamental domain $\cC_r^\pm$ for the action of $W$ on $\ft_r$ gives the surjective maps
\begin{gather*}
\psi^\pm\colon \ K/T\times \cC_r^\pm\longrightarrow \fH^\pm_r, \qquad \psi^\pm(kT,u) = kuk^{-1}.
\end{gather*}

\subsection{The local structure of the adjoint action}\label{section:localstructure}

In this section we describe the geometry `close' to a fixed Cartan subalgebra in the compact real form $\k$. We show that for $i\in \{1, \dots, {\ell}\}$, the adjoint action of the fundamental ${\rm SU}(2)_i$-subgroup of the compact real form $K$, is a rotation of the standard Cartan subalgebra around the hyperplane
\[L_{\ft,i}:=\{z\in \ft\, |\, \alpha_i(z)=0\}\]
fixed by the generator $w_i$ of $W$.

The standard Cartan subalgebra is given by $\ft= \R z_1 +\dots + \R z_{\ell}$.
For each $1\leq i \leq \ell$, ${\rm SU}(2)_i + \ft = {\rm SU}(2)_i \oplus L_{\ft,i}$ since $\a_i(z_i) = \bi$.
For any $s,t,u\in\R$, we have $\exp({\rm ad}_{sx_i+ty_i+uz_i})\in {\rm SU}(2)_i\leq K$ and for $z\in\ft$ we have
$ [sx_i+ty_i+uz_i,z] = -\bi\alpha_i(z) (tx_i-sy_i) = 0$ when $\alpha_i(z) = 0$ which means $\exp({\rm ad}_{sx_i+ty_i+uz_i})(z) = z$ when
$\alpha_i(z) = 0$. While $({\rm ad}_{sx_i+ty_i+uz_i})z = ({\rm ad}_{sx_i+ty_i})z$ for any $z\in\ft$, we see that
\begin{gather*} ({\rm ad}_{sx_i+ty_i+uz_i})^2 z = -\bi\alpha_i(z) ({\rm ad}_{sx_i+ty_i+uz_i})(tx_i-sy_i) = -\bi\alpha_i(z) \big(usx_i + uty_i - \big(s^2+t^2\big)z_i\big)
\end{gather*}
differs from $({\rm ad}_{sx_i+ty_i})^2 z$, and higher powers have increasingly complicated expressions involving~$u$, so that $\exp({\rm ad}_{sx_i+ty_i+uz_i})z$ certainly depends on $u$. Nevertheless, the factor $\alpha_i(z)$ tells us that $L_{\ft,i}$ is the fixed point set in $\ft$ of $\exp({\rm ad}_{sx_i+ty_i+uz_i})$, so the Cartan subalgebra $\exp({\rm ad}_{sx_i+ty_i+uz_i})\ft$ is
spanned by $L_{\ft,i}$ and the vector $\exp({\rm ad}_{sx_i+ty_i+uz_i})z_i$. It appears that this is a three-parameter family of Cartan subalgebras, each containing $L_{\ft,i}$, but, in fact, we can see that the entire family is obtained from the two-parameter family with $u = 0$.

While these operators are defined on the entire Kac--Moody algebra, $\fg$, we are really only interested in the orbit of $z_i$ inside $\fsu(2)_i$ under the action of ${\rm SU}(2)_i$, so we can do this calculation with $2\times 2$ matrices. An arbitrary matrix in ${\rm SU}(2)$ is
$A = \bm \a&\bet\\ -\bar\bet&\bar\a\ebm$ where $\a,\bet\in\C$ and $\det(A) = 1$. With $z = \frac{\bi}{2} \bm 1&0\\0&-1\ebm$ we have
\[AzA^{-1} = \frac{\bi}{2} \bm (\a\bar\a-\bet\bar\bet)&-2\a\bet \\-2\bar\a\bar\bet&-(\a\bar\a-\bet\bar\bet)\ebm.\]
So the stabilizer of $z$ consists of those $A$ such that $\a\bet = 0$ and $\a\bar\a-\bet\bar\bet = 1$, which implies that $\bet = 0$ and
the stabilizer of $z$ is the diagonal ${\rm U}(1) = \left\{\bm \a&0\\ 0&\bar\a\ebm \,|\, \a\bar\a = 1\right\}$. The orbit ${\rm SU}(2)\cdot z$ is in
bijection with ${\rm SU}(2)/{\rm U}(1)$ which is well-known to be the $2$-sphere. We will show below that the two-parameter family
$\exp({\rm ad}_{sx_i+ty_i})z_i$ gives a $2$-sphere, so it is enough to find all the Cartan subalgebras containing $L_{\ft,i}$.

For each $s,t\in\R$, $\exp({\rm ad}_{sx_i+ty_i})\ft = \ft^i(s,t)$ is either another Cartan subalgebra such that $\ft\cap\ft^i(s,t) = L_{\ft,i}$ or else
$\ft = \ft^i(s,t)$. With $v=\big(\frac{s}{2}+\frac{t\bi}{2}\big)\in\C$ we have $sx_i+ty_i = v e_i-\overline{v}f_i$, so
\[\exp({\rm ad}_{sx_i+ty_i}) = \exp({\rm ad}_{v e_i-\overline{v}f_i}) \in {\rm SU}(2)_i\]
gives another parameterization of $\{\ft^i(s,t)\, |\, s,t\in\R\} = \{ \exp({\rm ad}_{v e_i-\overline{v}f_i})\ft\, |\, v\in\C\}$. For fixed~$i$, distinct choices of $(s,t)\in\R$ can give the same Cartan subalgebra $\ft^i(s,t)$, so we would like to know exactly how to
parameterize the distinct Cartan subalgebras in this set. The first part of the following theorem shows that for each $1\leq i \leq \ell$ the family of distinct Cartan subalgebras obtained this way can be parameterized by a 2-sphere with antipodes identified, the real projective space~$P_2(\R)$. A more careful interpretation of this calculation will later give us information related to the set of chambers in the building which share a common panel.

\begin{Theorem} \label{theorem:expformulas} For any $s,t\in\R$ such that $0 < r^2 = s^2 + t^2$
and for any $z\in\ft$, we have
\begin{enumerate}\itemsep=0pt
\item[$1)$] $\exp({\rm ad}_{sx_i+ty_i}) z = z - {\bi}\alpha_i(z) (\cos(r)-1) z_i - {\bi}\alpha_i(z) \frac{\sin(r)}{r} (tx_i-sy_i)$,
\item[$2)$] $\exp({\rm ad}_{sx_i+ty_i}) x_i = x_i - \frac{t\sin(r)}{r} z_i - \frac{t}{r^2} (\cos(r) - 1) (tx_i-sy_i)$,
\item[$3)$] $\exp({\rm ad}_{sx_i+ty_i}) y_i = y_i + \frac{s\sin(r)}{r} z_i - \frac{s}{r^2} (\cos(r) - 1) (tx_i-sy_i)$.
\end{enumerate}
\end{Theorem}

\begin{proof}We prove only the first relation since the others are analogous. We have
\begin{gather*}
({\rm ad}_{sx_i+ty_i})^1 z = [sx_i+ty_i,z] = -\bi\alpha_i(z) (tx_i-sy_i) ,\\
({\rm ad}_{sx_i+ty_i})^2 z = [sx_i+ty_i,-\bi\alpha_i(z) (tx_i-sy_i)] = \bi\alpha_i(z) \big(s^2+t^2\big) z_i ,\\
({\rm ad}_{sx_i+ty_i})^3 z = \big[sx_i+ty_i,\bi\alpha_i(z) \big(s^2+t^2\big) z_i\big] = \bi\alpha_i(z) \big(s^2+t^2\big)(tx_i-sy_i), \\
({\rm ad}_{sx_i+ty_i})^4 z = \big[sx_i+ty_i,\bi\alpha_i(z) \big(s^2+t^2\big)(tx_i-sy_i)\big] = -\bi\alpha_i(z) \big(s^2+t^2\big)^2 z_i .
\end{gather*}
In the third step we have used that $\a_i(z_i) = \bi$.

Using $r^2 = s^2+t^2 \neq 0$, it is clear that for $n\geq 1$ we have
\[({\rm ad}_{sx_i+ty_i})^{2n} z = -\bi\alpha_i(z) (-1)^n \big(r^2\big)^{n} z_i\]
and for $n\geq 0$ we have
\[({\rm ad}_{sx_i+ty_i})^{2n+1} z = -\bi\alpha_i(z) (-1)^n \big(r^2\big)^{n} (tx_i-sy_i)\]
so we get
\begin{align*}
\exp({\rm ad}_{sx_i+ty_i}) z & = z-\bi\alpha_i(z) \sum_{n=1}^\infty \frac{(-1)^n r^{2n}}{(2n)!} z_i -\bi\alpha_i(z)
\sum_{n=0}^\infty \frac{(-1)^n r^{2n}}{(2n+1)!} (tx_i-sy_i)\\
& = z -\bi\alpha_i(z) (\cos(r)-1) z_i -\bi\alpha_i(z) \frac{\sin(r)}{r} (tx_i-sy_i).\tag*{\qed}
\end{align*}\renewcommand{\qed}{}
\end{proof}

\begin{Corollary} \label{corollary:hemispherical_parametrization}
For each $1\leq i \leq \ell$, the family of distinct Cartan subalgebras in $\{\ft^i(s,t)\, |\, s,t\in\R\}$, including $\ft$, is parametrized by the unit hemisphere
\[\big\{(\sin(r)\sin(\psi),-\sin(r)\cos(\psi),\cos(r))\in\R^3\, |\, 0\leq r < \pi,\, 0\leq\psi<\pi\big\},\]
where $r = \sqrt{s^2 + t^2}\geq 0$ and $\psi$ is defined when $r > 0$ by
$\sin(\psi) = \frac{t}{r}$ and $\cos(\psi) = \frac{s}{r}$. Also, $(s',t') = \frac{r-\pi}{r}(s,t)$ corresponds to the antipodal
point determined by $(\pi-r,\psi+\pi)$ and for any $z\in\ft$, we have
\[\exp( {\rm ad}_{s'x_i+t'y_i}) z = \exp( {\rm ad}_{sx_i+ty_i}) w_i(z).\]
\end{Corollary}

\begin{proof}
The Cartan subalgebra $\ft^i(s,t)$ is spanned by $L_{\ft,i}$ (for any choice of $(s,t)$) and the vector
\[\exp( {\rm ad}_{sx_i+ty_i}) z_i= \sin(r)\sin(\psi)x_i - \sin(r) \cos(\psi) y_i + \cos(r) z_i \in \fsu_2^i\]
from Theorem~\ref{theorem:expformulas}(1) and the fact that $\alpha_i(z_i) = \bi$.
With respect to basis $\{x_i,y_i,z_i\}$ of $\fsu_2^i$, the coordinates of this vector are
\[(\sin(r)\sin(\psi), - \sin(r) \cos(\psi), \cos(r))\]
so the vector is on a unit sphere. Each point on one hemisphere is uniquely determined by the choices
$0\leq r < \pi$ and $0\leq\psi<\pi$, and no two vectors of the above form are co-linear, so all of the corresponding subspaces $\ft^i(s,t)$ are distinct. The point on the sphere determined by $(r,\psi)$ corresponding to $(s,t)$ has antipodal point determined by $(\pi-r,\psi+\pi)$ corresponding to $(s',t') = \frac{r-\pi}{r}(s,t)$. The corresponding subspaces are the same since
$\exp( {\rm ad}_{s'x_i+t'y_i}) z_i = - \exp( {\rm ad}_{sx_i+ty_i}) z_i$.

We have $(s',t')$ determined by $sin(\psi+\pi) = \frac{t'}{\pi-r}$ and $\cos(\psi+\pi) = \frac{s'}{\pi-r}$ so $\frac{-t}{r} = -\sin(\psi) = \frac{t'}{\pi-r}$ and $\frac{-s}{r} = -\cos(\psi) = \frac{s'}{\pi-r}$, which gives $t' = \frac{r-\pi}{r} t$ and $s' = \frac{r-\pi}{r} s$. For $\alpha_i(z) = 0$, that is, for $z\in L_{\ft,i}$, we have $w_i(z) = z$ and so $\exp( {\rm ad}_{s'x_i+t'y_i}) z = z = \exp( {\rm ad}_{sx_i+ty_i}) w_i(z)$. For $z = z_i$ we have $w_i(z_i) = -z_i$, and we have shown above that $\exp( {\rm ad}_{s'x_i+t'y_i}) z_i = - \exp( {\rm ad}_{sx_i+ty_i}) z_i$. Since these operators are linear, that proves the claimed formula for any $z\in\ft$.
\end{proof}

\begin{Remark}We may find it useful later to associate $w\in\hat\C$ with a point $(x,y,z)$ on the unit sphere, and to record the relationship between antipodal points on the sphere and their associated complex numbers. The standard formula for a~stereographic projection from the north pole $(0,0,1)$ to a point in the complex plane $z = 0$ is
\[(x,y,z) \to \frac{x}{1-z} + \bi \frac{y}{1-z} = w\]
so for the antipodal point
\[(-x,-y,-z) \to \frac{-x}{1+z} + \bi \frac{-y}{1+z} = w^a.\]
Since $x^2 + y^2 + z^2 = 1$, we then have $w^a = -{\bar w}^{-1}$.
\end{Remark}

We wish to realize the twin apartments for a hyperbolic Kac--Moody group as a geometric object inside the lightcone of the compact real form of the Kac--Moody Lie algebra. It is not hard to see the Coxeter complex structure
in each side of the Tits cone of a standard compact real Cartan subalgebra. The tessellation of the forward and backward Tits cones by the action of the hyperbolic Weyl group provides a pair of opposite fundamental apartments, $(\Sigma^+,\Sigma^-)$, with tessellations $\Sigma^\pm = \bigcup\limits_{w\in W} w\big(\cC^\pm\big)$ where $\cC^\pm$ is a fundamental chamber in $\Sigma^\pm$. A~panel is a~nonempty intersection of maximal dimension of two $w_i$-adjacent chambers, that is, $\cC_1^\pm = w\cC^\pm$ and $\cC_2^\pm = w w_i \cC^\pm$, so that $w w_i w^{-1} \cC_1^\pm = \cC_2^\pm$. It means that the panel is contained in the wall of all fixed points of the Weyl group reflection $w w_i w^{-1}$ which is a $W$-conjugate of the simple reflection $w_i$. We may simplify the discussion by considering just the panels of a fundamental chamber, $\cC^\pm$, contained in the wall of fixed points of the simple reflection $w_i$. The hyperplane~$L_{\ft,i}$ intersects both sides of the Tits cone and determines a wall containing a panel of $\cC^\pm$ in each side of the twin apartment. The operators $\exp( {\rm ad}_{sx_i+ty_i})\in {\rm SU}_2^i$ for any $s,t\in\R$ fix $L_{\ft,i}$ pointwise, but take $z_i$ to a vector in $\fsu_2^i$ not in~$L_{\ft,i}$, giving a family of compact Cartan
subalgebras, $\exp( {\rm ad}_{sx_i+ty_i}) \ft = \ft^i(s,t)$, each of which contains a Tits cone, and all of which share the common subspace $L_{\ft,i}$. Consider the distinct chambers in this family of apartments which share a common panel with the fundamental chamber in the fundamental apartment (in either side of the twin fundamental apartment). For now, let us only think about $\Sigma^+$, and its images under this family of operators. The fundamental chamber, $\cC^+$, has a unique panel fixed by $w_i$, and that panel is contained in $L_{\ft,i}$. So the family of distinct chambers $\exp( {\rm ad}_{sx_i+ty_i}) \cC^+$ obtained as $s$ and $t$ vary, all contain that $w_i$-fixed panel. The same can be said of the chamber $w_i \cC^+$ in $\Sigma^+$, and we may wish to consider its orbit under this family of operators. But we know that with $(s,t) = (\pi,0)$, which corresponds to $(r,\psi) = (\pi, 0)$, that operator restricted to~$\ft$
equals $w_i$, so we should be able to understand that orbit starting from either chamber.

The point of the first part of Corollary~\ref{corollary:hemispherical_parametrization} is that for $(s,t)$ corresponding to $(r,\psi)$ giving distinct points on the hemisphere $0\leq r < \pi$, $0\leq\psi<\pi$, the real compact Cartan subalgebras~$\ft^i(s,t)$ are all distinct, so the chambers $\exp( {\rm ad}_{sx_i+ty_i}) \cC^+$ must all be distinct. But the second part of the corollary says that for $(s',t') = \frac{r-\pi}{r}(s,t)$ corresponding to the antipodal point determined by $(\pi-r,\psi+\pi)$,
and for any $z\in\ft$, we have $\exp( {\rm ad}_{s'x_i+t'y_i}) z = \exp( {\rm ad}_{sx_i+ty_i}) w_i(z)$. This means that $\exp( {\rm ad}_{s'x_i+t'y_i}) \cC^+ = \exp( {\rm ad}_{sx_i+ty_i}) w_i(\cC^+)$ is the chamber adjacent to $\exp( {\rm ad}_{sx_i+ty_i}) \cC^+$ in the apartment $\exp( {\rm ad}_{sx_i+ty_i}) \Sigma^+$, sharing the $w_i$-fixed panel. Thus, the complete set of all distinct
chambers $\exp( {\rm ad}_{sx_i+ty_i}) \cC^+$ is obtained when $(s,t)$ varies so that the corresponding $(r,\psi)$ gives all points on the unit sphere. Antipodal points give distinct chambers in the same apartment, and the set of all chambers sharing the common panel is parametrized by a real 2-sphere, exactly corresponding to the abstract building picture, where the answer is the Riemann sphere. From the remark above, if $w\in\hat\C$ is the label of a chamber sharing the common panel, the
label $w^a = -{\bar w}^{-1} \in\hat\C$ denotes the $w_i$ reflection of the first chamber.

\subsection{The action of the compact real form on the twin building}\label{subsection:action_of_compact_real_form}
To describe the action of $K$ on the twin building, we use the Iwasawa decomposition~\cite{deMedtsGramlichHorn09} which yields
\[G=KAU^{\pm},\]
where $G$ denotes a complex Kac--Moody group, $K$ denotes the compact real form of $G$, $A\cong \mathbb{R}^{{\rm rank}(G)}$ is an abelian subgroup and $U^{\pm}$ is the group generated by all positive (respectively negative) real root groups. Recall that $\ft$ denotes
a Cartan subalgebra in $\k$, and T$=exp(\ft)$ its torus. Using $T_\C(G) = TA$, the decomposition $B^{\pm}=T_\C(G) U^{\pm}$ and the Iwasawa decomposition, we have a bijection between coset spaces
 \begin{gather}\label{equation:GB=KT}
 G/B^{\pm} \leftrightarrow K/T.
 \end{gather}

Recalling the description of the building from equation~(\ref{equation:buildingforKacMoodygroup}) in Section~\ref{subsection:BNpair_Tits_building}
\begin{gather}\label{equation:buildingGequivariance}
\cB^{\pm}=\big(G/B^{\pm}\times \Delta^{\pm}\big)/{\sim},
\end{gather}
we obtain an equivalent new description by equation~(\ref{equation:GB=KT}):
\begin{gather}\label{equation:buildingk-invariance}
\cB^{\pm}=\big(K/T\times \Delta^{\pm}\big)/{\sim}.
\end{gather}

In the description of the twin building $\cB=\cB^+\cup \cB^-$ given by equation~(\ref{equation:buildingGequivariance}), the natural action of $G$ via left multiplication is apparent, while in the description of $\cB$ as in equation~(\ref{equation:buildingk-invariance}) the $G$-symmetry is broken to $K$-symmetry. Note that in equation~(\ref{equation:buildingGequivariance}) the cosets of the two opposite buildings $\cB^+$ and $\cB^-$ are defined with respect to different subgroups $B^+$ and $B^-$ respectively. Hence there are subgroups of~$G$ which act differently on $\cB^+$ and $\cB^-$. For example, for any $g\in G$ the coset $gB^+\in G/B^+$ is fixed by all elements in the subgroup
$B^+_g=gB^+g^{-1}$, but that subgroup acts transitively on $\big\{fB^-\in G/B^-\, |\, g^{-1}f\in B^+B^-\big\}$ which certainly contains~$gB^-$. But the description of $\cB^\pm$ in equation~(\ref{equation:buildingk-invariance}) shows that the action of $K$ is the same in both. In particular this means that the group $K$ does not act transitively on any apartment system. More precisely, for $A\in\cB^\pm$ an apartment, define the orbit $\cA_K(A) := K\cdot A$. Then for each chamber $c\in\cB^\pm$ there is exactly one apartment $A_c\in \cA_K(A)$ such that $c\in A_c$.

An apartment $A$ is called $\omega$-stable iff $\omega(A)=A$. If $A$ is $\omega$-stable, the set $\cA_K(A)$ contains all $\omega$-stable apartments.

Let $c = \big(fB^\pm, \Delta_\varnothing^\pm\big)$ be any chamber in $\cB^\pm=\big(G/B^\pm\times \Delta^\pm\big)/{\sim}$, and let $\cCham\big(\cB^\pm\big)$ denote the set of all chambers of the building $\cB^\pm$.
For $i\in I$, the $i$-panel of $c$ is $\big(fB^\pm, \Delta_i^\pm\big)$ and the $i$-residue of $c$ is defined to be
$R_i(c)=\big\{d = \big(gB^\pm, \Delta_\varnothing^\pm\big)\in\cCham\big(\cB^\pm\big)\, |\, \big(fB^\pm, \Delta_i^\pm\big) \sim \big(gB^\pm, \Delta_i^\pm\big) \big\}$. Then we have $R_i(c)\cong \P^1(\C)$. Identifying $\mathbb{P}^1(\C)$ with the Riemann sphere $\widehat{\C}=\C\cup \infty$, the action of the subgroup ${\rm SU}(2)_i$ on $R_i(c)$ can be identified with the action of ${\rm SU}(2)$ on $\widehat{\C}$ by M\"obius transformations. Additional details may be found in~\cite[Chapter~6]{AbramenkoBrown08},.

\section[Simplicial complex, distance and codistance on $\hfH_r$]{Simplicial complex, distance and codistance on $\boldsymbol{\hfH_r}$}\label{section:simplicial structure}

In this section we define for each $r<0$ a simplicial complex on the set $\hfH_r$, where we are using the notations in Definitions~\ref{definition:surface_notation} and~\ref{definition:lightlike_closure}. We begin in the standard Cartan subalgebra with $\ft_r\subset \ft$. Recall, that $\ft_{r} = \ft_{r}^+ \cup \ft_{r}^-$ is (up to rescaling) isometric to a pair of hyperbolic spaces (for $\ell>2$) which are both tessellated by the action of the Weyl group $W = \la S\ra$. We call an element $x\in \ft_r$ {\em singular} if it is fixed by a Weyl group element $w\neq 1$ and call $x$ {\em regular} otherwise. Let $\ft^{\sin}_r$ be the set of all singular elements in $\ft_r$ and let $\ft^{\rm reg}_r$ be the set of all regular elements, and similarly we have sets $\big(\ft^{\sin}_r\big)^\pm$ and $\big(\ft^{\rm reg}_r\big)^\pm$. Then we have the decomposition
\[\ft_r=\ft^{\sin}_r\cup \ft^{\rm reg}_r.\]
We denote by ${\rm Comp}_r^\pm$ the set of connected components of $\ft^{\rm reg}_r$.
In each sheet one connected component, $\cC^\pm_r = \cC^\pm\cap\big(\ft_{r}^{\rm reg}\big)^\pm$, corresponds to the fundamental domain $\cC^\pm$. For each sign $\pm$, the Weyl group acts simply transitively on that set ${\rm Comp}_r^\pm$. We use the Weyl group action to index the elements of ${\rm Comp}_r^\pm$ as follows: Let $1\in W$ denote the identity element. In each sheet~$\ft_r^\pm$ we index the
fundamental chamber~$\cC_r^\pm$ with $1$ such that $\cC_r^- = - \cC_r^+$.
Then we index the connected component $w\cC_r^\pm\in {\rm Comp}_r^\pm$ by
$w$ yielding
\[{\rm Comp}_r^\pm = \big\{w\cC^\pm_r \,|\, w\in W \big\}\subset \big(\ft^{\rm reg}_r\big)^\pm.\]
Let $\overline{w\cC^\pm_r}$ denote the closure of the component $w\cC^\pm_r$ and let $\cU^\pm$ denote the union
\[\cU^\pm = \bigcup_{w\in W} \overline{w\cC^\pm_r}\]
so $\cU^\pm$ covers $\ft_{r}^\pm$. It may happen that the closure $\overline{\cC^\pm_r}$ includes points at infinity, that is, ideal points not on the surface $\ft_{r}$, but which correspond to rays in the null cone as discussed in Definition~\ref{definition:lightlike_closure}. In that case, some vertices ($0$-simplices) will be in $\hft_r^\pm$.

For any subset $J\subsetneq S$ recall that $W_J = \la s\,|\, s\in J\ra$, and define the intersection
\[{\rm Simp}_J^\pm = \bigcap_{w\in W_J} w\overline{\cC^\pm_r} \subset \hft_r^\pm.\]
We identify ${\rm Simp}_J^\pm$ with a simplex of dimension $\ell-1-|J|$.

For example, when $J = \varnothing$, ${\rm Simp}_J^\pm = \overline{\cC^\pm_r}$ is an $\ell-1$ simplex, and when $J = S\backslash \{i\}$, ${\rm Simp}_J^\pm$ is a $0$-simplex.

\begin{Theorem}[simplicial structure on $\hfH_r$]\label{theorem:simplicial_structure}
For each $r<0$, with the notations above, we have a~Coxeter complex on $\big(k\hft k^{-1}\big)_r^\pm$ for each $k\in K$ and their union over all $k\in K$ forms a simplicial complex in $\hfH_r^\pm$.
\end{Theorem}

\begin{proof}It is straightforward to check that $\cS^\pm = \big\{{\rm Simp}_J^\pm\,|\, J\subsetneq S\big\}$ is a Coxeter complex~\cite{AbramenkoBrown08} for $(W,S)$. That is, $\cS^\pm$ admits a $W$-action which is simply transitively on simplices of maximal dimension. Thus for any two simplices of maximal dimension, there is a chain of maximal-dimensional simplices such that two consecutive ones share a common face of codimension~1. Since our compact real Kac--Moody algebra $\fk$ has rank $\ell$, all of its
Cartan subalgebras, $k\ft k^{-1}$, for $k\in K$, have dimension $\ell$ and each surface $\big(k\ft k^{-1}\big)_r$ of radius $r < 0$ has dimension $\ell-1$. We have now defined a Coxeter complex on $\big(k\hft k^{-1}\big)_r^\pm$ for each $k\in K$.

The union of all these Coxeter complexes forms a simplicial complex in $\hfH_r^\pm$ as follows. We must only check that the Weyl group tessellations on different Cartan subalgebras fit together in a well-defined way. Cartan subalgebras intersect exactly in hyperplanes fixed by Weyl group elements (as in Section~\ref{section:localstructure}), hence they intersect in singular elements. Thus each simplex of
maximal dimension lies in exactly one Cartan subalgebra. Thus the simplicial complexes in different Cartan subalgebras fit together, leading to a simplicial complex in~$\hfH_r^\pm$.

Each simplex of dimension $\ell-2$ lies in the intersection of a fixed point hyperplane $L_{k\ft k^{-1},i}$ with the surface $\big(k\ft k^{-1}\big)_r$ of radius $r < 0$ as calculated in Section~\ref{section:localstructure}, while simplices of dimension $\ell-n$ lie in the intersection of $n-1$ such hyperplanes with the surface. The stabilizer in~$K$ of~$\ft$ is the Weyl group~$W$ since $W=N/T$ is a quotient of the normalizer of~$T$. Hence simplicies of lower dimension are fixed by nontrivial subgroups of~$W$.
\end{proof}

Before we can state and prove the main embedding theorem in the next section, we must define the distance and codistance functions on $\fH_r$. For any two chambers $k_1 \cC_r^\pm k_1^{-1}$ and $k_2 \cC_r^\pm k_2^{-1}$ in $\hfH_{r}^\pm$, define the $W$-valued distance function
\[\delta^\pm\big(k_1 \cC_r^\pm k_1^{-1}, k_2 \cC_r^\pm k_2^{-1}\big) = w \in W\qquad\hbox{when}\qquad k_1^{-1} k_2\in B^\pm w B^\pm\]
and the codistance function
\[\delta^*\big(k_1 \cC_r^\pm k_1^{-1}, k_2 \cC_r^\mp k_2^{-1}\big) = w \in W\qquad\hbox{when}\qquad k_1^{-1} k_2\in B^\pm w B^\mp.\]
Note that $\big(\hfH_r^\pm,\delta^\pm\big)$ is a Tits building with an apartment system determined by $\delta^\pm$, and \linebreak $\big(\hfH_r^+,\hfH_r^-,\delta^*\big)$ is a twin building. We choose $\cC_r^\pm$ to be a fundamental chamber in $\ft_r^\pm$ and then $\{w\cC_r^\pm\,|\, w\in W\}$ is a fundamental apartment in $\hfH_r^\pm$.

\section{The main embedding theorem}

The main result of this paper, given in this section, is an embedding of the twin building $\cB=(\cB^+, \cB^-,\delta^*)$ of the Kac--Moody
group $G$ with twin $BN$-pair $(B^+,B^-,N)$ into the compact real form $\k$. This is a bijective simplicial map from $\cB^\pm$ onto the simplicial complex defined in $\hfH_{r}^\pm$ for each $r < 0$ defined in Section~\ref{section:simplicial structure}, such that the
$W$-valued distance and codistance functions are respected, and the map is $K$-equivariant, but not $G$-equivariant.

\begin{Theorem}\label{theorem:mainembedding theorem}Let $A$ be a symmetrizable hyperbolic generalized Cartan matrix, $\g$ its complex Kac--Moody algebra and $G$ its complex Kac--Moody group with twin $BN$-pair $(B^+,B^-,N)$, $W$-valued distance functions, $\delta^\pm$, $W$-valued codistance function, $\delta^*$, and let $r < 0$. Let $K$ be the compact real form of $G$ and let $\k$ be its Lie algebra. Let $\cB=(\cB^+, \cB^-,\delta^*)$ be the geometric realization $($Section~{\rm \ref{subsection:BNpair_Tits_building})} of the twin building of $G$ over $\C$, with Tits buildings, $(\cB^\pm,\delta^\pm)$ and codistance function $\delta^*$ in $\cB$. We use the same notations, $\delta^\pm$ and $\delta^*$, for distance functions between chambers in $\fH_r^\pm$ and codistance function in $\fH_r$.
\begin{enumerate}\itemsep=0pt
\item[$1.$] There is a $K$-equivariant simplicial map
$\Psi_r\colon \cB\hookrightarrow \hfH_r\subset \hfk$, that is, the following diagram commutes:
 \begin{displaymath}
 \begin{xy}
 \xymatrix{
 \cB \ar[rr]^{K} \ar[d]^{\Psi_r} & & \cB \ar[d]^{\Psi_r} \\
 \hfH_r \ar[rr]^{{\rm Ad}_K} & & \hfH_r
 }
 \end{xy}
 \end{displaymath}
where ${\rm Ad}_K$ denotes the adjoint action and $K\subset G$ acts on $\cB$ by left multiplication.
\item[$2.$] When $A$ is strictly hyperbolic the $K$-equivariant restrictions $\Psi_r^\pm\colon \cB^\pm\hookrightarrow \fH_r^\pm$ are bijective, otherwise, $\Psi_r^\pm\colon \cB^\pm\hookrightarrow \hfH_r^\pm$ are injective. For any chambers $C$, $D$ in $\cB^\pm$, we have $\delta^\pm(C,D) = \delta^\pm\big(\Psi_r^\pm(C), \Psi_r^\pm(D)\big)$, so $\Psi_r^\pm$ respects the distance functions
in $\cB^\pm$ and $\fH_r^\pm$.
\item[$3.$] There is a $W$-equivariant embedding of the fundamental apartment of $\cB^\pm$ into
$\hft^\pm_r\subset\hft_r$.
\item[$4.$] The map $\Psi_r$ is a twin building isomorphism. In particular, for any chamber $C$ in
$\cB^\pm$ and any chamber $D$ in $\cB^\mp$, we have $\delta^*(C,D) = \delta^*\big(\Psi_r^\pm(C), \Psi_r^\mp(D)\big)$, so $\Psi_r$ respects the codistance functions in $\cB$ and $\hfH_r$.
\end{enumerate}
\end{Theorem}

\begin{Remark} The image $\Im(\Psi_r)$ is contained in the interior of the lightcone if $A$ is strictly hyperbolic. Examples of this type are the rank $2$ hyperbolic Kac--Moody algebras (see Section~\ref{section:Lightcone construction of the Tits building in rank 2 hyperbolic type}). Otherwise it contains points on the boundary of the lightcone. Examples of this type are the algebras $\cF$, $\cI$ and $E_{10}$.
\end{Remark}

\begin{proof}To construct the embedding of the twin building we use the description given in equation~\eqref{equation:buildingk-invariance}
\[\cB^\pm=\big(K/T\times \Delta^{\pm}\big)/{\sim}. \]
Our construction is in three steps. We first define the embedding of the fundamental chamber. Then we make use of the $K$-action on the building and the adjoint action of $K$ on the Lie algebra. In a third step we establish the properties claimed in the theorem.

{\bf Part 1:} We choose the fundamental chambers in the buildings $\cB^\pm$ to be $c^\pm =\big(1 T, \Delta_\varnothing^\pm\big)\in \cB^\pm$. The $\ell$ panels of the fundamental chamber correspond to the simplices $\big(1 T, \Delta_i^\pm\big)$ for $1\leq i \leq \ell$, its vertices correspond to the simplices $(1 T, \Delta_{[i]}^\pm)$. Recall that $[i] = I\backslash\{i\}$. Two pairs $\big(f T, \Delta_i^\pm\big)$ and $\big(g T, \Delta_i^\pm\big)$ describe the same simplex if any only if $fK_i=gK_i$ where $K_i=K\cap P_i=K\cap (B\sqcup Bw_iB)$; similarly for any subset $J\subsetneq I$ two $J$-cells $\big(f T, \Delta_{J}^\pm\big)$ and $\big(f T, \Delta_{J}^\pm\big)$ are equivalent if and only $f K_J=g K_J$ for $K_J=K\cap P_J$.

For fixed $r<0$ we have chosen fundamental chambers $\cC_r^\pm$ such that
\begin{gather}\label{equation:opposite_fundamental_domain}
\cC_r^- = - \cC_r^+.
\end{gather}
 If $A$ is strictly hyperbolic the fundamental domain is contained in the interior of $\ft_r$; hence its closure is contained in $\ft_r$. On the other hand if $A$ is not strictly hyperbolic then it is unbounded and its closure is contained in $\hft_r$ but not in $\ft_r$.
 We identify the boundary hyperplanes $L_{\ft,i}$ of the fundamental chamber with the generators $w_i$, $1\leq i\leq \ell$, of the Weyl group $W$. We also identify intersections of hyperplanes $L_{\ft,i_1}\cap \dots \cap L_{\ft,i_n}$ for distinct indices $\{i_1, \dots, i_n\}\subset I$ with the subset of generators $\{w_{i_1}, \dots, w_{i_n}\}$. Vertices correspond to the intersection of $(\ell-1)$ hyperplanes and hence to subsets $[i] = I\backslash\{i\}$. For a subset $J\subsetneq I$ we define the boundary components \[\cC_J^\pm=\cC_r^\pm\cap\bigcap_{j\in J} L_{\ft, j}.\]

We define $\Psi_r\big(1T,\Delta_{[i]}^\pm\big)= \cC_{[i]}^\pm$ and extend this map to $\Delta^\pm$ by mapping a point $x$ in the geometric realization of the simplex $\big(1 T,\Delta_J^\pm\big)$ with normalized barycentric coordinates $x=[\lambda_1:\dots :\lambda_\ell]$ to the point with the same normalized hyperbolic barycentric coordinates.

Recall the definition of normalized hyperbolic barycentric coordinates: Let $\Delta$ be a simplex in $n$-dimensional hyperbolic space with vertices $(v_0, \dots, v_n)$. Some vertices may be on the boundary of hyperbolic space. We denote by $V= V(\Delta)$ the volume of $\Delta$. For any point $p\in \Delta$ we can define $n+1$ simplices $\Delta_{[i]}$ for $0\leq i\leq n$, spanned by the $(n+1)$-tuple of vertices $(v_0, \dots, p,\dots, v_n)$, where the vertex $v_i$ has been replaced by~$p$. Let $V_i=V(\Delta_{[i]})$ denote the volume of $\Delta_{[i]}$. Then the normalized hyperbolic barycentric coordinates of $p$ are given by
\[\left[\frac{V_0}{V}: \dots :\frac{V_{n}}{V}\right].\]

This yields a simplicial map $\Psi_r\colon \big(1 T, \Delta^\pm\big)\longrightarrow \cC_r^\pm$.

{\bf Part 2:} Let $x^\pm\in \Delta^\pm$. We extend the map defined in Part $1$ to a map \[\Psi_r\colon \ \cB\hookrightarrow \hfH_r\] by defining
\[\Psi_r\big(k T, x^\pm\big)={\rm Ad}_k \big[\Psi_r\big(1 T, x^\pm\big)\big].\]

We have to check that $\Psi_r$ is well-defined. Assume we have two equivalent elements $\big(k_1 T, x_1^\pm\big)\allowbreak \sim \big(k_2 T, x_2^{\pm}\big)$, hence $x_1^\pm=x_2^\pm$ and assuming $x_1^\pm\in \Delta^\pm_J$ for some subset $J\subsetneq I $, we have $k_1 K_J= k_2 K_J$. Then there is some $l\in K_J$ such that $k_1=k_2 l$ and we have
\[\Psi_r\big(k_1 T, x_1^\pm\big)={\rm Ad}_{k_1}\big[\Psi_r\big(1 T, x_1^\pm\big)\big] = {\rm Ad}_{k_2 l} \big[\Psi_r\big(1 T, x_1^\pm\big)\big].\]
Since $x_1^\pm\in \Delta_J^\pm$ we have \[\Psi_r\big(1T, x_1^\pm\big)\in \bigcap_{j\in J} L_{\ft,j}.\]
Hence ${\rm Ad}_l \Psi_r\big(1 T, x_1^\pm\big)=\Psi_r\big(1 T, x_1^\pm\big)$. Thus, from ${\rm Ad}_{k_2 l} = {\rm Ad}_{k_2}{\rm Ad}_{l}$ we get
\[{\rm Ad}_{k_2 l}\big[\Psi_r\big(1 T, x_1^\pm\big)\big]={\rm Ad}_{k_2}\big[\Psi_r\big(1 T, x_2^\pm\big)\big]=\Psi_r\big(k_2 T, x_2^{\pm}\big).\]

Note that $\Psi_r$ maps the simplex $(k T, \Delta_J^\pm)$ for $J\subsetneq I$ spanned by vertices $\big(k T, \Delta^\pm_{[i]}\big)$,
$i\in J$, onto the simplex spanned by $\Psi_r\big(kT, \Delta^\pm_{[i]}\big)$, $i\in J$.
If two simplices $\big(k_1 T, \Delta^\pm_J\big)$ and $\big(k_2 T, \Delta^\pm_J\big)$ share a common face $\big(lT, \Delta^\pm_{L}\big)$,
for $J\subsetneq L$ in $\cB$ then by definition we have $k_1 T\subset {\rm Ad}_l K_L $ and $k_2 T\subset {\rm Ad}_l K_L $. But then
$\Psi_r\big(k_1 T, \Delta^\pm_L\big)=\Psi_r\big(lT, \Delta^\pm_L\big)= \Psi_r\big(k_2 T, \Delta^\pm_L\big)$ is the commonly shared
face of the simplices $\Psi_r\big(k_1 T, \Delta^\pm_J\big)$ and $\Psi_r\big(k_2 T, \Delta^\pm_J\big)$. Hence $\Psi_r$ preserves the simplicial structure of $\cB$ and is thus a simplicial complex map.

{\bf Part 3:} We continue to check that $\Psi_r$ satisfies the properties stated in the theorem.

{\em Proof of~$(1)$:} Recall the left $K$-action on the building $\cB^\pm$:
\begin{displaymath}
K\colon \ \cB^\pm \longrightarrow \cB^\pm,\qquad k\cdot \big(fT, \Delta^\pm\big)\mapsto \big(kfT, \Delta^\pm\big).
\end{displaymath}

We need to verify that for $k_1,k_2\in K$ and $J\subsetneq I$
\[\Psi_r\big( k_1 \cdot\big(k_2T, \Delta^\pm_J\big)\big)={\rm Ad}_{k_1} \Psi_r\big(k_2 T, \Delta^\pm_J\big).\]
We have
\begin{align*}
\Psi_r (k_1\cdot \big(k_2 T, \Delta^\pm_J)\big)&=\Psi_r\big(k_1k_2T, \Delta^\pm_J\big)={\rm Ad}_{k_1k_2} \Psi_r\big(1 T, \Delta^\pm_J\big) \\
&={\rm Ad}_{k_1} \big({\rm Ad}_{k_2} \Psi_r\big(1 T,\Delta^\pm_J\big)\big)= {\rm Ad}_{k_1}\big(\Psi_r\big(k_2T, \Delta^\pm_J\big)\big).
\end{align*}

{\it Proof of~$(2)$:}
Assume two simplices $\big(k_1T, \Delta_J^\pm\big)$ and $\big(k_2T, \Delta_J^\pm\big)$ satisfy
\[\Psi_r\big(k_1T, \Delta_J^\pm\big)=\Psi_r\big(k_2T, \Delta_J^\pm\big).\]
Then we have
\[{\rm Ad}_{k_1} \Psi_r\big(1T, \Delta_J^\pm\big)={\rm Ad}_{k_2} \Psi_r\big(1T, \Delta_J^\pm\big),\]
which implies ${\rm Ad}_{k_1^{-1}k_2}\Psi_r\big(1T, \Delta_J^\pm\big)=\Psi_r\big(1T, \Delta_J^\pm\big)$.

As $\Psi_r\big(1T, \Delta_J^\pm\big)=\cC^\pm_J$ its stabilizer in $K$ is $K_J$. Hence $k_1^{-1}k_2\in K_J$.
Thus $\big(k_1T, \Delta_J^\pm\big)\sim\big(k_2T, \Delta_J^\pm\big)$ which proves injectivity.

Suppose we have an arbitrary element $x^\pm\in\fH_r^\pm$. Then there is some group element $k\in K$ such that
${\rm Ad}_k x^\pm\in \ft^\pm_r$ so it is in the closure of some chamber, $w\cC_r^\pm$, uniquely labeled by an element $w\in W$.
The action of the corresponding element $\wtil\in\Wtil^{\rm ad}\leq K$ matches the action of $w$ on $\ft$. So we have
${\rm Ad}_{\wtil^{-1}} {\rm Ad}_k x^\pm\in \cC^\pm_r$. Therefore
$x^\pm\in {\rm Ad}_{k^{-1}} {\rm Ad}_{\wtil} \Psi_r\big(1T, \Delta^\pm\big)=\Psi_r\big(k^{-1} \wtil T, \Delta^\pm\big)$.
This shows that in the case when~$A$ is strictly hyperbolic, $\Psi_r\colon \cB^\pm\to \fH_r^\pm$ is surjective. Otherwise there can be elements
$x^\pm\in B^\infty\big(\fH_r^\pm\big)$ which are in ${\rm Im}(\Psi_r)$, but some may not. It would be interesting to understand precisely which points in $B^\infty\big(\fH_r^\pm\big)$ are in ${\rm Im}(\Psi_r)$.

Let $C = \big(k_1T,\Delta^\pm\big)$ and $D = \big(k_2T,\Delta^\pm\big)$ be chambers in $\cB^\pm$. By definition, $\delta^\pm(C,D) = w\in W$ when $k_1^{-1} k_2\in B^\pm w B^\pm$.
We also have by definition,
\begin{gather*}
\Psi_r^\pm(C) = {\rm Ad}_{k_1}\big[\Psi_r^\pm\big(1T,\Delta^\pm\big)\big] = k_1 \cC_r^\pm k_1^{-1}, \qquad\hbox{and}\\
 \Psi_r^\pm(D) = {\rm Ad}_{k_2}\big[\Psi_r^\pm\big(1T,\Delta^\pm\big)\big] = k_2 \cC_r^\pm k_2^{-1} \qquad\hbox{so}\\
 \delta^\pm\big(\Psi_r^\pm(C),\Psi_r^\pm(D)\big) = \delta^\pm\big(k_1 \cC_r^\pm k_1^{-1}, k_2 \cC_r^\pm k_2^{-1}\big) = w,
\end{gather*}
since $k_1^{-1} k_2\in B^\pm w B^\pm$. This completes the proof of~(2).

The proof of~(3) is clear.

{\it Proof of~$(4)$:}
Let $C = \big(k_1T,\Delta^\pm\big)$ be a chamber in $\cB^\pm$ and let $D = \big(k_2T,\Delta^\mp\big)$ be a chamber in $\cB^\mp$. By definition,
$\delta^*(C,D) = w\in W$ when $k_1^{-1} k_2\in B^\pm w B^\mp$.
We also have by definition,
\begin{gather*}
\Psi_r^\pm(C) = {\rm Ad}_{k_1}\big[\Psi_r^\pm\big(1T,\Delta^\pm\big)\big] = k_1 \cC_r^\pm k_1^{-1}, \qquad\hbox{and}\\
\Psi_r^\mp(D) = {\rm Ad}_{k_2}\big[\Psi_r^\mp\big(1T,\Delta^\mp\big)\big] = k_2 \cC_r^\mp k_2^{-1} \qquad\hbox{so}\\
\delta^*\big(\Psi_r^\pm(C),\Psi_r^\mp(D)\big) = \delta^*\big(k_1 \cC_r^\pm k_1^{-1}, k_2 \cC_r^\mp k_2^{-1}\big) = w,
\end{gather*}
since $k_1^{-1} k_2\in B^\pm w B^\mp$.
\end{proof}

\begin{Remark} Based on our understanding of the examples $\cF$ and $E_{10}$, we believe that the intersection of the nullcone with ${\rm Im}(\Psi_r)$ can be characterized as those rays on the nullcone which contain roots of $\fg$. Each such ray of null vectors corresponds to a copy of an affine Kac--Moody subalgebra inside~$\fg$. Each such ray should be conjugate under the action of the Weyl group $W$ to one ray in the closure of the fundamental chamber, so the number of $W$ orbits is the number of ideal points in that fundamental domain. Classifying those rays then becomes a~significant problem involving the arithmetic properties of $W$, essentially understanding its cusps as a modular group. Such an analysis was carried out for $E_{10}$ in \cite[Lemma~5.2]{KleinschmidtPalmkvistNicolai12}.
\end{Remark}

Thus for each $r<0$ we have established an embedding $\Psi_r$ of the twin building $\cB=\big(\cB^+, \cB^-,\delta^*\big)$ into $\fH_r^\pm$ which is inside the lightcone
\[\cL_{\k}=\{x\in \k\,|\, (x,x)\leq 0\}\]
of the compact real form $\fk$ of a strictly hyperbolic Kac--Moody algebra $\fg$, and into its lightlike closure, $\hfH_r^\pm$, otherwise.

\begin{Remark}Recall that $\cB^\pm$ is contractible since each apartment is contractible (see~\cite{Ronan89}). Since $\Psi_r$ for $r\leq 0$ is a simplicial complex isomorphism onto its image, the image $\Psi_r\big(\cB^\pm\big)$ is also contractible.
\end{Remark}

This observation yields the corollary

\begin{Corollary}The spaces $\hfH_r^{\pm}$ are contractible.
\end{Corollary}

\begin{Remark}The analogs of Theorem~\ref{theorem:mainembedding theorem} for finite-dimensional compact Lie groups and affine Kac--Moody groups are well-known (see~\cite{Eberlein96, Freyn10a, Heintze06, Mitchell88}). It can be generalized in these cases to $s$-representations and relates in this way the (twin) building to the isotropy representations of finite-dimensional Riemannian symmetric spaces and affine Kac--Moody symmetric spaces respectively. The existence of hyperbolic Kac--Moody symmetric spaces has only recently been investigated in~\cite{FreynHartnickHornKoehl17}.
\end{Remark}

We conclude this section with a discussion of Cartan involutions and the embedding of the twin building.
Recall from Section~\ref{section:Campact_real_forms_of_KM_algebras_and_groups} that the Cartan involution is the automorphism ${\omega}_0\colon$ $\g\longrightarrow \g$ determined by ${\omega}_0(e_i)=-f_i$, $\omega_0(f_i)=-e_i$ and ${\omega}_0(h_i)= -h_i$, and $\omega$ is the conjugate linear automorphism defined as the composition of $\omega_0$ with complex conjugation. This gives
\[\omega_0 \big(\ft^\pm\big)=\ft^\mp \qquad \text{and}\qquad -\omega\big(\ft^\pm\big)=\ft^\mp \,\]
as well as
\[\omega_0 \big(\ft^\pm_r\big)=\ft_r^\mp \qquad \text{and}\qquad -\omega\big(\ft^\pm_r\big)=\ft^\mp_r \,\]
for each $r\in\R$, so that
\[\omega_0 \big(\hft^\pm_r\big)=\hft_r^\mp \qquad \text{and}\qquad -\omega\big(\hft^\pm_r\big)=\hft^\mp_r .\]

\begin{Proposition}The following diagram is commutative:
\begin{displaymath}
 \begin{xy}
 \xymatrix{
 \cB \ar[rr]^{\omega} \ar[d]^{\Psi_r} & & \cB \ar[d]^{\Psi_r} \\
 \k \ar[rr]^{-\omega} & & \k
 }
 \end{xy}
\end{displaymath}
Furthermore, for any chamber $C$ in $\cB^\pm$ and $D$ in $\cB^\mp$ we have $\delta^*(C,D) = \delta^*(\omega(C),\omega(D))$.
\end{Proposition}

\begin{proof} By definition $-\omega|_\k=-{\rm Id}$. As before, we set $\cB^\pm=\big(K/T\times \Delta^\pm\big)/{\sim}$. The action of the Cartan involution $\omega$ on $\cB$ is defined by
\[\omega\big(kT, \Delta^\pm_J\big)=\big(kT, \Delta^{\mp}_J\big).\]

Let $(kT, \Delta^\pm_J)\in \cB^\pm$. Then we have
 \begin{align*}
 -\omega \big(\Psi_r\big(kT, \Delta^\pm_J\big)\big)& =-\omega \big({\rm Ad}_k \Psi_r\big(1T, \Delta^\pm_J\big)\big)= -{\rm Ad}_{\omega(k)} \omega\big(\Psi_r\big(1T, \Delta^\pm_J\big)\big)\\
 & =- {\rm Ad}_k \Psi_r\big(\big(1T,\Delta^\pm_J\big)\big).
\end{align*}
On the other hand
\[\Psi_r \big(\omega \big(kT, \Delta^\pm_J\big)\big)=\Psi_r \big(kT, \Delta^\mp_J\big)= {\rm Ad}_k \Psi_r\big(\big(1T, \Delta^\mp_J\big)\big) = - {\rm Ad}_k \Psi_r\big(\big(1T,\Delta^\pm_J\big)\big),\]
where the last equality comes from equation~(\ref{equation:opposite_fundamental_domain}). The statement about $\delta^*$ is clear from the definitions.
\end{proof}

\begin{Remark}The choice of $-\omega$ on $\fk$ comes from the identification of $\k$ with the $\p$-component of the Cartan decomposition of $\fg$ via the relation $\p=i\k$.
\end{Remark}

We note that while $\Psi_r$ is a simplicial complex isomorphism, the apartment system of $\cB$ has no direct geometric interpretation in~$\Psi_r(\cB)$. Our embedding identifies $\omega$-stable twin apartments in $\cB$ with Cartan subalgebras in $\k$. The other apartments are hidden because $\Psi_r$ is only $K$-equivariant, but not $G$-equivariant.

\section{Special results for the twin building in rank 2}\label{section:Lightcone construction of the Tits building in rank 2 hyperbolic type}

In this section assume that $\fg = \fg(A)$ is a rank $2$ hyperbolic Kac--Moody Lie algebra, so that
\[A =\begin{bmatrix} 2 & -a\\ -b & 2\end{bmatrix},\qquad ab > 4.\]
Rank $2$ hyperbolic Kac--Moody algebras were studied intensively by Lepowsky and Moody~\cite{LepowskyMoody79}, by
Feingold~\cite{Feingold80} for the ``Fibonacci hyperbolic'' ($a=b=3$), $\cF{\rm ib}$, by Kang and Melville~\cite{KangMelville95},
by Carbone, Kownacki, Murray and Srinivasan \cite{CKMS19} and by Andersen, Carbone and Penta \cite{ACP11}. Twin trees and their relationship to Kac--Moody groups of rank~2 are studied in~\cite{Tits96}.
In these rank~2 cases the hyperboloids $\ft_r$ for $r\neq 0$ are hyperbolas, and there is no topological distinction between ``one-sheeted'' for $r > 0$ and ``two-sheeted'' for $r < 0$, as there is in higher rank.
In this section we will provide a detailed construction of the twin building, $\cB= \big(\cB^+, \cB^-,\delta^*\big)$, whose simplicial structure is a pair of trees, each equipped with a~$W$-valued distance function, $W$-valued codistance function, and give the apartment system and the action of the complex group~$G_\C$ on each tree. We will also give an explicit description of the apartments in the $K$-orbit of the fundamental twin apartment, which will
provide more details of the embedding of the twin building into the compact real form in this rank 2 case.
In the next section we show how this allows the embedding of a spherical building at infinity.

The eigenvalues of $A$ are $\lambda_{\pm}=2\pm \sqrt{ab}$ so $\lambda_+ > 0$ and $\lambda_- < 0$. This means that
the signature of the bilinear form determined by $A$ on the split real Cartan subalgebra $\fh_\R$ is $(1,1)$. By equation~(\ref{equation:compactmetric}), the bilinear form on the compact real Cartan subalgebra $\ft=i\fh_\R$ also has signature $(1,1)$.

The Weyl group of $\fg$ is the infinite dihedral group
\[W=\big\langle w_1,w_2\,|\, w_1^2=1=w_2^2\big\rangle\cong {\Z}/2{\Z} \ast{\Z}/2{\Z} \cong{\Z}\rtimes \{\pm 1\}, \]
which acts on $\fh\subset\fg$ as well as on $\fh_\R\subset\fg_\R$ and on $\ft\subset \k$.
Denote the index $2$ maximal infinite cyclic subgroup of $W$ by
\[W^{\rm even} = \big\{(w_2 w_1)^m\, |\, m\in\Z\big\}\cong \Z\]
so that for $i = 1,2$,
\[W^{\rm odd} = \big\{(w_2 w_1)^m w_i\, |\, m\in\Z\big\}\]
is the other coset, and we have the relations
\[w_i (w_2 w_1)^m w_i^{-1} = (w_2 w_1)^{-m}\qquad\hbox{for} \quad m\in\Z, \quad i = 1,2.\]
With $t = w_2 w_1$ and $r = w_i$ this gives the presentation of $W$ as the infinite dihedral group
$W = \big\langle t, r\, |\, r^2 = 1, \, rtr^{-1} = t^{-1}\big\rangle$. It can be useful to index the elements of $W$ by the integers so that
$W^{\rm even}$ and $W^{\rm odd}$ correspond to even and odd integers, respectively (after making a definite choice for $i$, say $i=2$):
\[w(n) = \begin{cases} (w_2 w_1)^m& \hbox{if }n = 2m,\\ (w_2 w_1)^m w_2& \hbox{if }n = 2m+1.\end{cases}\]
Then we have $w_1 w(n) = w(-1-n)$ and $w_2 w(n) = w(1-n)$ for $n\in\Z$, and
\[w(n)^{-1} = \begin{cases} (w_2 w_1)^{-m}& \hbox{if }n = 2m,\\ w_2 (w_2 w_1)^{-m}& \hbox{if }n = 2m+1\end{cases}
= \begin{cases} w(-n)& \hbox{if }n = 2m,\\ w(n)& \hbox{if }n = 2m+1.\end{cases}\]
It is straightforward to check that for any $n,k\in\Z$, we have
\[w(n)w(k) = \begin{cases} w(n+k)& \hbox{if }n = 2m,\\ w(n-k)& \hbox{if }n = 2m+1.\end{cases}\]
We will also find it useful to similarly label certain elements of
$\Wtil$ by the integers in the same way:
\[\wtil(n) = \begin{cases} (\wtil_2 \wtil_1)^m& \hbox{if }n = 2m,\\ (\wtil_2 \wtil_1)^m \wtil_2& \hbox{if }n = 2m+1.\end{cases}\]

\begin{figure}[t]\centering
\includegraphics[scale=0.2]{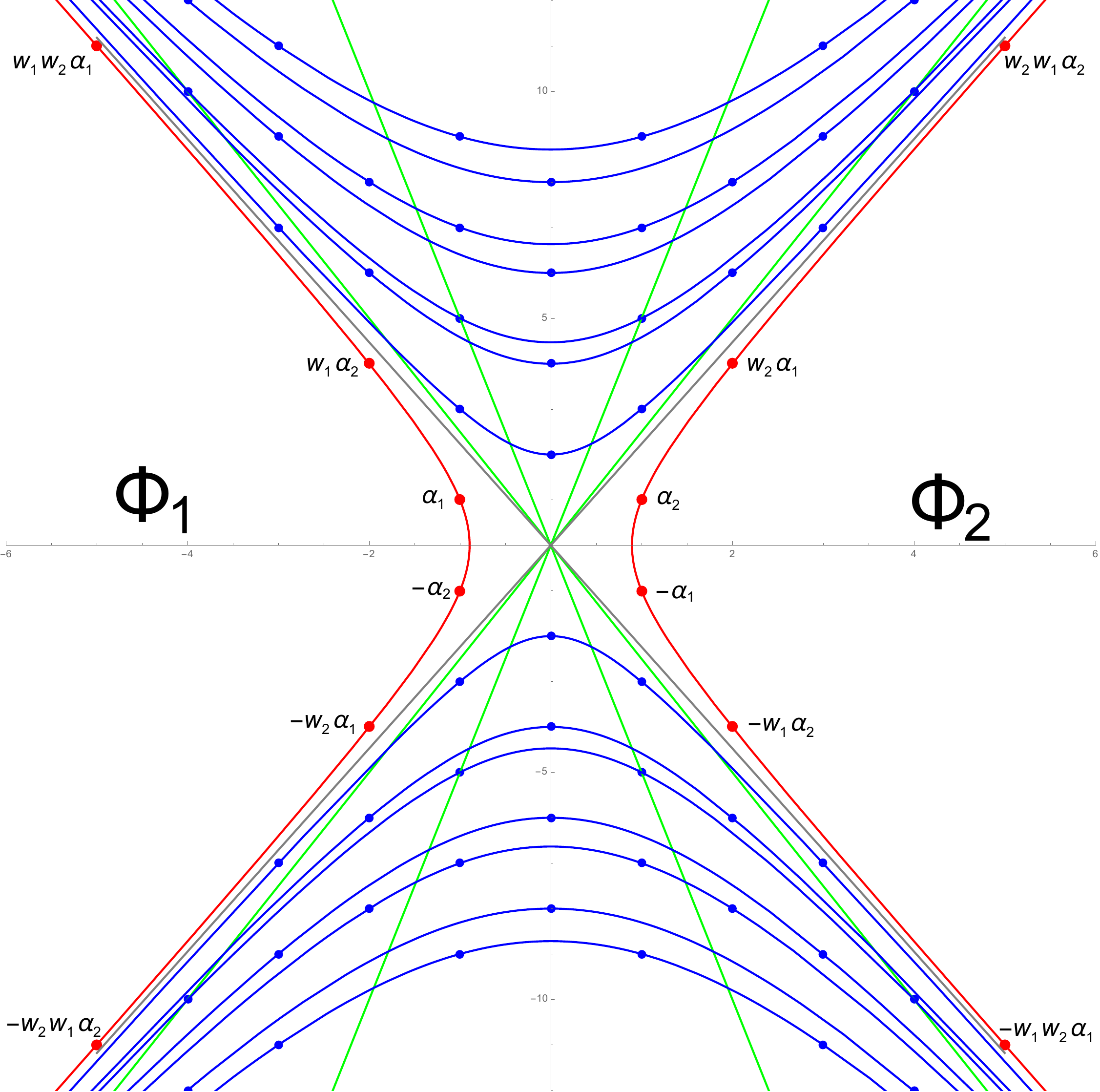}
\caption{The Fibonacci root system and non-standard partition of real roots
$\roots^{\rm re} = \roots_1\cup \roots_2$ where $\roots_1 = W^{\rm even}\{\a_1,-\a_2\}$ and
$\roots_2 = W^{\rm even}\{-\a_1,\a_2\}$.}\label{figure:Figure1}
\end{figure}

\begin{figure}[t]\centering
\includegraphics[scale=0.25]{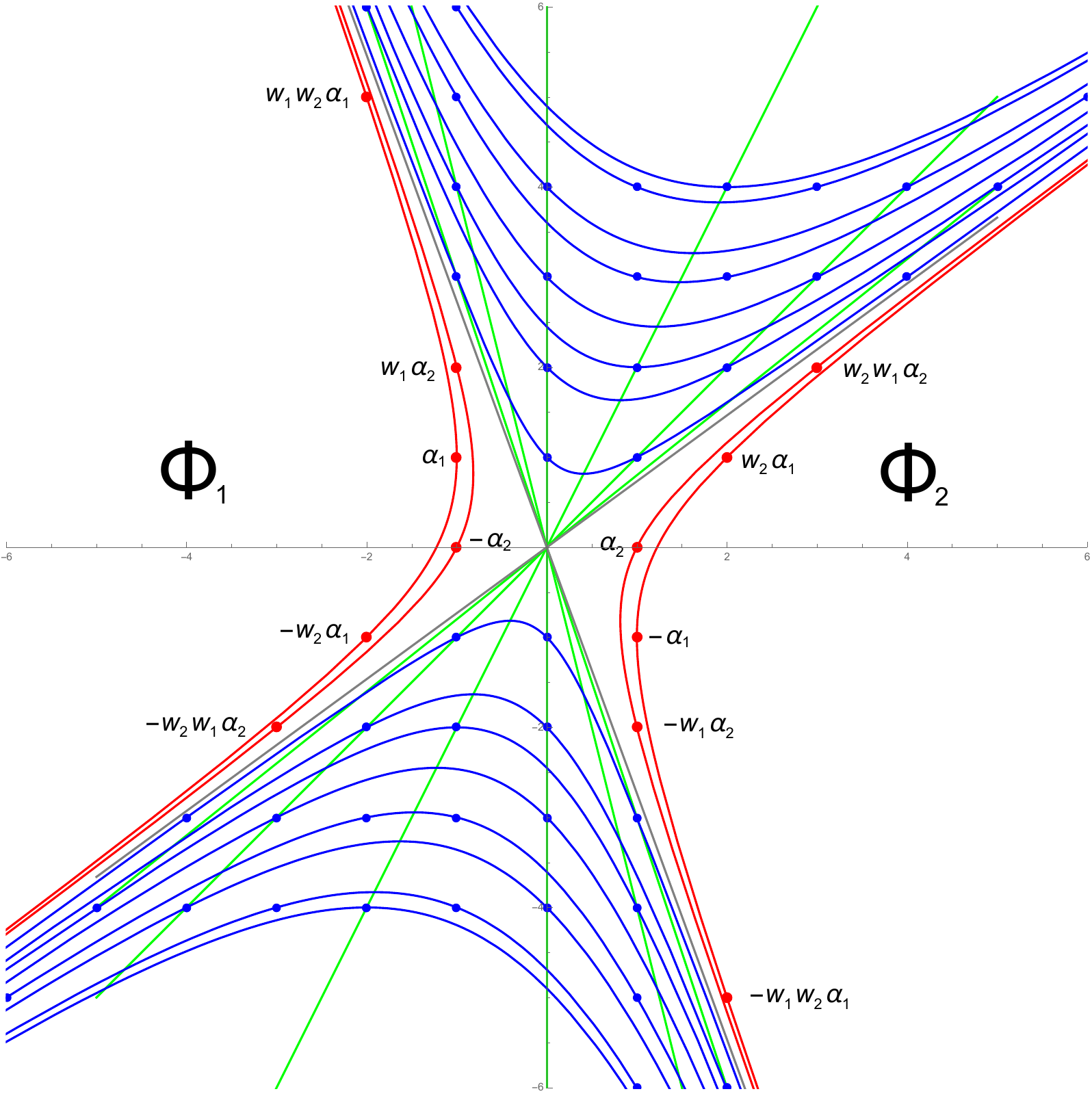}
\caption{The root system and non-standard partition of real roots for a rank 2 hyperbolic
with unequal root lengths, $a = 2$, $b = 3$.}\label{figure:Figure2}
\end{figure}

\begin{figure}[th!]\centering
\includegraphics{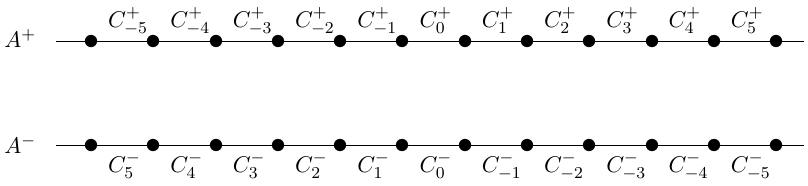}
\caption{Twin apartment $\big(A^+,A^-\big)$ in twin tree $\big(\cB^+,\cB^-\big)$ with chambers
$C^\pm_n = w(n) C^\pm_0$ and vertices $v^\pm_{n+\frac{1}{2}} = C^\pm_n \cap C^\pm_{n+1}$.}\label{figure:Figure3}
\end{figure}

We define the non-standard partition of the real roots of $\fg$, $\roots^{\rm re} = \roots_1\cup \roots_2$ where
\[\roots_1 = W^{\rm even}\a_1\cup W^{\rm odd} \a_2 = W^{\rm even}\{\a_1,-\a_2\}\]
and
\[\roots_2 = W^{\rm even}\a_2\cup W^{\rm odd} \a_1 = W^{\rm even}\{-\a_1,\a_2\} .\]
See Fig.~\ref{figure:Figure1} as an illustration for the ``Fibonacci'' rank~2 hyperbolic root system,
and see Fig.~\ref{figure:Figure2} for the rank 2 hyperbolic
root system coming from the Cartan matrix $A$ given at the beginning of
Section~\ref{section:Lightcone construction of the Tits building in rank 2 hyperbolic type}, with $a = 2$ and $b = 3$, which has
simple roots of different lengths. The real roots are on the red hyperbolas, and the non-standard partition is according to
whether a root is on a left branch or a right branch.
Imaginary roots of $\cF{\rm ib}$
are the blue dots inside the Tits cone (gray line asymptotes). The two green lines closest to the vertical $y$-axis are the lines fixed by $w_1$ ($\a_1(h) = 0$) and fixed by $w_2$ ($\a_2(h) = 0$), so the fundamental domains $C^+$ and $C^-$ (previously denoted by $\cC^\pm$) above and below the $x$-axis, respectively, are the wedges bounded by those green lines. Then $\bigcup\limits_{w\in W} w\cdot C^\pm$ fills up each component of the Tits cone. In contrast, the action of $W$ on the real roots has two orbits, and there are elements of each orbit on both left and right branches. The intersection of each timelike branch $\ft_r$, $r<0$, with the tessellation of the Tits cone gives the tessellation of each timelike branch into intervals, so we can map
a pair of timelike branches,~$\ft_r$, onto a pair of real lines,~$A^\pm$,
each tessellated by unit intervals centered on the integers. Then the fundamental domain (chamber) in $A^\pm$ is the interval $C^\pm_0 = \big[\frac{-1}{2}, \frac{1}{2}\big]$ centered on~$0$, and we label the chamber centered on $n\in\Z$ by
$C^\pm_n = w(n)C^\pm_0 = \big[n-\frac{1}{2}, n+\frac{1}{2}\big]$ as shown in Fig.~\ref{figure:Figure3}.
We may label the vertices $v^\pm_{n+\frac{1}{2}} = C^\pm_n \cap C^\pm_{n+1}$ (chamber walls) by elements in $\Z+\frac{1}{2}$.
The action of the generators of $W$ on the real lines $A^\pm$
are easily given by the formulas
\[w_1(x) = -1-x \qquad\hbox{and}\qquad w_2(x) = 1-x\]
so their action on the chambers and vertices of $\Sigma^\pm$ is given by
\begin{gather*} w_1 C^\pm_n = C^\pm_{(-1-n)},\qquad w_2 C^\pm_n = C^\pm_{(1-n)} ,\qquad
w_1 v^\pm_{n+\frac{1}{2}} = v^\pm_{-n-\frac{3}{2}},\qquad w_2 v^\pm_{n+\frac{1}{2}} = v^\pm_{-n+\frac{1}{2}}\end{gather*}
for any $n\in\Z$, so that $(w_2 w_1)^m C^\pm_n = C^\pm_{(n+2m)}$ and $(w_2 w_1)^m v^\pm_{n+\frac{1}{2}} =
v^\pm_{n+2m+\frac{1}{2}}$ for any $m,n\in\Z$.
We use this notation to describe a (left or right) {\it ray} of chambers starting with $C^\pm_n$ and including all boundary vertices:
\[L^\pm_{\rm ray}(n) = \big\{C^\pm_m\, |\, m\leq n\big\} \qquad\hbox{and}\qquad R^\pm_{\rm ray}(n) = \big\{C^\pm_m\, |\, m\geq n\big\} .\]
The action on rays is
\[w_1 L^\pm_{\rm ray}(n) = R^\pm_{\rm ray}(-1-n) \qquad\hbox{and}\qquad w_2 L^\pm_{\rm ray}(n) = R^\pm_{\rm ray}(1-n) \]
so for each $m,n\in\Z$ we have
\[(w_2 w_1)^m L^\pm_{\rm ray}(n) = L^\pm_{\rm ray}(n+2m) \qquad\hbox{and}\qquad (w_2 w_1)^m R^\pm_{\rm ray}(n) = R^\pm_{\rm ray}(n+2m) .\]
Note that left and right rays in $A^+$ are consistent with the usual orientation of the real line, but are reversed in
$A^-$, so we should think of viewing $A^-$ from above, facing down.

Let ${\rm Cham}^\pm = \big\{C^\pm_m\,|\, m\in\Z\big\}$ be the set of chambers of $A^\pm$, and define a $W$-distance function
$\delta^\pm\colon {\rm Cham}^\pm \times {\rm Cham}^\pm \to W$ by $\delta^\pm\big(w C^\pm,w' C^\pm\big) = w^{-1} w'$, so
$\delta^\pm\big(w' C^\pm,w C^\pm\big) = (w')^{-1} w = (w^{-1} w')^{-1}$ and
\[\delta^\pm\big(C^\pm_m,C^\pm_n\big) = w(m)^{-1} w(n)
= \begin{cases} w(n-m)& \hbox{if $m$ is even},\\ w(m-n)& \hbox{if $m$ is odd}.\end{cases}\]
Now define a $W$-codistance function
$\delta^*\colon {\rm Cham}^\pm \times {\rm Cham}^\mp \to W$ by $\delta^*\big(w C^\pm,w' C^\mp\big) = w^{-1} w'$, so
\[\delta^*\big(C^\pm_m,C^\mp_n\big) = w(m)^{-1} w(n)
= \begin{cases} w(n-m)& \hbox{if $m$ is even},\\ w(m-n)& \hbox{if $m$ is odd}.\end{cases}\]
Recall that a codistance function gives an opposition relation by $C^\pm_m\ {\rm op}\ C^\mp_n$ when $\delta^*\big(C^\pm_m,C^\mp_n\big) = 1$, but this is true iff $w(m) = w(n)$ which means $m = n$. Thus, each chamber $C^\pm_m$ in $A^\pm$ is opposite to precisely one chamber, $C^\mp_m$ in $A^\mp$, making $A = \big(A^+,A^-\big)$ a twin apartment in $\cB$.

Let ${\rm Vert}^\pm = \big\{v^\pm_{m-\frac{1}{2}}\,|\, m\in\Z\big\}$ be the set of vertices of $A^\pm$, and define $\N$-valued distance
functions and codistance function
\begin{gather*}
d^\pm\colon \ {\rm Vert}^\pm \times {\rm Vert}^\pm \to \N\qquad \hbox{by}\qquad d^\pm\big(v^\pm_{m-\frac{1}{2}},v^\pm_{n-\frac{1}{2}}\big) = |m-n|,\\
d^*\colon \ {\rm Vert}^\pm \times {\rm Vert}^\mp \to \N\qquad \hbox{by}\qquad d^\pm\big(v^\pm_{m-\frac{1}{2}},v^\mp_{n-\frac{1}{2}}\big) = |m-n|.\end{gather*}
Note that these functions are consistent with the $W$-valued functions if we define $|w(n)| = |n|$ so that
$\big|w(m)^{-1} w(n)\big| = |w(\pm(m-n))| = |m-n|$.

For $i = 1,2$ recall that
\[L_{\ft,i}= \{x\in \ft\, |\, w_i(x) = x\} = \{x\in \ft\, |\, \a_i(x) = 0\}\]
are the lines in $\ft$ fixed by the simple reflections $w_i$, respectively. Then
\[L_{\ft,1} = \{t(a z_1 + 2 z_2)\in \ft\, |\, t\in\R\}\qquad\hbox{and}\qquad
L_{\ft,2} = \{t(2 z_1 + b z_2)\in \ft\, |\, t\in\R\}.\]

For $i=1,2$ the line $L_{\ft,i}$ is the intersection of the family of Cartan subalgebras,
\[\big\{\exp({\rm ad}_{sx_i+ty_i})\ft = \ft^i(s,t)\, |\, s,t\in\R\big\}\]
parametrized by the $2$-sphere with antipodes identified, the real projective space $P_2(\R)$ (Corollary~\ref{corollary:hemispherical_parametrization}).

\begin{figure}[th]\centering
\includegraphics{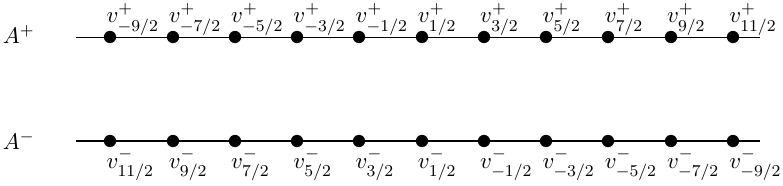}

\caption{Twin apartment $\big(A^+,A^-\big)$ in twin tree $\big(\cB^+,\cB^-\big)$ with vertices
$v^\pm_{m-\frac{1}{2}}$ with the co-distance function $d^*\big(v^+_{m-\frac{1}{2}},v^-_{n-\frac{1}{2}}\big) = |m-n|$.}\label{figure:Figure4}
\end{figure}

We also have the corresponding statements for the split real form.
For $i = 1,2$ recall that in the split real Cartan, $\fh_\R$, we have
\[L_{\fh_\R,i}= \{x\in \fh_\R\, |\, w_i(x) = x\} = \{x\in \fh_\R\, |\, \a_i(x) = 0\}\]
are the lines in $\fh_\R$ fixed by the simple reflections $w_i$, respectively. Then{\samepage
\[L_{\fh_\R,1} = \{t(a h_1 + 2 h_2)\in \fh_\R\, |\, t\in\R\}\qquad\hbox{and}\qquad
L_{\fh_\R,2} = \{t(2 h_1 + b h_2)\in \fh_\R\, |\, t\in\R\}.\]
In Figs.~\ref{figure:Figure1} and~\ref{figure:Figure2} those fixed lines are the inner green lines.}

For $i=1,2$ the line $L_{\fh_\R,i}$ is the intersection of the family of Cartan subalgebras,
\[\big\{\exp({\rm ad}_{se_i+tf_i})\fh_\R = \fh_\R^i(s,t)\, |\, s,t\in\R\big\}\]
parametrized by pairs of antipodal points on a $1$-sheeted hyperboloid.

Let $\big(B^+,B^-,N\big)$ be a twin $BN$-pair for the complex Kac--Moody group $G = G_\C(A)$.
The standard parabolic subgroups $P^\pm_J$ for $J\subsetneq \{1,2\}$ are
\[P^\pm_{\varnothing} = B^\pm, \qquad
P^\pm_1= B^\pm\sqcup B^\pm w_1 B^\pm,\qquad \hbox{and} \qquad
P^\pm_2=B^\pm\sqcup B^\pm w_2 B^\pm.\]
As a simplicial complex, the twin building $\cB = \big(\cB^+,\cB^-,\delta^*\big)$ associated to $\big(B^+,B^-,N\big)$ is a pair of homogenous $\P^1(\C)$-trees. The vertices of $\cB$ are in bijection with the conjugates of $P^\pm_1$ and $P^\pm_2$ in $G$, while the set of edges (chambers) are in bijection with conjugates of $B^\pm$. We identify the sets of vertices with the disjoint union of cosets
\[V\big(\cB^\pm\big) = G/P^\pm_1\sqcup G/P^\pm_2.\]
The set of edges is given by
\[E\big(\cB^\pm\big) = G/B^\pm.\]
The group $G$ acts by left multiplication on cosets. There are natural projections
on cosets induced by the inclusion of $B^\pm$ in $P^\pm_1$ and $P^\pm_2$:
\[\pi_i\colon \ G/B^\pm\longrightarrow G/P^\pm_i,\qquad i=1,2.\]
If $v_i^\pm\in G/P^\pm_i$ is a vertex, and
${\rm St}\big(v_i^\pm\big)=\pi_i^{-1}(v_i^\pm)$ is the set of edges
with origin $v_i^\pm$, then we may index ${\rm St}\big(v_i^\pm\big)$ by $P^\pm_i/B^\pm\subseteq G/B^\pm$, $i=1,2$. It can be seen that $P^\pm_i/B^\pm \cong \P^1(\C) = \{\infty\}\cup\C$.
The $W$-valued distance and codistance functions defined on the chambers of a twin building can be converted into $\N$-valued functions by composing with the function $|w(n)| = |n|$, which coincides with the length function. So for any two chambers, $C^\pm = \big(g_1B^\pm,\Delta^\pm\big)$ and $D^\pm = \big(g_2B^\pm,\Delta^\pm\big)$ in $\cB^\pm$ we define
\[d^\pm\big(C^\pm, D^\pm\big) = \big|\delta^\pm(C,D)\big|\qquad \hbox{and}\qquad d^*\big(C^\pm, D^\mp\big) = |\delta^*(C,D)|.\]
These distance and codistance functions on chambers may be converted to functions on the vertices of the
twin tree. Ronan and Tits~\cite{RonanTits94} assume that the distance and codistance functions on vertices exist and obey their axioms. In this description some care must be taken because a~vertex can be in many chambers.

{\it Apartments} in $\cB^\pm$ are infinite lines, but not every infinite line in the tree is an apartment.
The fundamental apartment $A^\pm$ in $\cB^\pm$ is a union of chambers, $C^\pm_n = w(n)C^\pm$, $n\in\Z$, which are Weyl group translates of the fundamental chamber $C^\pm = \big(1B^\pm, \Delta_{\varnothing}^\pm\big)$, which is fixed under the action of $B^\pm\supset U^\pm$. The boundary vertices of that fundamental chamber are denoted
by $v_{i-\frac{3}{2}}^\pm= \big(1B^\pm, \Delta_{i}^\pm\big)$ for $i = 1,2$, corresponding to the maximal parabolic subgroups $P_i^\pm$. By definition, the vertex $v_{i-\frac{3}{2}}^\pm$ is stabilized by the generator $w_i\in W$.
We can see that chambers $C^\pm_m$ and $C^\pm_n$
are adjacent with a common vertex iff $|m-n|=1$ because $\big(w(m)B^\pm, \Delta_{i}^\pm\big)\sim \big(w(n)B^\pm, \Delta_{i}^\pm\big)$ when
$w(\pm(m-n)) = w(m)^{-1} w(n)\in P^\pm_i$ iff $|m-n|=1$ ($w(1) = w_2\in P^\pm_2$ and $w(-1) = w_1\in P^\pm_1$). Thus we see the
apartment $A^\pm$ as an infinite line:
\[\big\{\dots, w_1w_2w_1\cC^\pm, w_1w_2\cC^\pm, w_1\cC^\pm, \cC^\pm, w_2\cC^\pm, w_2w_1\cC^\pm,
 w_2w_1w_2\cC^\pm, \dots\big\}.\]
The apartment system $\cA^\pm$ in $\cB^\pm$ is the set of apartments obtained from fundamental apartment~$A^\pm$ by the action of~$G$.

We can describe the action of the real root groups $U_{\pm(w_1 w_2)^m\a_i}$, $i=1,2$, $m\in\Z$, on the fundamental apartments~$A^\pm$.

Using the notation ${\hat i} = 3-i$, we label the roots in $\roots_i$ by the integers as follows:
\[\roots_i(n) = \begin{cases} w(n) \a_i & \hbox{if }n = 2m, \\ w(n) \a_{\hat i} & \hbox{if }n = 2m+1\end{cases}
= \begin{cases} (w_2 w_1)^m \a_i& \hbox{if }n = 2m, \\ (w_2 w_1)^m w_2 \a_{\hat i}& \hbox{if }n = 2m+1,\end{cases}\]
so that the labels of roots in both branches are the integers in order, negative to positive going in $\roots_1$ from
top to bottom, but going in $\roots_2$ from bottom to top:
\[\roots_i(n) \in \roots^+ \quad \text{for} \quad \begin{cases} n\leq 0 & \hbox{if }i = 1, \\ n\geq 0 & \hbox{if }i=2. \end{cases}\]
We also have $w_1\roots_1(n) = \roots_2(1-n)$, $w_2\roots_1(n) = \roots_2(-1-n)$ and $-\roots_2(n) = \roots_1(n+1)$.

Here is another useful application of the labeling of the real roots by $\roots_i(n)$, $n\in\Z$, as shown in the paragraph above.
A consistent choice of real root vectors given by
\[e_{\roots_i(n)} = \begin{cases} e_{w(n) \a_i}& \hbox{if }n = 2m, \\ e_{w(n) \a_{\hat i}}& \hbox{if }n = 2m+1\end{cases}
= \begin{cases} \wtil(n) e_{\a_i} & \hbox{if }n = 2m, \\ \wtil(n) e_{\a_{\hat i}}& \hbox{if }n = 2m+1,\end{cases}\]
and then we would have $\wtil_j e_\a = e_{w_j \a}$ for all real roots $\a$ and $j=1,2$. $L^\pm_{\rm ray}(n)$ is fixed by $U_{\pm\roots_2(n)}$ and $R^\pm_{\rm ray}(n)$ is fixed by $U_{\pm\roots_1(n)}$ for all
$n\in\Z$, but $U_{\pm\roots_2(n)} R^\pm_{\rm ray}(n-1)$ is a distinct family of rays in $\cB^\pm$ indexed by $\C$ whose intersection with $L^\pm_{\rm ray}(n)$ is a unique vertex, and
$U_{\pm\roots_1(n)} L^\pm_{\rm ray}(n+1)$ is a distinct family of rays in $\cB^\pm$ indexed by $\C$ whose intersection with $R^\pm_{\rm ray}(n)$ is a unique vertex.

\begin{Proposition}\label{prop:root_group_actions_on_Fund_Apt}
The action of real root groups $U_{\pm\roots_i(k)}$, $k\in\Z$, on the chambers $\cC^\pm(n)$, in the fundamental apartments
$A^\pm$ satisfies the following.
\begin{enumerate}\itemsep=0pt
\item[$1.$] The chambers of $L^\pm_{\rm ray}(k) = \big\{\cC^\pm(n)\, |\, n\leq k\big\}$ are each fixed by $U_{\pm\roots_2(k)}$ but $U_{\pm\roots_2(k)} R^\pm_{\rm ray}(k+1)$ is a family of distinct rays in $\cB^\pm$ indexed by $\C$ whose intersection with $L^\pm_{\rm ray}(k)$ is a unique vertex.
\item[$2.$] The chambers of $R^\pm_{\rm ray}(k) = \big\{\cC^\pm(n)\, |\, n\geq k\big\}$ are each fixed by $U_{\pm\roots_1(k)}$ but $U_{\pm\roots_1(k)} L^\pm_{\rm ray}(k-1)$ is a family of distinct rays in $\cB^\pm$ indexed by $\C$ whose intersection with $R^\pm_{\rm ray}(k)$ is a unique vertex.
\end{enumerate}
\end{Proposition}
\begin{proof} (1) We know that $\cC^\pm(n) = w(n) \cC^\pm$ and that for real $\a\in\roots$, $U_\a \cC^\pm = \cC^\pm$
when $\a\in\roots^\pm$. So $U_\a \cC^\pm(n) = \cC^\pm(n)$ when $w(n)^{-1} U_\a w(n) \cC^\pm = \cC^\pm$, that is,
when $U_{w(n)^{-1} \a} \cC^\pm = \cC^\pm$. But that happens when $w(n)^{-1} \a \in\roots^\pm$. Recall that
\[w(n)^{-1} = \begin{cases} w(-n)& \hbox{if }n = 2r,\\ w(n) & \hbox{if }n = 2r+1\end{cases}\qquad\hbox{and}\qquad
w(n)w(k) = \begin{cases} w(n+k)& \hbox{if }n = 2r,\\ w(n-k)& \hbox{if }n = 2r+1,\end{cases}\]
so for fixed $k\in\Z$, any $n\leq k$, and
\[\a = \pm\roots_2(k) = \begin{cases} \pm w(k) \a_2& \hbox{if }k = 2m, \\ \pm w(k) \a_{1}& \hbox{if }k = 2m+1, \end{cases}\]
we have
\begin{align*}
w(n)^{-1} \a & = \pm w(n)^{-1} \roots_2(k) =
\pm \begin{cases} w(n)^{-1} w(k) \a_2& \hbox{if }k = 2m, \\
w(n)^{-1} w(k) \a_{1}& \hbox{if }k = 2m+1 \end{cases}\\
 & = \pm \begin{cases} w(-n) w(k) \a_2& \hbox{if }n = 2r,\ k = 2m, \\
w(-n) w(k) \a_{1}& \hbox{if }n = 2r,\ k = 2m+1, \\
w(n) w(k) \a_2& \hbox{if }n = 2r+1,\ k = 2m, \\
w(n) w(k) \a_{1}& \hbox{if }n = 2r+1,\ k = 2m+1
\end{cases}\\
& = \pm \begin{cases} w(-n+k) \a_2&\hbox{if }n = 2r,\ k = 2m, \\
w(-n+k) \a_{1}& \hbox{if }n = 2r,\ k = 2m+1, \\
w(n-k) \a_2& \hbox{if }n = 2r+1,\ k = 2m, \\
w(n-k) \a_{1}& \hbox{if }n = 2r+1,\ k = 2m+1
\end{cases}\\
& = \pm \begin{cases} \roots_2(k-n)& \hbox{if }n = 2r,\ k = 2m, \\
\roots_2(k-n)& \hbox{if }n = 2r,\ k = 2m+1, \\
\roots_1(n-k)& \hbox{if }n = 2r+1,\ k = 2m, \\
\roots_1(n-k)& \hbox{if }n = 2r+1,\ k = 2m+1
\end{cases} \in \pm \roots^+ = \roots^\pm,
\end{align*}
since
\[\roots_i(s) \in \roots^+ \quad \hbox{for} \quad \begin{cases} s\leq 0 & \hbox{if }i = 1, \\ s\geq 0 & \hbox{if }i=2. \end{cases}\]
For $\a = \pm\roots_2(k)$ and $t\in\C$ let $g(t) = \exp({\rm ad}_{te_\a})\in U_\a$ and consider the family of rays
\[U_{\a} R^\pm_{\rm ray}(k+1) = \big\{g(t) R^\pm_{\rm ray}(k+1)\, |\, t\in\C \big\}\]
with chambers $\big\{g(t) \cC^\pm(n)\, |\, n\geq k+1, t\in\C \big\}$. Two such rays are certainly distinct if their first chambers are distinct, so suppose that $g(t_1) \cC^\pm(k+1) = g(t_2) \cC^\pm(k+1)$, that is,
$g(t_1) w(k+1) \cC^\pm = g(t_2) w(k+1) \cC^\pm$. Let $w = w(k+1)$ and $g^w = w^{-1} g w$, so we have
$g(t_1)^w \cC^\pm = g(t_2)^w \cC^\pm$ which gives $\cC^\pm = (g(t_1)^w)^{-1} g(t_2)^w \cC^\pm =
g(-t_1)^w g(t_2)^w \cC^\pm = g(-t_1+t_2)^w \cC^\pm$. Therefore, $g(-t_1+t_2)^w \in B^\pm$.
But for any $t\in\C$, we have $g(t)^w\in U_{w^{-1}\a}$.
Since
\[w^{-1} = w(k+1)^{-1} = \begin{cases} w(-k-1)& \hbox{if }k\hbox{ is odd}, \\
w(k+1)& \hbox{if }k\hbox{ is even}, \end{cases}\]
we get
\begin{align*}
w^{-1} \a & = \pm \begin{cases} w(-k-1) w(k)\a_1& \hbox{if }k\hbox{ is odd}, \\
w(k+1)w(k)\a_2& \hbox{if }k\hbox{ is even} \end{cases}
= \pm \begin{cases} w(-1) \a_1& \hbox{if }k\hbox{ is odd}, \\
w(1) \a_2\quad\hbox{if }k\hbox{ is even} \end{cases}\\
& = \pm \begin{cases} w_1 \a_1& \hbox{if }k\hbox{ is odd}, \\
w_2 \a_2& \hbox{if }k\hbox{ is even} \end{cases}
= \pm \begin{cases} -\a_1& \hbox{if }k\hbox{ is odd}, \\
-\a_2& \hbox{if }k\hbox{ is even}, \end{cases}
\end{align*}
so $w^{-1} \a \in -\roots^\pm$ which means $U_{w^{-1}\a} \leq U^\mp$.
Since $g(-t_1+t_2)^w \in B^\pm\cap U^\mp = \{1\}$ we get $g(-t_1+t_2)^w = 1$ so $g(-t_1+t_2) = 1$ so
$t_1 = t_2$.We have shown that this family consists of distinct rays indexed by $\C$.

The proof of (2) is similar.
\end{proof}

\section{Embedding the spherical building at infinity in rank 2}\label{section:Embedding the spherical building at infinity}

Let $A$ be a rank~2 hyperbolic Cartan matrix as in
Section~\ref{section:Lightcone construction of the Tits building in rank 2 hyperbolic type}, but with the additional conditions
that $a>1$ and $b>1$ because otherwise the real root groups~$U_i$ defined below
(in Definition~\ref{definition:nonstandard_root_groups}) may not be abelian (see \cite{CKMS19,Morita88}).
Let $G = G_\C(A)$ be a~rank~2 hyperbolic Kac--Moody group with compact real form $K$. Let $\cB= \big(\cB^+, \cB^-,\delta^*\big)$ denote its twin building, associated with a~twin $BN$-pair $\big(B^+,B^-,N\big)$, whose simplicial structure is a~pair of trees,
and whose apartments are lines which are infinite in both directions.
We will now define the {\em spherical twin building at infinity}, $\cB_\infty = \big(\cB_\infty^+,\cB_\infty^-,\delta_\infty^*\big)$
even though we do not obtain a twin $BN$-pair at infinity.
Our results should be compared to those of Carbone and Garland~\cite{CG2003}, where they obtain a spherical
building at infinity, a $BN$-pair and Bruhat decomposition for complete rank $2$ Kac--Moody groups over finite fields.

A {\em ray} in a tree is a sequence of incident vertices and edges $(v_1,e_1,v_2,e_2,\dots)$ with an initial vertex, $v_1$, infinite in only one direction. Two rays are called equivalent if their intersection is a~ray. An equivalence class of rays is called an {\em end} of the tree, and we denote the end of a ray $X_{\rm ray}$ by~$[X_{\rm ray}]$. A {\em line} in a tree is a sequence of incident vertices and edges infinite in both directions, so it can be expressed as the union of two rays having a finite nonempty intersection, which can always be taken to be exactly one vertex. So each line has two ends, the classes of those two rays. For each apartment, $X = X_{\rm ray1}\cup X_{\rm ray2}$ in $\cB^\pm$ there is an apartment $[X] = [X_{\rm ray1}]\cup [X_{\rm ray2}]$ in $\cB_\infty^\pm$ whose two chambers $[X_{\rm ray1}]$ and $[X_{\rm ray2}]$ uniquely determine the line $X$ as follows.

Without loss of generality, we may assume that $X_{\rm ray1}\cap X_{\rm ray2} = \{X_0\}$ is exactly one vertex.
Suppose another apartment $Y = Y_{\rm ray1}\cup Y_{\rm ray2}$ with $Y_{\rm ray1}\cap Y_{\rm ray2} = \{Y_0\}$ also exactly one vertex, has
$[X_{\rm ray1}] = [Y_{\rm ray1}]$ and $[X_{\rm ray2}] = [Y_{\rm ray2}]$, so the intersections
$I_1 = X_{\rm ray1}\cap Y_{\rm ray1}$ and $I_2 = X_{\rm ray2}\cap Y_{\rm ray2}$ are both rays. The intersection
$I_1\cap I_2\subseteq X_{\rm ray1}\cap X_{\rm ray2}$ is either one vertex $\{X_0\}$ or empty, and
$I_1\cap I_2\subseteq Y_{\rm ray1}\cap Y_{\rm ray2}$ is either one vertex $\{Y_0\}$ or empty. So if $I_1\cap I_2$ is non-empty, then it is one vertex
and that vertex is $X_0 = Y_0$, so $X = Y$. Assume now that $I_1\cap I_2$ is empty and write $I_1 = (v_1,e_1,v_2,e_2,\dots)$ and
$I_2 = (w_1,f_1,w_2,f_2,\dots)$ so $v_1$ is a vertex in $X_{\rm ray1}\cap Y_{\rm ray1}$ and $w_1$ is a vertex in $X_{\rm ray2}\cap Y_{\rm ray2}$ but
the connected path in $X$ from $v_1$ to $w_1$ goes through $X_0$ while the connected path in $Y$ from $v_1$ to $w_1$ goes
through $Y_0$. If those two paths were distinct, there would be a loop in the tree, so they must be identical, giving $X = Y$.

Roughly speaking, a spherical building at infinity of a tree $\cB^\pm$ consists of a `sufficiently large' set of ends. The minimal spherical
building at infinity may be constructed as follows. For each sign $\pm$ let $A^\pm$ be the fundamental apartment of $\cB^\pm$ with
fundamental chamber $\cC^\pm$, and let $\cA^\pm = \big\{g\cdot A^\pm\, |\, g\in G\big\}$ be the set of all $G$-translates of $A^\pm$.
Then $\cA^\pm$ is the minimal apartment system of $\cB^\pm$ containing $A^\pm$.
(See~\cite[Section 4.5]{AbramenkoBrown08} for a discussion of complete apartment systems, and Section~11.8.4 for a discussion of incomplete apartment systems.)
A chamber in $\cB_\infty^\pm$ is an end of
an apartment in $\cA^\pm$, and ${\rm Cham}_\infty^\pm$ denotes the set of all chambers of $\cB_\infty^\pm$.
Of course, this building is highly degenerate. As a simplicial complex it consists of only $0$-simplices, and its apartments are
the subsets of two distinct points corresponding to the two ends of an apartment in $\cA^\pm$.
So in this rank~2 case, each apartment is the Coxeter complex consisting of those two points with order $2$ Coxeter group
$W_{\infty} = W/W^{\rm even} \cong \{\pm 1\} \cong \Z/2\Z$.
Since $W^{\rm even}$ acts as translations in each apartment of $\cB^\pm$, it stabilizes the ends, while $W^{\rm odd}$ switches them.

We need to define the distance $\delta_\infty^\pm$ and codistance $\delta_\infty^*$.
Any two distinct chambers $c^\pm = \big[X_{\rm ray1}^\pm\big]$ and $d^\pm = \big[X_{\rm ray2}^\pm\big]$ in ${\rm Cham}_\infty^\pm$ come from a~unique line
$X^\pm = X_{\rm ray1}^\pm\cup X_{\rm ray2}^\pm$ in $\cB^\pm$ written as a~union of two rays whose intersection is a~vertex $X_0^\pm$.
Let $C_1^\pm\in X_{\rm ray1}^\pm$ and $C_2^\pm\in X_{\rm ray2}^\pm$ be the initial chambers in these rays whose intersection is
$X_0^\pm$.
Then we define the $W_{\infty}$-valued distance function by
\[\delta_\infty^\pm\big(c^\pm, d^\pm\big) = \delta^\pm\big(C_1^\pm,C_2^\pm\big)W^{\rm even} = \begin{cases} -1& \hbox{if }c^\pm\neq d^\pm, \\
\hphantom{-}1& \hbox{if }c^\pm = d^\pm. \end{cases} \]

In order to use the codistance function $\delta^*$ in $\cB$ to get a codistance function $\delta_\infty^*$,
we need the concept of {\it twin ends}.

Let $x$, $y$, $z$ be chambers in $\cB^\pm$. Then the $\N$-valued distance function defined in
Section~\ref{section:Lightcone construction of the Tits building in rank 2 hyperbolic type}
satisfies the following properties.
\begin{enumerate}\itemsep=0pt
\item $d^\pm(x,y) = 0$ iff $x=y$.
\item If $d^\pm(y,z)=1$ we have $d^\pm(x,z) = d(x,y)^\pm \pm 1$ so $d^\pm(x,z) \neq d^\pm(x,y)$.
Also, there is at most one $z$ with $d^\pm(y,z) = 1$ and $d^\pm(x,z) = d^\pm(x,y) - 1$.
\item For any two distinct $x,y$, there exists $z$ with $d^\pm(y,z) = 1$ and $d^\pm(x,z) = d^\pm(x,y) + 1$.
\end{enumerate}
Let $x_\pm$ be a chamber in $\cB^\pm$ and let $y_\mp$ and $z_\mp$ be chambers in $\cB^\mp$. Then the $\N$-valued codistance
function defined in Section \ref{section:Lightcone construction of the Tits building in rank 2 hyperbolic type}
satisfies the following properties.

\begin{enumerate}\itemsep=0pt
\item If $d^\mp(y_\mp,z_\mp) = 1$ then $d^*(x_\pm,z_\mp) = d^*(x_\pm,y_\mp) \pm 1$.
\item Additionally, if $d^*(x_\pm,y_\mp)\neq 0$ then there exists a unique $z_\mp$ with $d^\mp(y_\mp,z_\mp) = 1$ and
$d^*(x_\pm,z_\mp) = d^*(x_\pm,y_\mp) + 1$.
\end{enumerate}

Let $A^\pm$ be an apartment in $\cB^\pm$, and let ${\rm Cham}^\pm\big(A^\pm\big)$ be the set of chambers of $A^\pm$. $\big(A^+,A^-\big)$ is a twin apartment in the twin tree $\big(\cB^+,\cB^-\big)$ iff for every $x_\pm\in {\rm Cham}^\pm(A^\pm)$ there exists a~unique $x_\mp\in {\rm Cham}^\mp\big(A^\mp\big)$ such that $d^*(x_\pm,x_\mp) = 0$.

\begin{figure}[t]\centering
\includegraphics{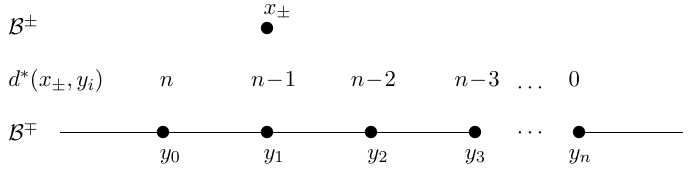}
\caption{Co-distances $d^*(x_\pm,y_i)=n-i$ between $x_\pm\in\cB^\pm$ and a path of adjacent vertices in $\cB^\mp$.} \label{figure:Figure5}
\end{figure}

Observations: Fix some $x_\pm\in {\rm Cham}^\pm\big(A^\pm\big)$ and let $(y_0,y_1,y_2,\dots,y_n)$ be a path of adjacent chambers in
$\cB^\mp$ with $d^*(x_\pm,y_1) = d^*(x_\pm,y_0) - 1$ then $d^*(x_\pm,y_i) = d^*(x_\pm,y_0) - i$ for $0\leq i\leq d^*(x_\pm,y_0)$.
(See Fig.~\ref{figure:Figure5}.) Consequently, for a twin apartment, $\big(A^+,A^-\big)$, if $x_\pm\in {\rm Cham}^\pm(A^\pm)$ and
$(y_0,y_1,y_2,\dots,y_n)$ is a path of adjacent chambers in $A^\mp$ with $d^*(x_\pm,y_1) = d^*(x_\pm,y_0) + 1$ for $d^*(x_\pm,y_0) > 0$ then $d^*(x_\pm,y_i) = d^*(x_\pm,y_0) + i$ for $0\leq i\leq n$.

For $x\in {\rm Cham}^\pm(A^\pm)$ and $y\in {\rm Cham}^\mp\big(A^\mp\big)$ with $d^*(x,y) = n > 0$, there is a {\it twin end} determined as follows. Using $x = x_\pm$ and $y = y_0$ as above, there is a~ray $r_\mp$ of adjacent chambers $(y_0,y_1,y_2,y_3,\dots)$ in ${\rm Cham}^\mp\big(A^\mp\big)$ with
$d^*(x,y_i) = n + i$ for $0\leq i$, and similarly there is also a~ray~$r_\pm$ of adjacent chambers $(x_0,x_1,x_2,x_3,\dots)$ in ${\rm Cham}^\pm\big(A^\pm\big)$ with $x_0 = x$ and $d^*(x_i,y) = n + i$ for $0\leq i$. These two rays determine two ends, $[r_\pm]$ and $[r_\mp]$, and we
call the pair $([r_+],[r_-])$ the {\it twin end} associated with $(x,y)$.
In this twin tree case, it is known \cite[Chapter~11, twin buildings]{Ronan89} that there exists a twin apartment $\big(A^+,A^-\big)$ with $x\in A^+$ and $y\in A^-$ and then we have $r_\pm\subset A^\pm$.

For $c^\pm \in {\rm Cham}_\infty^\pm$ and $d^\mp \in {\rm Cham}_\infty^\mp$ we define the $W_{\infty}$-valued codistance function by
\[\delta_\infty^*(c^\pm, d^\mp) = \begin{cases} -1& \hbox{if }c^\pm\hbox{ and }d^\mp \hbox{ are twin ends},\\
\hphantom{-}1& \hbox{otherwise}. \end{cases}\]

This means that for each chamber $c\in {\rm Cham}_\infty^\pm$ there is exactly one chamber $d\in {\rm Cham}_\infty^\mp$
with codistance $\delta^*(c,d)= -1$, while all other chambers $d'\in {\rm Cham}_\infty^\mp$ have codistance $\delta^*(c,d') = 1$. This is consistent with Lemma~\ref{lemma:longer_element_chamber_unique}.

We denote all objects associated to the spherical building at infinity, i.e., apartments, chambers or the Weyl group, with the subscript~`$\infty$'. Thinking of apartments of $\cB_\infty$ as the ends of apartments of~$\cB$, our goal now is to prove that the embedding of the twin building $\big(\cB^+, \cB^-, \delta^*\big)$, given in Theorem~\ref{theorem:mainembedding theorem} induces an embedding of the spherical building at infinity $\big(\cB^+_\infty, \cB^-_\infty, \delta^*_\infty\big)$ into a set of rays on the null cones of $\ft_0$.

In the rank~2 case, the compact real form of the Cartan subalgebra, $\ft$,
has a bilinear form with signature $(1,1)$. So for each $0\neq r\in\R$, $\ft_r$ is a hyperbola with two connected components (branches) $\ft_r^\pm$, and $\ft_0 = \partial\cL_\ft$ is a pair of lines (the asymptotes). We have $\cL_\ft = \cup_{r\leq 0} \ft_r$ and the real roots of
$\fg$ are on the hyperbolas $\ft_{(\a_i,\a_i)}$ for $i=1,2$. For symmetric Cartan matrices, $A$, all real roots have the same length, so
this is just one hyperbola. All imaginary roots of $\fg$ are in $\cL_\ft^0 = \{x\in\ft\, |\, (x,x)<0\}$, which is the Tits cone. We define an
equivalence relation on the nonzero vectors in $\ft_0$ by saying that two nonzero points, $x$ and $y$, are equivalent when $x = ry$ for
some $0<r\in\R$. There are exactly four equivalence classes of such points, corresponding to the four rays in $\ft_0$, which we denote
by $\big\{x_i^+, x_i^-\, |\, i=1,2\big\}$. This corresponds to the ``lightlike closure'' $B^\infty(\partial\cL_\ft)$ in Definition~\ref{definition:lightlike_closure}.
To be more precise,
$x_i^\pm$ denotes the ray in $\ft_0^\pm$ such that $x_2^- = - x_1^+$, $x_1^- = - x_2^+$ and
\[\big(x_1^+,\a_1\big) < 0,\qquad \big(x_1^+,\a_2\big) > 0,\qquad \big(x_2^+,\a_1\big) > 0,\qquad \big(x_2^+,\a_2\big) < 0.\]
This means that in Fig.~\ref{figure:Figure1}, $x_1^+$ is the ray going up to the right,
and $x_2^+$ is the ray going up to the left, while $x_1^-$ is the ray going down to the right, and $x_2^-$ is the ray going down to the left. We think of $\big\{x_1^\pm, x_2^\pm\big\}$ as the ends of the fundamental apartment $A^\pm$ in $\cB^\pm$, so as the fundamental apartment $A^\pm_\infty$ in $\cB^\pm_\infty$. As a Coxeter complex for $W_\infty$, each apartment consists of just two points which are exchanged by the nontrivial element of $W_\infty$, and the $W_\infty$-valued distance function is obvious.

Based on the labeling of the chambers of the fundamental twin apartment $\big(A^+,A^-\big)$ shown in Fig.~\ref{figure:Figure3}, and the
codistance function $\delta^*$, we see that the codistance function on the fundamental twin apartment $(A^+_\infty,A^-_\infty)$ is
\[\delta^*_\infty\big(x_i^+,x_j^-\big) = \begin{cases} -1& \hbox{if }i = j, \\
\hphantom{-}1& \hbox{if }i \neq j. \end{cases}\]

Define the {\em halos}, positive and negative, of $\fg$ to be the union of all $K$ conjugates of $\big\{x_i^\pm\, |\, i=1,2\big\}$
\[\Xi_\infty^\pm = \big\{k x_i^\pm k^{-1}\, |\, k\in K, i=1,2\big\}\]
and let the twin halo of $\fg$ be $\Xi_\infty = \Xi_\infty^+ \cup \Xi_\infty^-$.
This will be where we embed the spherical twin building at infinity $\cB_\infty$ in the next theorem.

Before we can state and prove the embedding theorem, we must define the distance function~$\delta_\infty^\pm$ on $\Xi_\infty^\pm$ and the codistance function $\delta_\infty^*$ on $\Xi_\infty$.

For any two chambers $k_1 x_i^\pm k_1^{-1}$ and $k_2 x_j^\pm k_2^{-1}$ in $\Xi_\infty^\pm$
define the $W_\infty$-valued distance function
\[\delta^\pm_\infty(k_1 x_i^\pm k_1^{-1}, k_2 x_j^\pm k_2^{-1}) =
\begin{cases} \hphantom{-}1& \hbox{if } i = j \hbox{ and }k_1= k_2,\\
-1& \hbox{otherwise}. \end{cases}
\]
and the codistance function
\[\delta_\infty^*(k_1 x_i^\pm k_1^{-1}, k_2 x_j^\mp k_2^{-1}) = \begin{cases}
-1& \hbox{if}\ i = j \hbox{ and }k_1= k_2,\\
\hphantom{-}1& \hbox{otherwise}.\end{cases}\]

Note that $\big(\Xi_\infty^\pm,\delta^\pm_\infty\big)$ is a Tits building with an apartment system determined by $\delta^\pm_\infty$, and $\big(\Xi_\infty^+, \Xi_\infty^-,\delta_\infty^*\big)$ is a spherical twin building. We choose $x_1^+$ to be the fundamental chamber in~$A_\infty^+$ and~$x_2^-$ to be the opposite fundamental chamber in~$A_\infty^-$ with codistance $\delta_\infty^*\big(x_1^+,x_2^-\big) = 1$, consistent with the behavior of the codistance $\delta^*$ in $\cB$.

\begin{Theorem}\label{theorem:embeddingofsphericalbuilding}
 There is a $K$-equivariant bijective twin building map $\Psi_\infty\colon \cB_{\infty}\to \Xi_\infty$ respecting distance and codistance functions $\delta_\infty^\pm$ and $\delta_\infty^*$, such that the following diagram commutes:
\begin{displaymath}
\begin{xy}
 \xymatrix{
 \cB_{\infty} \ar[rr]^{g\in K} \ar[d]^{\Psi_\infty} 		& & \cB_{\infty} \ar[d]^{\Psi_\infty} \\
 \Xi_\infty \ar[rr]^{{\rm Ad}_g} 	& & \Xi_\infty
 }
\end{xy}
\end{displaymath}
\end{Theorem}

\begin{Remark}
Such a spherical building at infinity exists {\em only} in the highly degenerate case of rank 2. The analogous construction for higher
rank hyperbolic Kac--Moody groups would have to involve a new kind of structure beyond the theory of buildings, perhaps replacing
Coxeter groups with some other class of groups. \end{Remark}

Our strategy for the proof of this result is to follow the proof of Theorem~\ref{theorem:mainembedding theorem}, but since chambers consist only of points, there are some simplifications.
We use results of Ronan--Tits~\cite{RonanTits94} about the stabilizers of twin ends.

Recall that $A^{\pm}$ denotes the fundamental apartment in $\cB^\pm$ so it is a line in that tree with two ends. (See Figs.~\ref{figure:Figure3} and~\ref{figure:Figure4}.)
Let $e_{i,\infty}^\pm$, $i=1,2$, denote the two ends of that line, that is, the equivalence classes of certain rays defined as follows.
As shown in Proposition~\ref{prop:root_group_actions_on_Fund_Apt}, for each $i=1,2$, $m\in\Z$, the
subgroup $U_{(w_1 w_2)^m \a_i}$ fixes the end $e_{i,\infty}^\pm$ if we define
\begin{alignat*}{3}
& e_{1,\infty}^+ = [R_{\rm ray}^+(m)],\qquad && e_{2,\infty}^+ = [L_{\rm ray}^+(m)],& \\
& e_{1,\infty}^- = [L_{\rm ray}^-(m)],\qquad && e_{2,\infty}^- = [R_{\rm ray}^-(m)].&
\end{alignat*}
Furthermore, we have shown that $U_{-(w_1 w_2)^m \a_1}$ fixes $e_{2,\infty}^\pm$ and
$U_{-(w_1 w_2)^m \a_2}$ fixes $e_{1,\infty}^\pm$.
This means that $\omega(e_{1,\infty}^+) = e_{2,\infty}^-$ and $\omega(e_{2,\infty}^+) = e_{1,\infty}^-$.

Let us fix some notation, defining the stabilizer of $e_{i,\infty}^\pm$
\begin{Definition}\label{definition:end_stabilizers}
The stabilizer of an end $e_{i,\infty}^\pm$ is the group
\begin{gather*}
B_{i,\infty}^\pm = \big\{g\in G \,|\, g\cdot e_{i,\infty}^\pm=e_{i,\infty}^\pm \big\} .
\end{gather*}
\end{Definition}

Recall the non-standard partition of the real roots of $\fg$, $\roots^{\rm re} = \roots_1\cup \roots_2$ where
\[\roots_1 = W^{\rm even}\{\a_1,-\a_2\}\qquad\hbox{and}\qquad \roots_2 = W^{\rm even}\{-\a_1,\a_2\} .\]
\begin{Definition}\label{definition:nonstandard_root_groups}
For $i=1,2$, we set $U_i=\langle U_{\alpha}\,|\, \alpha\in \roots_i\rangle$ and $B_i = T W^{\rm even} U_i$.
\end{Definition}
We now have the following result.

\begin{Proposition}\label{prop:gens_fixing_ends}
For $i = 1,2$ we have $B_i \leq B_{i,\infty}^\pm$.
\end{Proposition}
\begin{proof}
For any $m\in\Z$, we have established that
$U_{(w_1 w_2)^m \a_1}$ and $U_{-(w_1 w_2)^m \a_2}$ are in $B_{1,\infty}^\pm$, and that
$U_{(w_1 w_2)^m \a_2}$ and $U_{-(w_1 w_2)^m \a_1}$ are in $B_{2,\infty}^\pm$, so $U_i \leq B_{i,\infty}^\pm$.
It is clear that $T W^{\rm even} \leq B_{i,\infty}^\pm$, so we get $B_i \leq B_{i,\infty}^\pm$.
\end{proof}

We would like to make some remarks about why $(B_1,B_2,N)$ is not a twin $BN$ pair, and therefore why we have no
Bruhat decomposition for the twin building at infinity. There are three requirements to be a $BN$-pair,
T1--T3, in Definition \ref{definition:BN-pair}, and three more requirements to be a twin $BN$-pair, TW1--TW3, in
Definition~\ref{definition:Twin-BN-pair}. The condition T1 requires that
$G=\la B_i,N\ra$, $H=B_i \cap N \lhd N$, $W_\infty =N/H = \la S\ra$ for $i = 1,2$. Since $w_i\in N$ and $w_i B_i w_i^{-1} = B_{3-i}$
and $G=\la B_1,B_2\ra$, we have $G=\la B_i,N\ra$. We also have $H = T W^{\rm even} U_i\cap TW = TW^{\rm even}$
so $N/H = (TW)/(TW^{\rm even}) = W_\infty$.
Condition T3 requires that the nontrivial element in $W_\infty$, represented by either $w_1$ or $w_2$, satisfies
$w_iB_j w_i^{-1}\not\subseteq B_j$, and that is true since $B_{3-j}\not\subseteq B_j$.
However, condition T2 requires that for $i,j,k\in\{1,2\}$ we have $w_iB_k w_j\subseteq B_k w_i w_jB_k \cup B_k w_j B_k$,
so in particular, T2 requires that $B_2 = w_1B_1w_1 \subseteq B_1 \cup B_1 w_1 B_1$.
If this were true it would mean that $U_2 = w_1U_1w_1 \subseteq U_1 \cup U_1 w_1 U_1$
which would imply $U_1 \subseteq U_2 \cup w_1 U_1 U_2$. Since $U_1\cap U_2 = 1$, this would require
$U_1 \subseteq w_1 U_1 U_2$ and then $w_1 U_1 \subseteq U_1 U_2$ so $w_1 \in U_1 U_2$. Certainly we know from
the usual formula that $w_1\in U_1 U_2 U_1$, but there is no expression for $w_1$ in $U_1 U_2$.
The only case of concern in TW2 is when $s = w_1 = w$, so that $\ell(sw) = \ell(1) = 0 < \ell(w_1) = 1$, and the requirement
is that $B_1 w_1 B_1 w_1 B_2 = B_1 B_2$, which is true since $w_1 B_1 w_1 = B_2$.
Concerning TW3, since we know $w_1 B_1 = B_2 w_1$ we get that $B_2 w_1 \cap B_1 = w_1 B_1 \cap B_1$ is the empty set
since $w_1\notin B_1$ so these are distinct left cosets of $B_1$.

We thank Peter Abramenko for his exposition of the following results of Ronan-Tits about twin trees, and how they can be
applied to achieve our goal in this section.

We now discuss real {\it roots} and real {\it root groups} from the point of view of twin buildings. Let $\big(\cB^+,\cB^-,d^*\big)$ be a twin building with codistance $d^*$, and fix a fundamental twin apartment $A_0 = \big(A_0^+,A_0^-\big)$. A {\it root} $\a^\pm$ of $\cB^\pm$ is a half-apartment of $A_0^\pm$, and $\Phi^\pm$ is the set of all roots of $A_0^\pm$. For each choice of a root $\a^\pm$, denote by $-\a^\pm$ the other half-apartment such that $-\a\cup\a = A_0^\pm$ and $-\a\cap\a$ is exactly one panel. The use of superscript $\pm$ to distinguish roots (half-apartments) in the two buildings is therefore not always consistent with the choice of a~factor of $\pm 1$ according to choice of building. From the Lie algebra point of view, there is just one set of roots, $\roots$, but the action of a root group $U_\a$ for $\a\in\roots$ depends
on the building $\cB^\pm$. For the twin tree $\cB = (\cB^+,\cB^-,\delta^*)$ we are studying in this section,
if $\a = \pm\roots_2(k)$, the corresponding half-apartment in $\cB^\pm$ would be $L^\pm_{\rm ray}(k)$ which is fixed by the real root group $U_\a$ according to Proposition~\ref{prop:root_group_actions_on_Fund_Apt}, and $-\a = \pm\roots_1(k+1)$ corresponds to the half-apartment $R^\pm_{\rm ray}(k+1)$ in $\cB^\pm$ which is fixed by the real root group $U_{-\a}$. We have $L^\pm_{\rm ray}(k)\cup R^\pm_{\rm ray}(k+1) = A_0^\pm$, and
the intersection $L^\pm_{\rm ray}(k)\cap R^\pm_{\rm ray}(k+1)$ is the single vertex $v^\pm_{k+\frac{1}{2}} = C^\pm(k)\cap C^\pm(k+1)$.

\begin{figure}[t]\centering
\includegraphics{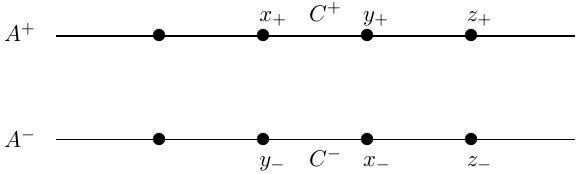}
\caption{Twin apartment $(A^+,A^-)$ in twin tree $(\cB^+,\cB^-)$ with chambers $C^+\op C^-$ and vertices~$x_+\op x_-$ and $y_+\op y_-$.}\label{Figure6}
\end{figure}

As shown in Fig.~\ref{Figure6}, where $C^\pm = (x_\pm,y_\pm)$ and $C^+\op C^-$ with $d^*(x_+,x_-) = 0 = d^*(y_+,y_-)$,
so $d^*(x_+,y_-) = 1 = d^*(x_-,y_+)$. Then there is a vertex $z_-$ adjacent to $x_-$ with $d^*(y_+,z_-) = 2$ and continuing
along a ray in that direction the co-distance function $d^*(y_+,\cdot)$ is increasing. There is also a vertex $z_+$ adjacent to
$y_+$ with $d^*(x_-,z_+) = 2$ and continuing in a ray in that direction the co-distance function $d^*(x_-,\cdot)$ is increasing.
For any pair of opposite chambers, $C^+\op C^-$, there exists a twin apartment $(A^+,A^-)$ with $C^\pm\in A^\pm$.

Let $V^\pm = V^\pm(\cT^\pm)$ be the set of vertices of the tree $\cT^\pm = \cB^\pm$.
Let $x\in V^+$ and $y\in V^-$ with $d^*(x,y) = n > 0$, so we get rays
\begin{gather*} r_+\ \hbox{ of adjacent vertices }(x=x_0,x_1,x_2,x_3,\dots)\ \ \hbox{ with end } e_+ \hbox{ and }\\
r_-\ \hbox{ of adjacent vertices }(y=y_0,y_1,y_2,y_3,\dots)\ \ \hbox{ with end } e_-,
\end{gather*}
where $d^*(x,y_j) = n+j$ for all $j\geq0$ and $d^*(x_i,y) = n+i$ for all $i\geq0$ uniquely determines the vertices $x_i$ and~$y_j$. This gives $d^*(x_i,y_j) = n+i+j$ for all $i,j\geq0$. The pair of ends $(e_+,e_-)$ determined by this process is called a {\it twin end}. Using the $\N$-valued codistance between vertices in twin trees and the $\N$-valued codistance between chambers (edges), it seems clear that for any chambers $C_i^+ = [x_i,x_{i+1}]$ in the ray $[x,e_+)$ and $C_j^- = [y_j,y_{j+1}]$ in the ray $[y,e_-)$, the chamber codistance $\delta^*(C_i^+,C_j^-) > 0$, so no such pair of chambers is opposite.

The following fact was proved by Ronan--Tits \cite[Proposition~3.4]{RonanTits94}.

\begin{Proposition}[\cite{RonanTits94}] \label{prop:Ronan-Tits} Let $x\in V^+$ and $y\in V^-$ with $d^*(x,y) = n > 0$ determine rays~$r_+$ and~$r_-$ with ends $e_+$ and $e_-$. Let $x'\in V^+$ and $y'\in V^-$ be another pair of non-opposite vertices $($so $d^*(x',y')>0)$, with another pair of rays $r'_+$ and $r'_-$ with ends $e'_+$ and $e'_-$. If $e'_+ = e_+$ then $e'_- = e_-$.
\end{Proposition}

\begin{Corollary} \label{corollary:twin-ends} Let $\Gamma$ be any group acting on the twin tree $(\cT^+,\cT^-,d^*)$ preserving
$d^*$. If we have a twin end $(e_+,e_-)$ and if $g\in\Gamma$ with $g\cdot e_+ = e_+$ then $g\cdot e_- = e_-$.
\end{Corollary}
\begin{proof} With notation as above, let $g\cdot (x,y) = (g\cdot x,g\cdot y) = (x',y')$. We can apply $g$ to the rays~$r_+$ and~$r_-$, to get rays $r'_+$ and $r'_-$, which consist of the vertices $g\cdot x_i$ and $g\cdot y_i$ uniquely determined by~$x'$ and $y'$ because $g\in\Gamma$ preserves $d^*$. But then we must have $g\cdot (e_+,e_-) = (g\cdot e_+,g\cdot e_-)$. We are given that $g\cdot e_+ = e_+$, so Proposition~\ref{prop:Ronan-Tits} gives $g\cdot e_- = e_-$.
\end{proof}

\begin{figure}[t]\centering
\includegraphics{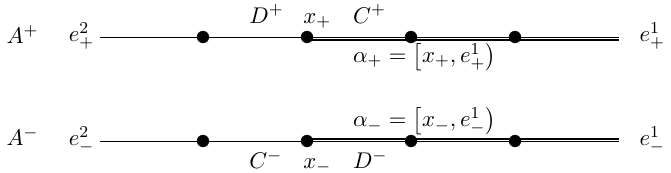}
\caption{Twin root $\a = (\a_+,\a_-)$ in twin apartment $\big(A^+,A^-\big)$ determined by adjacent chambers $C^+$ and $D^+$ with $C^+\cap D^+ = x_+$ and $C^+\in\a_+$, so $\a_- = -{\rm op}\, \a_+$. Also showing twin ends $\big(e^1_+,e^1_-\big)$ of $\a$.}\label{Figure7}
\end{figure}

We have $G = T\la U_\a\,|\,\a\in\Phi\ra$ is generated by its (real) root groups, $U_\a$ and the torus
$T = {\rm Fix}_G\big(A_0^+,A_0^-\big)$ that fixes (chamber-wise) the fundamental twin apartment
$A_0 = \big(A_0^+,A_0^-\big)$ in the twin tree. In the tree, $\cB^\pm$, a root $\a_\pm$ means a ray in the
fundamental apartment $A_0^\pm$, and each root $\a_+$ (ray) in $A_0^+$ determines a root $\a_-$ (ray)
in $A_0^-$, giving a {\it twin root} $\a = (\a_+,\a_-)$, a~pair of rays. That root $\a_-$ is the {\it negative} of the root (ray) in $A_0^-$ whose edges are the opposites of the edges of $\a_+$. This means that $\a_-$ is the maximal ray in $A_0^-$ such that no edge of $\a_-$ is opposite to any edge of $a_+$.
(See Fig.~\ref{Figure7}.)

$U_\a$ acts transitively on the set of chambers (edges) of $\cT^+$ containing $x_+$ different from the chamber $C^+$ in fixed ray $\a_+$.

We have $T = \cap_{\a\in\Phi} N_G(U_\a)$, and for $i = 1,2$, we have defined ends $e^\pm_{i,\infty}$ as equivalence classes of rays, and for the following proof we use the notation $e^i_\pm$ for $e^\pm_{i,\infty}$, so that
\[U_i = \la U_\a\,|\, \a =(\a_+,\a_-),\ \a_+ \hbox{ represents end }e^i_+\ra\]
and
\[T \leq D = TW^{\rm even} = \big\{g\in G\,|\, g\ \hbox{acts by translation on } A_0^+\big\} \leq {\rm Stab}_G\big(A_0^+\big) = N.\]
The following theorem gives a precise description of $B_{i,\infty}^\pm = {\rm Stab}_G\big(e^i_\pm\big)$, but is written only for the case of $i=1$ and $\pm = +$. The result shows that the answer is independent of the choice of~$\pm$.

\begin{figure}[t]\centering
\includegraphics{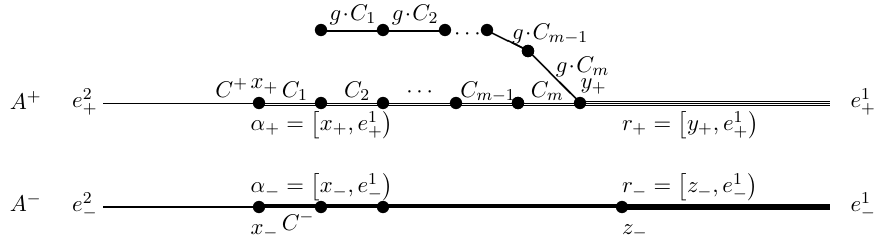}
\caption{Diagram for Case (2) in proof of Theorem~\ref{theorem:Stabilizer} showing: Twin root $\a = (\a_+,\a_-)$ determined by adjacent chambers $C^+$ and $C_1$. Also showing pair of rays $(r_+,r_-)$ contained in $\a$, where $r_+ = \big[y_+,e^1_+\big)$ and $r_- = \big[z_-,e^1_-\big)$. The path of edges $(C_1,C_2,\dots,C_{m-1},C_m)$ forms the interval $[x_+,y_+]$ which contains the sequence of adjacent vertices $(x_+ = x_0,x_1,\dots,x_{m-1},x_m=y_+)$ so $C_i = [x_{i-1},x_i]$ for $1\leq i\leq m$. The path of edges $(g\cdot C_1,g\cdot C_2,\dots,g\cdot C_{m-1},g\cdot C_m)$ forming the interval $[g(x_+), g(y_+)=y_+]$ is also shown.}\label{Figure8}
\end{figure}

\begin{Theorem}\label{theorem:Stabilizer}We have ${\rm Stab}_G\big(e^1_+\big) = D U_1 = B_1$.
\end{Theorem}
\begin{proof} We already have the containment $D U_1 \leq {\rm Stab}_G\big(e^1_+\big)$. Let $g\in G$ with
$g\big(e^1_+\big) = e^1_+$, so for any ray $r_+$ in $A_0^+$ with end~$e^1_+$, the ray $g(r_+)$ also has end $e^1_+$. It means there is a~subray $r'_+\subset r_+$ such that $g(r'_+)\subseteq r_+\subset A_0^+$. If necessary, we may adjust the choice of $r_+$ by adding or deleting one chamber at its beginning, so that there is a translation $d\in D$ with $dg(r'_+) = r'_+$. If we replace~$r_+$ by~$r'_+$ and replace $g$ by $g' = dg$ then $g'$ fixes $r'_+$ pointwise. Without loss of generality, we now assume that $g\in {\rm Stab}_G\big(e^1_+\big)$ and we have a ray $r_+ = \big[y_+,e^1_+\big)$ in $A_0^+$ such that
$g(x) = x$ for every vertex $x\in r_+$. That means the end $e^1_+$ of $r_+$ is fixed by $g$ so $g\big(e^1_-\big) = e^1_-$ by Corollary~\ref{corollary:twin-ends} since $g$ preserves the twin structure. This now implies the existence of a ray $r_- = [z_-,e^1_-)$ in $A_0^-$ with end $e^1_-$ which is pointwise fixed by $g$. If $r_-$ were not pointwise fixed by $g$,
then $z_- \neq g(z_-) \in r_-$ so $\delta^*(y_+, g(z_-)) > \delta^*(y_+, z_-)$, but $\delta^*(y_+, z_-) = \delta^*(g(y_+), g(z_-)) =
\delta^*(y_+, g(z_-))$ since $g(y_+) = y_+$.

Case (1): The pair $(r_+,r_-)$ is not contained in any twin root $(\a_+,\a_-)$ of $A_0 = \big(A_0^+,A_0^-\big)$ so there are chambers $C^\pm\in r_\pm$ with $C^+\op C^-$. Then we have $g\big(C^+\big) = C^+$ and $g\big(C^-\big) = C^-$ which implies that $g$ fixes the unique twin apartment $\big(A_0^+,A_0^-\big)$ containing both $C^+$ and $C^-$, and that implies $g\in {\rm Fix}_G\big(A_0^+,A_0^-\big) = T$.

Case (2): The pair $(r_+,r_-)$ is contained in some twin root $\a = (\a_+,\a_-)$ of $A_0 = \big(A_0^+,A_0^-\big)$. Let $x_\pm$ be the starting point of ray (root) $\a_\pm$ where $x_+\op x_-$ (that is, $d^*(x_+,x_-) = 0$). We have $r_+\subseteq \a_+$ and $r_-\subseteq \a_-$. Let $y_+$ be the starting point of ray $r_+$ and suppose $y_-\op y_+$. Then we have the twin root $\beta = (\beta_+,\beta_-)$ where $\beta_+ = r_+$ and $\beta_-$ has starting point $y_-$ so $r_-\subseteq\beta_-$. Here we have $C^-$ is the first chamber of $\a_-$ with starting vertex $x_-$, and $C^+$ is the opposite chamber which is adjacent to the first chamber of $\a_+$ with starting vertex $x_+$. (See Fig.~\ref{Figure8}.)

Our goal is to modify $g$ to $g'$ by multiplying with elements from $U_1$ such that $g'\cdot C^- = C^-$ and $g'\cdot C^+ = C^+$, which would imply that $g'\in T$. Let $(C_1, C_2,\dots,C_{m-1},C_m)$ be a sequence of
adjacent edges in $A_0^+$ with $C_i = [x_{i-1},x_i]$ for $1\leq i\leq m$ and $x_0 = x_+$, forming the interval $[x_+,y_+]$ so that $C_m\cap r_+ = y_+ = x_m$. Then $(g\cdot C_1,g\cdot C_2,\dots,g\cdot C_{m-1}, g\cdot C_m)$ is a sequence of adjacent edges in $\cT^+$ and $g\cdot C_m$ is attached to vertex $y_+=x_m$ since $g$ fixes $r_+$ pointwise. For $0\leq i\leq m$ define the ray (root) $\beta_i = \big[x_i,e_+^1\big)$ so that $U_{\beta_i}\leq U_1$ fixes $\beta_i$ pointwise and acts transitively on the set of all edges attached to $x_i$. There is an element $u_m\in U_{\beta_m}$ such that $(u_m g)\cdot C_m = u_m\cdot(g\cdot C_m) = C_m$ so $u_m g$ fixes the ray $\beta_{m-1} = C_m\cup r_+$ pointwise.
Then there is an element $u_{m-1}\in U_{\beta_{m-1}}$ such that $(u_{m-1} u_m g)\cdot C_{m-1} = C_{m-1}$ so
$u_{m-1} u_m g$ fixes the ray $\beta_{m-2} = C_{m-1}\cup C_m\cup r_+$ pointwise. Continuing this way, we find a~sequence of elements $u_i\in U_{\beta_i}$ for $1\leq i\leq m$ such that $g' = u_1 u_2\cdots u_m g$ fixes the ray $\a_+ = C_1\cup C_2\cup\cdots \cup C_{m-1}\cup C_m\cup r_+$ pointwise. Since $G$ acts on the twin tree preserving $d^*$, it preserves the opposition relation, so $\a_- = -\op \a_+$ means $g'$ also fixes the ray $\a_-$ pointwise. Thus, $g'$ fixes the twin end $(e^1_+,e^1_-)$ and by Corollary~\ref{corollary:twin-ends}, $g'$ fixes the other twin end $\big(e^2_+,e^2_-\big)$ of the twin apartment $A_0 = \big(A_0^+,A_0^-\big)$. That implies that $g'$ stabilizes $A_0$, so $g'\in T\leq D$ and then $g \in D U_1$.
\end{proof}

Theorem~\ref{theorem:Stabilizer} gives a complete characterization of the stabilizers of an end, yielding
\[B_{i,\infty}^\pm = T W^{\rm even} U_i = B_i q\qquad \textrm{for}\quad i =1,2\]
is independent of the choice of $\pm$.

Let $B_{i,K} = B_i \cap K$. In our next step towards the proof of Theorem~\ref{theorem:embeddingofsphericalbuilding},
we show that $B_{1,K} = B_{2,K} = B_K$ is independent of the choice of $i$, and establish a $K$-equivariant bijection between
$G/B_i$ and $K/B_K$.

\begin{Lemma}\label{lemma:U-orbit}
For $i=1,2$ we have $U_i \cdot \fh = \fh\oplus \bigoplus_{\a\in\roots_i} \fg_\a$.
\end{Lemma}

\begin{proof}With our assumptions on the rank $2$ Cartan matrix, $U_i$ is abelian since $[e_{\a_j}, e_{\a_k}] = 0$ for any two roots
$\a_j,\a_k\in\roots_i$ and for any root vectors $e_{\a_j}\in\fg_{\a_j}$ and $e_{\a_k}\in\fg_{\a_k}$. A general element in $U_i$ can be
written as $u = \prod\limits_{j=1}^r \exp({\rm ad}_{e_{\a_j}})$, so for any $h\in\fh$ we have
\begin{gather*}
u(h) = h - \sum_{j=1}^r \a_j(h) e_{\a_j} \in \fh \oplus \bigoplus_{j=1}^r \fg_{\a_j}.\tag*{\qed}
\end{gather*}\renewcommand{\qed}{}
\end{proof}

\begin{Lemma}{\label{Bcomponent}}We have $B_{i,K} =T W^{\rm even} = B_K$ is independent of $i=1,2$.
\end{Lemma}

\begin{proof}By Theorem~\ref{theorem:Stabilizer} we have $B_i = T W^{\rm even} U_i$. It is clear that $T W^{\rm even} \subset K$. For $i=1,2$, set $(U_i)_K = U_i\cap K=\{u\in U_i \,|\, \omega(u)=u\}$ where $\omega$ denotes the Cartan involution. Since $\omega(\fg_\a) = \fg_{-\a}$ for any $\a\in\roots$, we get $\omega\circ U_\a \circ \omega ^{-1} = U_{-\a}$ from Lemma~\ref{lemma:auto_action_on_real_root_groups}.
Suppose $1\neq u = \prod\limits_{j=1}^r \exp({\rm ad}_{e_{\a_j}})\in(U_i)_K$ for $\a_j\in\roots_i$, $1\leq j\leq r$, where $e_{\a_j}\in\fg_{\a_j}$. Then $u = \omega(u) = \prod\limits_{j=1}^r \exp({\rm ad}_{e_{-\a_j}})$ for some basis vector $e_{-\a_j}\in\fg_{-\a_j}$.
Let $h\in\fh$ be an element such that $\a_j(h)\neq 0$ for $1\leq j\leq r$. From Lemma~\ref{lemma:U-orbit}, we get $u(h) = h - \sum\limits_{j=1}^r \a_j(h) e_{\a_j}$. Using the formula above for $\omega(u)$ we similarly calculate $\omega(u)(h) = h + \sum\limits_{j=1}^r \a_j(h) e_{-\a_j}$, but then $u(h) = \omega(u)(h)$ gives the non-trivial dependence relation $- \sum\limits_{j=1}^r \a_j(h) e_{\a_j} = \sum\limits_{j=1}^r \a_j(h) e_{-\a_j}$ among the independent basis vectors from the root spaces
$\oplus_{j=1}^r \fg_{\a_j}$ and $\oplus_{j=1}^r \fg_{-\a_j}$. This contradiction proves that $u = 1$ so $(U_i)_K=\{1\}$.
\end{proof}

\begin{Corollary}\label{lemma:K-orbit} The maps ${\rm Ch}^\pm\colon K/B_K \to {\rm Cham}\big(\cB_{\infty}^\pm\big)$ given by ${\rm Ch}^\pm(kB_K) = k\big(e_{i,\infty}^\pm\big)$ are bijective for either choice of $i=1,2$ and either choice of~$\pm$.
\end{Corollary}

\begin{proof}For each choice of $\pm$, $K$ acts transitively on ${\rm Cham}\big(\cB_{\infty}^\pm\big)$, which is the $K$ orbit of $e_{i,\infty}^\pm$, with stabilizer $B_i\cap K = B_K$, so the orbit-stabilizer theorem gives the result.
\end{proof}

\begin{Lemma}\label{lemma:stabilizer} For $\Xi_\infty = \Xi_\infty^+ \cup \Xi_\infty^-$ the twin halo of $\fg$, each element the set $\big\{x^\pm_i\, |\, 1=1,2\big\}\subset\Xi_\infty^\pm$ is fixed by $B_K$, which is the pointwise stabilizer in $K$ of that set.
\end{Lemma}

\begin{proof}The action of the even Weyl group $W^{\rm even}$ fixes the four rays of $\partial\cL_\ft$, so its induced action on $\big\{x^\pm_i\, |\, 1=1,2\big\}$ is trivial. The same is true of the adjoint action of $T=\exp(\ft)$ which fixes $\ft$ pointwise, so this is true for $B_K$, which is contained in the pointwise stabilizer of this set. Now suppose that $k\in K$ stabilizes the set pointwise. Then $k$ sends two linearly independent points on the two lines of $\partial\cL_\ft$ to linearly independent points, a basis for $\ft$, so $k$ normalizes $\ft$. But $N_K(\ft) = TW$, and the elements in the set $TW^{\rm odd}$ exchange $x^\pm_1$ and $x^\pm_2$, so $k\in TW^{\rm even} = B_K$.
\end{proof}

We now have all technical ingredients to prove Theorem~\ref{theorem:embeddingofsphericalbuilding}.

\begin{proof}[Proof of Theorem~\ref{theorem:embeddingofsphericalbuilding}]
We begin to define $\Psi_\infty$ by setting
\[\Psi_\infty\big(e_{1,\infty}^+\big) = x_1^+ \qquad\hbox{and}\qquad
\Psi_\infty\big(e_{1,\infty}^-\big) = x_1^-\]
and we extend this to all of $\cB_\infty$ by
\[\Psi_\infty\big((kB_K)e_{1,\infty}^\pm\big) = k \Psi_\infty\big(e_{1,\infty}^\pm\big) k^{-1}\]
for any $k\in K$. This map is well-defined and injective by Lemma~\ref{lemma:stabilizer}.
It is bijective by Corollary~\ref{lemma:K-orbit}, and is clearly $K$-equivariant by its definition.
\end{proof}

The results in this section suggest that on the Lie algebra level one could study a {\em non-standard Cartan decomposition}
\[\fg = \fh \oplus \bigoplus_{\a\in P_1} \fg_\a \bigoplus_{\a\in P_2} \fg_\a \]
based on a non-standard partition $\Phi = P_1 \cup P_2$ such that $\Phi_i \subset P_i$ for $i=1,2$, and such that
\[\fn_1 = \bigoplus_{\a\in P_1} \fg_\a \qquad \hbox{and}\qquad \fn_2 = \bigoplus_{\a\in P_2} \fg_\a\]
are subalgebras. Then $\fb_i = \fh\oplus\fn_i$, $i=1,2$, would be non-standard Borel subalgebras corresponding to the subgroups
$B_i$. That would be accomplished if each subset $P_i$ were closed under addition.
Since $\Phi$ is a subset of the dual of the split real Cartan subalgebra, $\fh_\R^*$, one might consider using a timelike line in the interior of the lightcone $\cL_{\fh_\R}^0$ to partition all of $\Phi$. If such a~line does not contain any roots, every root is on one side or the other, and~$\Phi_1$ and~$\Phi_2$ are on opposite sides for any choice of the line. But if the line contains
roots, one would have to decide which ones go in which part of the partition, and it must be done in such a way that each~$P_i$ is closed under addition. The solution would be to divide up the line into two rays from the origin, and divide up the roots on the line according to which ray they are in.

There are two obvious partitions determined by the two lightcone lines themselves. For the line determined by $x^+_2$ the partition would be $P_1 = \Phi_1 \cup \big(\Phi^{\rm im}\big)^+$ and $P_2 = \Phi_2 \cup \big(\Phi^{\rm im}\big)^-$, while for the line determined by $x^+_1$ the partition would be $P_1 = \Phi_1 \cup \big(\Phi^{\rm im}\big)^-$ and $P_2 = \Phi_2 \cup \big(\Phi^{\rm im}\big)^+$. We have found only two distinct non-standard Borels at infinity, $B_i$, $i = 1,2$ which may correspond to these two partitions, but we have not yet seen a family of non-standard Borels corresponding to
other choices of partitions.

Having found a non-standard Cartan decomposition as above, one gets a corresponding non-standard decomposition of the universal enveloping algebra $\cU(\fg) = \cU(\fn_1)\cU(\fh)\cU(\fn_2)$ and can construct induced Verma modules of two types,
${\rm Verma}^i(\lambda) = \cU(\fn_i)v^i_\lambda$, $i=1,2$, where $v^i_\lambda$ are vectors such that
$h\cdot v_\lambda^i = \lambda(h) v_\lambda^i$ for any $h\in\fh$, $\fn_1\cdot v^2_\lambda = 0$ and $\fn_2\cdot v^1_\lambda = 0$.
The quotient of such a Verma module by its maximal proper submodule would be an irreducible module, ${\rm Irred}^i(\lambda)$ generated by $v^i_\lambda$. These would be integrable $\fg$-modules for $\lambda = n_1\lambda_1 + n_2\lambda_2$ in the weight lattice of $\fg$ ($\lambda_1$, $\lambda_2$ are the fundamental weights of $\fg$ such that $\lambda_i(h_j) = \delta_{ij}$) but outside of the lightcone, so that~$n_1$ and~$n_2$ have opposite signs. Examples of such integrable modules, in addition to the adjoint representation, were mentioned by
Borcherds in~\cite[Section~6]{Borcherds86}.
Such modules have been found to occur in the decomposition of the rank 3 hyperbolic algebra,
$\cF$, with respect to its rank~2 hyperbolic Fibonacci subalgebra, recently studied by Penta~\cite{Penta16}.

\section{Conclusion and further directions}

Our lightcone embedding of the twin building of a hyperbolic Kac--Moody group is motivated by the conjectural existence of hyperbolic Kac--Moody symmetric spaces. There have been some efforts recently to develop the geometry of hyperbolic Kac--Moody symmetric spaces, building on work of the third author on the construction of affine Kac--Moody symmetric spaces~\cite{Freyn09, Freyn10b, Freyn10d,Freyn14, Freyn15b, FreynHartnickHornKoehl17}.

Recalling the well-known finite-dimensional theory (see for example~\cite{Eberlein96, Helgason01}), the boundary of a symmetric space $M$ of non--compact type $IV$ corresponding to a complex simple Lie group~$G$, can be identified with the building over $\mathbb{C}$ associated to~$G$. Via the isotropy representation at a point $p\in M$, the building can be embedded into the unit sphere of the tangent space~$T_pM$. In this way, points in the building get identified with directions in the tangent space of the symmetric space. Via the duality between the compact type and non-compact types, we can identify spaces of type~$IV$ with spaces of type~$II$ and thus also obtain an embedding of the building into the tangent space of a compact symmetric space.

In the absence of hyperbolic Kac--Moody symmetric spaces, important properties of the local geometry are captured via the embedding of the building into the Lie algebra. Our embedding of the building into the Lie algebra gives local pictures of the tangent spaces of conjectural hyperbolic Kac--Moody symmetric spaces of types~$II$ and~$IV$. We note however, that since our twin building embeds into the lightcone of the compact form of the Lie algebra, it captures only the timelike directions in the tangent space.

Via Proposition~\ref{proposition:orbits orthogonal}, we have also obtained a hyperbolic analog of the notion of a {\em polar representation} by the group~$K$ on the $K$-conjugacy class of a Cartan subalgebra. Recall that a~group representation $G\colon V\longrightarrow V$ on a vector space $V$ is called {\em polar} if there exists a~subspace $\Sigma\subset V$, called a {\em section}, such that each orbit $G\cdot v$ for $v\in V$ intersects~$\Sigma$ orthogonally. Finite-dimensional polar representations are orbit equivalent to isotropy representations of finite-dimensonal Riemannian symmetric spaces~\cite{Dadok85}. Similar observations were made in the affine case~\cite{Gross00, HPTT}. Our Proposition~\ref{proposition:orbits orthogonal} shows that the action of~$K$ on~$\fH$ is a polar action with section~$\ft$.

Further questions about polar representations for hyperbolic Kac--Moody groups remain open and the full differential geometry need to develop hyperbolic Kac--Moody symmetric spaces remains elusive. We hope to take this up elsewhere.

\appendix
\section[A new formula for the Weyl group generators on any integrable module]{A new formula for the Weyl group generators\\ on any integrable module}

It is remarkable that the proof of the following theorem, starting from a simple calculation in the two-dimensional representation of
$\fsl(2,\C)_i$, generalizes to a result in any integrable module for any Kac--Moody algebra $\fg$. Special cases of this formula appeared
in some physics papers, for example, in \cite{DamourHillmann09}, where Damour and Hillmann found it in a representation of
$K(E_{10})$.

It is known \cite{Kac90} that for any integrable representation $\phi\colon \fg\to {\rm End}(V)$, the group $\tW$ generated by the operators
\[\tw_i^\phi = \exp(\phi(e_i)) \exp(\phi(-f_i) \exp(\phi(e_i)),\qquad 1\leq i\leq n\]
is a subgroup of the Kac--Moody group $G$ that acts on $V$ as follows. If $V_\mu$ is the $\mu$-weight space of $V$, then $\tw_i^\phi (V_\mu) = V_{w_i(\mu)}$, where $w_i$ is the reflection with respect to the simple root $\alpha_i$ in the dual space, $\fh^*$, of the standard Cartan subalgebra, $\fh$. The Weyl group, $W$, is the Coxeter group generated by those simple reflections, and is also defined on $\fh$ by the formula $w_i(h) = h - \alpha_i(h) h_i$. There is a surjection from $\tW$ onto $W$ whose kernel is given by Remark 3.8 in \cite{Kac90}. In the case when $\phi$ is the adjoint representation, the restriction of each $\tw_i^\phi$ to $\fh$ equals $w_i$.

For any integrable representation $\phi\colon \fg\to {\rm End}(V)$, define the ``compact" operators from $G^\theta$,
\[s_i^\phi = \exp(\phi(\pi x_i)) = \exp(\phi(\pi(e_i - f_i)/2)),\qquad 1\leq i\leq n.\]
We prove that $\tw_i^\phi = s_i^\phi$. This provides a new description of $\tW$ as generated by operators from the real compact form. This generalizes a result discovered from a physics point of view by Damour and Hillmann in a representation of~$K(E_{10})$.

\begin{Theorem}\label{theorem:newWeylFormula} Let $\fg$ be any Kac--Moody algebra of rank $n$ with the usual Chevalley generators $e_i$, $f_i$, $h_i$ and let $\phi\colon \fg\to {\rm End}(V)$ be any integrable representation. Then for $1\leq i\leq n$, we have
\[\exp(\phi(e_i)) \exp(\phi(-f_i) \exp(\phi(e_i)) = \exp(\phi(\pi(e_i - f_i)/2)).\]
\end{Theorem}

\begin{proof}Using the notations as above, for $1\leq i\leq n$, let $\fg_i = \fsl(2,\C)_i$ be the Lie subalgebra of~$\fg$ with basis~$\{e_i,f_i,h_i\}$. Then $V$ has a direct sum decomposition into irreducible $\fg_i$-modules
\[V = \bigoplus_{j\in J} V_j(m),\]
where $\dim(V_j(m)) = m+1$ and the index set $J$ includes information about the $\fg_i$ highest weight vector in $V_j(m)$ that locates it in~$V$.
It is clear that $\tw_i^\phi = s_i^\phi$ on $V$ if and only if their restrictions to each $V_j(m)$ are equal. They are certainly equal for all trivial one-dimensional modules, $V_j(0)$. On any two-dimensional module, $V_j(1)$, the computation comes down to the simple fact that
\[\exp\left( \bm 0&1\\0&0\ebm \right) \exp\left( \bm 0&0\\-1&0\ebm \right) \exp\left( \bm 0&1\\0&0\ebm \right)
= \bm 1&1\\0&1\ebm \bm 1&0\\-1&1\ebm \bm 1&1\\0&1\ebm = \bm 0&1\\-1&0\ebm \]
matches
\[\exp\left( \bm 0&t\\-t&0\ebm \right) = \bm \cos(t)&\sin(t) \\-\sin(t)&\cos(t) \ebm,\]
when $t = \pi/2$. The well-known tensor product decomposition of irreducible $\fsl(2,\C)$-modules says that for any integer $m\geq 1$, we have
\[V(m)\otimes V(1) = V(m+1) \oplus V(m-1).\]
So if the operators are equal on modules $V(m)$ and $V(1)$ then they are equal on the tensor product, so they are equal on the component
$V(m+1)$. This proves by induction that they are equal on any $V_j(m)$ that occurs in $V$.
\end{proof}

\subsection*{Acknowledgements}
This material is based upon work supported by the National Science Foundation under Grant No.~1002477. The first author was supported in part by the Simons Foundation, Mathematics and Physical Sciences-Collaboration Grants for Mathematicians, Award Number: 422182. All authors wish to thank the IH\'ES for support during various visits during 2013--2019. The second and third authors wish to thank the Max-Planck Institute for Gravitational Physics (Albert Einstein Institute), Potsdam, Germany, for support during various visits during 2013--2019.

The authors wish to thank Peter Abramenko for his helpful comments on an earlier draft of the manuscript, and for more recent comments on twin tree structures. They would also like to thank Victor Kac for helpful comments in May 2015 at IH\'ES. The second author would like to thank Kai-Uwe Bux, Max Horn, Tobias Hartnick, Ralf K\"ohl and Peter Abramenko for helpful discussions at the June 2015 conference on ``Generalizations of Symmetric Spaces'' in Israel. Finally, the authors wish to express their thanks to the anonymous referees for many valuable suggestions to improve this paper.

\pdfbookmark[1]{References}{ref}
\LastPageEnding

\end{document}